\documentclass[aos]{imsart}

\RequirePackage{amsthm,amsmath,amsfonts,amssymb}
\RequirePackage[numbers,sort&compress]{natbib}
\RequirePackage[colorlinks,citecolor=blue,urlcolor=blue]{hyperref}
\RequirePackage{graphicx}
\RequirePackage{subcaption} 
\usepackage{enumitem}
\usepackage[ruled,linesnumbered]{algorithm2e} 
\usepackage{booktabs,multirow}
\usepackage{float}
\usepackage{threeparttable}
\usepackage{tikz}
\usetikzlibrary{calc}
\usepackage{fix-cm}
\usepackage{comment}

\startlocaldefs
\theoremstyle{plain}

\newtheorem{theorem}{Theorem}

\newtheorem{proposition}{Proposition}
\newtheorem{corollary}{Corollary}
\theoremstyle{definition}

\newtheorem{condition}{Condition}
\newtheorem{remark}{Remark}
\newcommand{\EXPT}{\mathbb{E}}

\newcommand{\mytrans}{\top}
\newcommand{\F}{\mathrm{F}}
\newcommand{\op}{\operatorname{op}}

\DeclareMathOperator*{\argmin}{arg\,min}

\endlocaldefs

\begin{document}

\begin{frontmatter}
\title{Bridging Theory and Practice:\\Statistical Inference for  Latent Space Models of  Networks}
\runtitle{} 

\begin{aug}
%%%%%%%%%%%%%%%%%%%%%%%%%%%%%%%%%%%%%%%%%%%%%%%
%% Only one address is permitted per author. %%
%% Only division, organization and e-mail is %%
%% included in the address.                  %%
%% Additional information such as            %%
%% identifying the corresponding author must %%
%% be included in in the Acknowledgments     %%
%% section if necessary.                     %%
\author[A]{\fnms{Yuang}~\snm{Tian}\ead[label=e1]{yatian@ust.hk}}
\author[B]{\fnms{Jiajin}~\snm{Sun}\ead[label=e2]{jsun5@fsu.edu}}
\author[C]{\fnms{Yinqiu}~\snm{He}\ead[label=e3]{yinqiu.he@wisc.edu}}
%%%%%%%%%%%%%%%%%%%%%%%%%%%%%%%%%%%%%%%%%%%%%%
%% Addresses                                %%
%%%%%%%%%%%%%%%%%%%%%%%%%%%%%%%%%%%%%%%%%%%%%%
\address[A]{Department of Mathematics, Hong Kong University of Science and Technology \printead[presep={ ,\ }]{e1}}

\address[B]{Department of Statistics, Florida State University \printead[presep={,\ }]{e2}}

\address[C]{Department of Statistics, University of Wisconsin-Madison \printead[presep={,\ }]{e3}}
\end{aug}

\begin{abstract}
Latent space models have been widely adopted in modeling network data. 
Developing statistical inference for  estimated model parameters enables quantifying associated uncertainty and is pivotal for downstream tasks. Despite recent progress on statistical inference of maximum likelihood estimation,   crucial gaps remain  between asymptotic theoretical guarantees and practical use. 
Specifically, how are the oracle maximum likelihood estimators related to the solutions produced by algorithms in practice? 
Can rigorous guarantees be established for existing algorithms without unnecessary restrictions? 
% What is the relationship  between theoretical maximum likelihood estimators and solutions from  algorithms in practice? Can theoretical guarantees be obtained without unnecessary restrictions in existing algorithms? 
To address these fundamental questions, 
we  develop a unified analytical framework that bridges theory and practice of statistical inference for latent space models. 
First, for the maximum likelihood estimation, we relax  the spectral-multiplicity constraint  in the existing asymptotic theory to broaden the applicability. 
Second, we overcome the dependence on unknown true parameters in  prior  algorithmic analyses by developing   novel adaptive criteria and theoretical tools. 
For the widely used algorithm based on the projected gradient descent and  the  singular value thresholding,  we explicitly  connect their outputs  to the maximum likelihood estimator without relying on unknown information. Our results provide a solid foundation for practically useful and statistically principled statistical inference in network analysis.  
\end{abstract}



\end{frontmatter}
%%%%%%%%%%%%%%%%%%%%%%%%%%%%%%%%%%%%%%%%%%%%%%
%% Please use \tableofcontents for articles %%
%% with 50 pages and more                   %%
%%%%%%%%%%%%%%%%%%%%%%%%%%%%%%%%%%%%%%%%%%%%%%
%\tableofcontents

\section{Introduction}

Networks encode relational information between entities in  complex systems and have become ubiquitous across  various scientific domains \citep{sengupta2025statistical}, including social science \citep{park2023social,lu2025zero}, economics \citep{brauning2020dynamic},  neuroscience \citep{wang2023locus,wang2024sex}, and biomedical studies  \citep{nakis2025signed}. In most settings,  a network can be represented as a graph consisting of a set of nodes and a set of edges, which encode the entities in the system and the interactions  between them, respectively.

 To capture complex edge dependence  and other structural characteristics in networks, latent space modeling is one of the widely used approaches    \citep{hoff2002latent,matias2014modeling,athreya2018survey,smith2019geometry}. The central  idea  is to associate each node with a  low-dimensional vector in the latent space. The  relational structure between two nodes, such as their connecting probability or weighted edge value, can then be characterized by a kernel function of  latent vectors.  One common  choice of the kernel function is the inner product between latent vectors, which can capture important network  features, such as transitivity, reciprocity, and  community structures \citep{zhao2012consistency,sussman2013consistent,athreya2018survey,athreya2021estimation,ma2020universal,li2023statistical,li2025high,macdonald2025latent,fang2025transfer}. 

In this work, we study statistical inference of  the latent vectors. 
Estimated latent vectors   are often used for interpretation and visualization  to investigate fundamental structures underlying networks.  Quantifying their  uncertainty   is therefore essential for drawing uncertainty-aware  conclusions in downstream tasks that rely on these estimates, including vertex clustering \citep{rohe2011spectral}, network-based regression \citep{hayes2022estimating,nath2022identifying},  testing for networks  \citep{tang2017nonparametric}, and link prediction \citep{pan2025latent}. Nevertheless, inference for  latent vectors poses unique challenges that the number of unknown vectors grows with respect to the network size, and these vectors themselves are inherently unidentifiable. 
% This makes the nalawhich differs from classical parametric estimation problems.  
% Therefore, this question has  drawn increasing attention.  

% Under the random dot product graphs (RDPG) \citep{young2007random}, 
In the existing literature, one research line focuses on  the random dot product graphs (RDPG) \citep{young2007random} and studies the  asymptotic distributions of spectral  estimators  \citep{athreya2016limit,tang2017semiparametric,tang2018limit}.  Meanwhile, \cite{xie2023efficient} show  that spectral embeddings do not fully exploit  the Bernoulli sampling likelihood information and  could therefore be suboptimal in terms of asymptotic covariance matrices. To address this, \cite{xie2023efficient} propose  an efficient one-step update based on solving likelihood-based  estimating equations and demonstrate  improved efficiency under the RDPG model. 

Another line of research considers    
a broad class of latent space models with general link functions. 
% and inner-product kernels,  
In this setting, 
\cite{li2023statistical} establish  uniform consistency and asymptotic distributions for the maximum likelihood estimators of the latent vectors. Technically, to handle the non-identifiability of latent vectors, they introduce a Lagrange-type  penalty function. %singularity of the  Hessian matrix. 
Computationally, 
% To solve the maximum likelihood estimators, 
\cite{li2023statistical} solve the maximum likelihood estimators by the two-stage approach proposed in \cite{ma2020universal}, which utilizes universal singular value thresholding \citep{chatterjee2015matrix} for initialization, followed by the projected gradient descent for  optimization. 
Despite these advances, an important gap remains between the available inferential theory and the algorithms used in practice. 
% However, for conducting statistical inference in practice,  
%  critical gaps exist between existing theoretical and algorithmic frameworks.  
% it remains   an open question  
In particular, it is unclear whether the estimator produced by the algorithm in \cite{ma2020universal} coincides with the constrained maximum likelihood estimator defined in \cite{li2023statistical}, putting a question on applying the inferential results in practice. 

This  gap  cannot be directly addressed by existing developments and poses unique challenges. 
First, the analysis in \cite{li2023statistical}  assumes that the matrix of latent vectors has  non-vanishing  eigengaps among all non-zero  eigenvalues, and that  the likelihood is maximized within a bounded region containing the true parameters.  However, these assumptions can be difficult to verify or 
implement because the ground truth is unknown in practice. 
% \citep{cheng2021tackling,agterberg2023distributional}.   
Second,  \cite{li2023statistical} examine a Lagrange-type penalized log-likelihood function, 
whereas the gradient descent algorithm in \cite{ma2020universal} directly operates on the unpenalized log-likelihood. 
% Empirically,  
% % observations show that 
% the algorithm in    \cite{ma2020universal} performs well without the additional penalties or restrictions. The existing theoretical analyses do not directly explain such observations.
This mismatch is especially important because the original likelihood for latent vectors is nonconvex. Consequently, unlike in convex problems,  there is no default link between the ideal maximum likelihood estimator and the output of the optimization algorithm. 
% so classical results for parametric models with convex log-likelihoods are not directly applicable. 
% Notably, the original likelihood function of latent vectors is non-convex, so  classical results under  parametric models with convex log-likelihood do not apply. 
Third, the original algorithmic formulation  in \cite{ma2020universal} includes  impractical  or seemingly unnecessary steps to reach theoretical guarantee of convergence. %in theory.  
% \textit{Practical puzzle in the algorithm: } 
For example, \cite{ma2020universal}  point out    that the step of projecting parameters onto unknown bounded sets is unnecessary in practice, even though it is required in its theoretical analysis.  In addition, 
convergence guarantees require tuning parameters, such as the projection set used in singular value thresholding and the step size in gradient descent, to be chosen appropriately. 
% to guarantee  the convergence of the algorithm,   appropriate tuning parameters, such as the projecting set in the singular value thresholding and the step size in gradient descent, need to be chosen. 
In general, however, suitable choices of these quantities depend on unknown model parameters and are therefore infeasible in practice.

% But a good choice for convergence  depends on unknown true model parameters  and is therefore infeasible in general. 

% \noindent \textbf{Our contributions:} 
In this work, we examine the   maximum likelihood estimation of latent vectors under the  class of latent space models considered  in \cite{ma2020universal} and  \cite{li2023statistical}. 
Our goal is to build a rigorous bridge between the practical computation of the constrained maximum likelihood estimator (MLE) and the inferential theory developed for its idealized formulation. 
% We make novel contributions by building foundational connections   between the practical computation of the constrained maximum likelihood estimator (MLE) and its   inferential theory under the idealized formulation. 
Specifically, our main contributions are summarized as follows. 
\begin{itemize}
    \item[1.]\textit{MLE inferential theory.}    We develop a unified  inferential framework for the constrained MLE of latent vectors. A key  feature of our theory is that it  surpasses the spectral restrictions   commonly imposed in the existing literature. 
    % on the spectrum of the matrix of latent vectors  in the literature. 
    \item[2.]\textit{Data-adaptive algorithms.} For the widely-used two-stage algorithm, we develop novel adaptive adjustments that eliminate their dependence on unknown true model parameters. In particular, we develop new adaptive line search conditions for selecting appropriate step sizes and rigorously justify skipping the  projection onto unknown constraint sets in the gradient descent algorithm. Moreover, we develop a range-adaptive singular value thresholding for the initialization that can be applied to both bounded and unbounded settings. 
    \item[3.]\textit{Bridging MLE and algorithms.} We build explicit connections between the constrained  MLE and the outputs of  practical algorithms.  Notably, we show that although the definition of the constrained MLE relies on the true parameters, the practical algorithm can approximate this estimator    under suitable conditions without access to   those true parameters. Our analysis directly examines   solutions of the non-penalized likelihood maximization and does not rely on  additional transformations. 
\end{itemize}

The rest of this paper is organized as follows.  Section \ref{sec:model_homo} introduces the   class of latent space models under consideration and the problem setup. Section \ref{sec:gap_homo} explains the  gaps between existing theoretical and algorithmic frameworks in detail. Section \ref{sec:mle_theory_1} presents our new asymptotic theory for the maximum likelihood estimators, while  guarantees on practical algorithms are given  in Section \ref{sec:gd_theory_1}. Sections \ref{sec:simulations} and \ref{sec:data} present simulation studies and a data analysis, respectively. 
Additional numerical results and all proofs are deferred to the Supplementary Material. 

\smallskip  

 \textit{Notation.} 
 Let $\mathbb{R}$ and $\mathbb{N}$ denote the sets of the real numbers and the natural numbers, respectively. 
 Given two sequences of real numbers $\{g_n\}$ and $\{h_n\}$, the notation $g_n \lesssim h_n$ means that there exists a constant $c > 0$ such that $g_n \leqslant c h_n$;  $g_n \asymp h_n$ means $h_n\lesssim g_n$ and $g_n\lesssim h_n$ simultaneously; 
$g_n \ll h_n$ indicates that $\lim_{n \to \infty} g_n / h_n = 0$. For vectors $x = (x_i)_{i=1}^n$ and $y=(y_i)_{i=1}^n \in \mathbb R^n$, define their inner product as 
$\langle x, y\rangle = \sum_{i=1}^n x_iy_i$,  the two norm as $\|x\|_2 = \sqrt{\langle x, x\rangle}$, and the infinity norm as $\|x\|_{\infty} = \max_{1\leqslant i \leqslant n}|x_i|$. For matrices $X = (x_{ij})_{1\leqslant i\leqslant n, 1\leqslant j\leqslant m}$ and $Y = (y_{ij})_{1\leqslant i\leqslant n, 1\leqslant j\leqslant m} \in \mathbb R^{n \times m}$, define their inner product as $\langle X,Y \rangle = \sum_{i=1}^n \sum_{j=1}^m x_{ij} y_{ij}$, the Frobenius norm as 
 $\|X\|_{\mathrm{F}} = \sqrt{\langle X, X \rangle}$, the operator norm as $\|X\|_{\operatorname{op}} = \sup_{\|v\|_2 = 1} \|Xv\|_2$, and the two-to-infinity norm as $\|X\|_{2 \to \infty} = \sup_{\|v\|_2 = 1} \|Xv\|_\infty$. %, and 
 % , and the entrywise maximum norm as $\|X\|_{\max} = \max_{1\leqslant i\leqslant n, 1\leqslant j\leqslant m} |x_{ij}|$.  
 For matrices $X_1,\dots,X_p$, let $\operatorname{blkdiag}(X_1,\dots,X_p)$ be the block‑diagonal matrix with $X_1,\dots,X_p$ on its diagonal.
Let $\mathrm{I}_k$ denote $k\times k$ identity matrix, and $\mathcal{O}(k)=\left\{Q \in \mathbb{R}^{k \times k}: Q Q^{\top}=\mathrm{I}_k\right\}$ represent the set of $k \times k$ orthogonal transformation matrices.
 For a square matrix $X \in \mathbb R^{n \times n}$, let $\operatorname{vech}(X) = (X_{12}, X_{13}, \ldots, X_{n-1,n})^\mytrans \in \mathbb R^{n(n-1)/2}$ be the vectorization of its strictly upper-triangular entries.

\smallskip  

\section{Latent Space Model} \label{sec:model_homo}

% \subsection{Model}
% For concreteness, we next assume that we observe 
% Assume we observe 
% Assume that we collect 
Consider 
an undirected network of $n$ nodes, encoded by a symmetric  adjacency matrix $A = (A_{ij})_{1 \leqslant i,j \leqslant n}$ 
where  
each $A_{ij}$ represents  the connection between two nodes $i$ and $j \in \{1,\ldots, n\}$ 
and $A_{ii} = 0$ indicating no self loops. 
For binary networks, $A_{ij}=1$ means that two nodes $i$ and $j$ are connected, and zero otherwise \citep{ma2020universal}. 
For weighted networks, $A_{ij}$ represents the weight of the edge connecting two nodes $i$ and $j$ \citep{he2023semiparametric,tian2024efficient}. 
% For concreteness and to match common  scenarios in applications, we assume that the network does not contain self loop, represented by  $A_{ii} = 0$ for all $1\leqslant i\leqslant n$, whereas the analysis can be straightforwardly generalized when the network contains self loops. 
% To accommodate flexible practical use, we allow $A_{ij}$ . 
In applications,  networks may contain binary, continuous, or count-valued edge weights. Latent space modeling assumes that each node $i$ is associated with a latent vector $z_i \in \mathbb R^k$ and a degree heterogeneity parameter $\alpha_i \in \mathbb R$. Then for any $1 \leqslant i < j \leqslant n$, 
 \begin{align} \label{eq:model}
    A_{ij} = A_{ji}\ \sim  \  p( \cdot  \mid  \Theta_{ij} ) \quad \text{ independently with } \quad \Theta_{ij} = \alpha_i + \alpha_j + \langle  z_i, z_j \rangle, 
\end{align} 
where $p( \cdot  \mid \theta )$ denotes the probability density/mass function  with a one-dimensional parameter $\theta$.  Let $Z = [z_1, \ldots, z_n]^\mytrans \in \mathbb R^{n \times k}$ and $\alpha = [\alpha_1, \ldots, \alpha_n]^\mytrans \in \mathbb R^{n}$ represent the parameters in matrix form. We can also write a matrix form $\Theta:=(\Theta_{ij})_{1\leqslant i, j\leqslant n} = ZZ^{\top}+\alpha 1_n^{\top} + 1_n\alpha^{\top}$ if we extend the definition of  $\Theta_{ij}$ to the case $i\geqslant j$ in the same way as in \eqref{eq:model}, 
% define $\Theta_{ij}$ with $i\geqslant j$ similarly as in \eqref{eq:model}, 
where $1_n$ is the all-ones vector in $\mathbb{R}^n$. Such a matrix form has also been discussed in \cite{ma2020universal}.

In this work, we assume that the observed data are generated by fixed true parameters $Z^\star=[z_1^{\star},\ldots, z_n^{\star}]^{\top}\in \mathbb{R}^{n\times k}$ and $\alpha^\star=[\alpha_1^{\star},\ldots, \alpha_n^\star]^{\top}\in \mathbb{R}^n$. 
% there exist fixed and true parameters $Z^\star=[z_1^{\star},\ldots, z_n^{\star}]^{\top}\in \mathbb{R}^{n\times k}$ and $\alpha^\star=[\alpha_1^{\star},\ldots, \alpha_n^\star]^{\top}\in \mathbb{R}^n$   that generate the observed data, and 
Our goal is to estimate these true parameters and quantify  the associated uncertainty. However, $(Z^\star, \alpha^\star)$ is not uniquely identifiable under the model \eqref{eq:model}. 
Given $(Z^\star, \alpha^\star)$, we can choose  any $Q\in \mathcal{O}(k)$ and $c\in \mathbb{R}^k$ to construct transformed parameters 
\begin{align}\label{eq:reparameterization}
 \tilde{z}_i=Q(z_i^{\star}-c)\quad  \text{and} \quad \tilde{\alpha}_i =\alpha_i^{\star}+(c^{\top} z_i^{\star} -\|c\|_2^2/2) \quad \text{for all}\ \   1\leqslant i\leqslant n, 
\end{align}
satisfying $\operatorname{vech}(\Theta^{\star}) = \operatorname{vech}(\tilde \Theta)$, i.e., $\alpha_i^{\star} + \alpha_j^{\star} + \langle  z_i^{\star}, z_j^{\star} \rangle= \tilde{\alpha}_i + \tilde{\alpha}_j + \langle  \tilde{z}_i, \tilde{z}_j \rangle$  for all  $1\leqslant i<j\leqslant n$.  
% Since the data distribution in  \eqref{eq:model} is determined by $\Theta_{ij}=\alpha_i + \alpha_j + \langle  z_i, z_j \rangle$, 
By the model formulation in \eqref{eq:model}, 
we know $(Z^\star, \alpha^\star)$ and $(\tilde{Z},\tilde{\alpha})$ yield  the same data distribution and therefore cannot be distinguished from each other. 
But generically,  this translation-rotation ambiguity is the only source of non-identifiability, shown by Proposition \ref{prop:identifiability} below. 
% But fortunately, the parameters are identifiable up to translation and rotation of latent vectors, which is formally established in Proposition \ref{prop:identifiability} below.  

% $[Z,\alpha]$ parameters giving the same data distribution can not be uniquely identified, as the model \eqref{eq:model} is determined by $\Theta_{ij}$. 

% Let $\Theta^{\star}=Z^{\star}Z^{\star\top} + \alpha^{\star}1_n^{\top}+1_n\alpha^{\star\top}$. 
% Given any $\tilde{Z}\in \mathbb{R}^{n\times k}$ and $\tilde{\alpha}\in \mathbb{R}^{n}$ such that the induced $\tilde{\Theta}=\Theta^{\star}$,   there exist $Q\in \mathcal{O}(k)$ and  $c\in \mathbb{R}^k$ such that 
% % $Z^{\star}=\tilde{Z}Q$
% $\tilde{Z}=(Z^{\star}-1_nc^{\top})Q$ 
%  and $\tilde{\alpha}=\alpha^{\star} + (Z^{\star} c-1_nc^{\top}c 1_n/2)$

% When $k+1\leqslant n+1/2-(8n+1)^{1/2}/2$, the model is generically globally identified in the sense that 
\begin{proposition} \label{prop:identifiability}
Assume $n>2(k+2)$. Then there exists a set $\mathcal{N}\subset \mathbb{R}^{n\times (k+1)}$ with zero Lebesgue measure such that, for any $[Z^\star, \alpha^\star]$ and $[\tilde{Z},\tilde{\alpha}] \notin \mathcal{N}$ satisfying $\operatorname{vech}(\Theta^{\star}) = \operatorname{vech}(\tilde \Theta)$, there exist  $Q\in \mathcal{O}(k)$ and  $c\in \mathbb{R}^k$ such that \eqref{eq:reparameterization} holds. 
\end{proposition} 

% Proposition \ref{prop:identifiability} only establishes generic identifiability because 
% If self-loops of networks also follow the model 
If self-loops $A_{ii}$ are observed and also follow \eqref{eq:model}, i.e., $A_{ii}\sim p(\cdot \mid \Theta_{ii})$, then the same identifiability holds without the need to exclude a zero-measure set. Proposition \ref{prop:identifiability} considers a more challenging scenario where the  diagonal $\Theta_{ii}$ information is unavailable, and therefore only generic identifiability can be established, which is consistent with the existing literature  \citep{bekker1997generic}. We focus on the no-self-loops scenario, since meaningful self-loops are often unavailable in real-world datasets, and the corresponding  theoretical  developments need to address additional technical caveats. 
% is more involved. 
% Though we highlight  that the extensions to self-loops is easy as adding more technical discussions 
% developments can be easily extended with self-loop, as  
More details are discussed in Section C of the Supplementary Material.
Proposition \ref{prop:identifiability} implies that true parameters yielding the same distribution form an equivalence class. 
To simplify the presentation, we introduce the following regularity conditions on $(Z^{\star}, \alpha^{\star})$, which will specify a convenient  representative in the equivalence class.

\begin{condition}[True parameters]\label{cond:truevalue}
Assume  $(Z^\star, \alpha^\star)$ satisfy: 
\begin{enumerate} [label=(\roman*), leftmargin=*]
%[ label=(\roman*),   widest=iii,  leftmargin=*] 
\item \ \  $1_n^{\mytrans}Z^{\star}=0$ and $Z^{\star \mytrans} Z^\star$ is diagonal. 

\item \ There exists a positive constant $M_{1}$ such that $\left\|Z^\star\right\|_{2 \to \infty} \leqslant  M_{1}$ and
 $\|\alpha^\star\|_{\infty} \leqslant M_{1}$.

\item  There exists a positive constant $M_{2}$ such that $\sigma_{\min}[Z^{\star \mytrans} Z^\star/n]
\geqslant M_{2}$, where $\sigma_{\min}(\cdot)$ represents the minimum singular value of the input matrix. 
 
\end{enumerate}
\end{condition} 

Condition~\ref{cond:truevalue} 
% imposes regularity on the true parameters, and  it
is comparable  to Assumptions II--IV in \citet{li2023statistical}, except that  we do not require that the limit of $Z^{\star \mytrans} Z^\star /n$ has  unique eigenvalues. % We next explain details in the condition. The two constraints in 
% Note Condition \ref{cond:truevalue} (i)  can be assumed without loss of generality.  
In  Condition~\ref{cond:truevalue}~(i), the constraint  $1_n^{\mytrans}Z^{\star}=0$  fixes the translational indeterminacy by centering the latent vectors, while the diagonality of $Z^{\star \mytrans} Z^\star$ specifies a particular coordinate system for the latent space. These two constraints  can be imposed without loss of generality, since the model is invariant to translations and rotations of the latent vectors by Proposition \ref{prop:identifiability}. 
% the center and coordinate of the latent space can be changed without impacting  the induced data distribution. 
Specifically, given any  $(Z^{\star},\alpha^{\star})$, we can construct $(\tilde{Z},\tilde{\alpha})$  following \eqref{eq:reparameterization} with $c^{\top}=1_n^{\top}Z^{\star}/n$ and $Q$ being an orthogonal matrix formed by right singular vectors of $Z^{\star}-1_nc^{\top}$. 
Then $(\tilde{Z},\tilde{\alpha})$ satisfies  Condition~\ref{cond:truevalue}~(i) and yields  the same data distribution as $(Z^{\star},\alpha^{\star})$ does. When $Z^{\star}-1_nc^{\top}$ has repeated singular values, the choice of $Q$ is not unique, but any choice gives an equivalent representative. 
% the diagonalization is not unique, but any such choice of $Q$ yields an equivalent representative.
% restricts that row sum of $Z^{\star}$ to satisfy a  row sum value and It . 
% If they are violated, one can center and rotate parameters by updating $z_i^{\star}$ and $\alpha_i^{\star}$ to $z_i^{\star}-\bar{z}$
% It as discussed above such that 
% The transformed parameters would satisfy (i) and yield the same data distribution. 
Moreover, 
Condition~\ref{cond:truevalue}~(ii)--(iii) assume  the parameters are  bounded and $Z^{\star\top}Z^{\star}/n$ is well-conditioned. These assumptions are common and standard in the related  literature \citep{young2007random,ma2020universal,zhang2020flexible,he2023semiparametric} and can facilitate theoretical analysis. 

\section{Gaps Between Existing Theory and Practice} \label{sec:gap_homo}

\subsection{Constrained Maximum Likelihood Estimator} \label{sec:gap_constrained_mle}
Under the latent space model \eqref{eq:model}, \citet{li2023statistical} establish the asymptotic distribution of 
 the  maximum likelihood estimator  under 
constraints.
% boundedness constraints of parameters. 
In particular, they consider  maximizing the log likelihood function $\sum_{1 \leqslant i < j \leqslant n} \ell(\Theta_{ij}; A_{ij}) $
 % where $M$ is a pre-specified constant, and 
with $ \ell(\theta; x) = \log p(x \mid \theta)$ and $\Theta_{ij}=\alpha_i + \alpha_j + \langle z_i,z_j \rangle$ under the constraints 
\begin{align}
    \|Z\|_{2 \to \infty } \leqslant M,\quad  \|\alpha\|_{\infty} \leqslant M, \quad  1_n^\mytrans Z = 0, \quad  \text{ and } Z^{\top}Z \text{ is diagonal}. \label{eq:constraint_original}
\end{align}
The first two inequalities in \eqref{eq:constraint_original} constrain the parameters to be bounded. 
% fall within bounded regimes. 
The zero-mean constraint on $Z$ removes the translational non-identifiability  as described in \eqref{eq:reparameterization}. 
Further restricting  $Z^{\top}Z$  to be diagonal sets a particular coordinate system for the latent space, under which $Z$ is identifiable up to sign flip when $Z^{\top}Z/n$ has distinct  eigenvalues. More technical explanations are provided in Section~\ref{sec:implicit_regular}. 
% makes $Z$ identifiable when  $Z^{\top}Z/n$ has distinct limiting eigenvalues. 
% To further address the unidentifiability from the orthogonal transformation, \cite{li2023statistical} proposed to restrict  $Z^{\top}Z$ to be diagonal. 
% the zero-mean and diagonal constraints ensure that the latent vectors are identifiable, up to sign flip. 
% The constraint resolves the identifiability of $Z$  up to an orthogonal transformation and 

% Directly applying 
Classical inferential theory is not directly applicable to analyze the log likelihood under the latent space model, because the likelihood function  is non-convex in $[Z,\alpha]$,  and the number of parameters grows with respect to the network size $n$ at an  unconventional rate (see Remark~\ref{rm:score_expansion_idea} for more details). 
To address these difficulties,  
\cite{li2023statistical} propose  augmenting the likelihood with a  
Lagrange multiplier penalty proportional to 
\begin{align}\label{eq:mutlipler_li}
\|\operatorname{vech}(Z^{\mytrans}Z)\|_{2}^2 +  \|Z^{\top} 1_n\|_{2}^2, 
\end{align} 
% \red{(YT: this is already defined in Notation)}, 
% and \eqref{eq:mutlipler_li}  is 
which is motivated from constraints in \eqref{eq:constraint_original}. 
% Intuitively, the two terms in \eqref{eq:mutlipler_li} aim to address the unidentifiability of $Z$ with respect to $\mathcal{O}(k)$ transformation and mean shift, respectively. 
Then \cite{li2023statistical} show the   Lagrange-adjusted likelihood function is strongly convex with high probability, thereby facilitating the proofs.  
More broadly, auxiliary Lagrange multiplier terms  have long been used in constrained maximum likelihood estimation and in likelihood analyses with singular information matrices, dating back to 
% the use of 
% auxiliary Lagrange multiplier terms in constrained maximum likelihood estimation 
% or likelihood with singular information matrices dates back to 
\cite{aitchison1958maximum} and \cite{el1994wald} 
% echo with but generalize classical studies 
in low-dimensional problems.  
% \citep{aitchison1958maximum,el1994wald}. 
\citet{wang2022maximum}   extends a  related idea to high-dimensional generalized factor models. 
Although the augmenting idea is shared, different model properties induce distinct technical challenges. 
Under the latent space networks, 
\cite{li2023statistical} further tackle with unique symmetric formulations  and carefully control residual errors with the specific penalty in \eqref{eq:mutlipler_li}.   
% address  under  the penalty in \eqref{eq:mutlipler_li}. 
% associated with the constraints in \eqref{eq:constraint_original}.  

% Nevertheless, 
Despite these advances for network models, 
the theoretical analysis  in \cite{li2023statistical} has several limitations from a practical perspective.   
 First, their  asymptotic theory requires  $Z^{\star\mytrans}Z^{\star}/n$ to converge  to a diagonal matrix with distinct eigenvalues, while additional constraints are needed to handle repeated eigenvalues. 
 % while aruging that extra constraints are needed to tackle with repeated eigenvalues.  
 % Such a strategy could make  practical  rules out practical situation where no eigengaps or eigen gaps diminish with respect to $n$.  
 Such a case-specific treatment could be   unsatisfactory in practice, 
 % hinders our understanding of practical situations 
 where prior spectral information is often unavailable, necessitating  a unified approach that accommodates various scenarios. 
% Such a case-by-case analysis could hinder our understanding about practical situation without prior knowledge, calling for a uniform conclusion embracing different cases. 
Second, the  upper bound $M$ in \eqref{eq:constraint_original} 
is required to be sufficiently large so that the underlying true parameters $[Z^{\star}, \alpha^{\star}]$ fall within the feasible regime. 
% depends on the underlying true parameters, \blue{(to add)} and is therefore generally unknown in practice. 
Since true parameters are typically unknown in practice, 
the implications of this constraint on practical usage and interpretation remain  ambiguous. 
Third, \cite{li2023statistical}   exclusively analyzes the constrained maximum likelihood estimator, whereas its  relationship to the actual outputs of practical algorithms  is  unclear.  

\vspace{-0.6em} 
\begin{remark}[Challenges in classical analysis]\label{rm:score_expansion_idea}
To make this paper self-contained, we briefly describe  the technical difficulties that are also mentioned in \cite{li2023statistical}.
% by the non-convexity of loss function.  
% Let $y\in \mathbb{R}^{n(k+1)}$ denote a vectorization of $[Z,\alpha]$, and  let $S_L(y)$ denote the gradient of $L(y) $. 
 % Let $S_L(Z,\alpha)$ and $S_{P^{\star}}(Z,\alpha)$ denote the gradient of $L(Z,\alpha)$ and $P^{\star}(Z,\alpha)$, respectively.  Intuitively,  
 % Let $S_L(y)$ and $S_{P^{\star}}(y) $ denote the gradient of $L(y) $ and $P^{\star}(y) $, respectively.  Intuitively,  
 For the simplicity of notation, let $y\in \mathbb{R}^{n(k+1)}$ denote the vectorization of $[Z,\alpha]$, and let $L(y)$ denote the negative log likelihood function, with its gradient denoted as $S_L(y)$.
The classical route to   deriving the distribution of a maximum likelihood estimator  utilizes  the first-order optimality condition and the Taylor expansion \citep{van2000asymptotic}. 
In our problem, that is, 
  $0 = S_L(\hat{y} ) = S_L(y^{\star} ) + H_L(y^{\star} )(\hat{y}-y^{\star})+ R(\hat{y}, y^{\star})$, where $R(\hat{y}, y^{\star})$ denotes the residual term. If  $R(\hat{y}, y^{\star})$ is uniformly small and the Hessian matrix  $H_L(y^{\star})$ is invertible, one might  obtain $\hat{y}-y^{\star}= - \{H_L(y^{\star})\}^{-1}S_L(y^{\star})+\text{small-order residuals}$,  serving as a starting point for  deriving asymptotic distributions of the  entries in  $\hat{y}-y^{\star}$ in canonical settings. However,  the leading term in $H_L(y^{\star})$ turns out to have exactly $k(k+1)/2$ zero eigenvalues, and the dimension of $y$, namely  $n(k+1)$,  is large compared to the  effective sample size $n(n-1)/2$, making it difficult to apply the standard argument. More details on the characterization of the null space of $H_L(y^{\star})$ are provided in Remark B.5 in the Supplementary Material. 
  % After augmenting the $L(y)$ with the penalty term in \eqref{eq:mutlipler_li}, \cite{li2023statistical} shows that the 
  %corresponding to the number of constraints for $Z^{\top}Z$ being diagonal.  
 \end{remark} 
 
% The first-order optimality condition $S_L(\hat{y})=0$, where $S_L(y)$ represents the gradient of $L(y)$. 
% our proof exploits the first-order optimality condition $S_L(\hat{y})=0$. 
% As argued in \cite{li2023statistical}, 
% classical routine of deriving the distribution of maximum likelihood estimator under convex likelihood function  \citep{van2000asymptotic}   expands
% $0 = S_L(\hat{y} ) = S_L(y^{\star} ) + H_L(y^{\star} )(\hat{y}-y^{\star})+ R(\hat{y}, y^{\star})$ 
% and relies on that   $R(\hat{y}, y^{\star})$ can be uniformly small and negligible and the Hessian matrix  $H_L(y^{\star})$ is full rank. 
% % $  S_L(\hat{Z},\hat{\alpha})$ leads to a singular Hessian matrix in proofs so that classical routine of expanding $0=S_L(\hat{Z},\hat{\alpha})$ at $(Z^{\star},\alpha^{\star})$ to study maximum likelihood estimator under convex likelihood function does not apply \citep{van2000asymptotic}. 
% % Technically, we utilize 
% However, one can show the leading term in $H_L(y^{\star})$ has exactly $k(k+1)/2$ zero eigenvalues, corresponding to the number of constraints for $Z^{\top}Z$ being diagonal. 

 \subsection{Projected Gradient Descent and Universal Singular Value Thresholding} \label{sec:pgdreview}
 Practically, to compute the constrained maximum likelihood estimator, \cite{li2023statistical} adopt a common strategy in the existing literature: a singular value
thresholding procedure as in Algorithm 3 of \cite{ma2020universal} and a projected gradient descent method as in Algorithm 1 of \cite{ma2020universal} for initialization and
optimization, respectively. Such a two-stage approach   has  been widely used and generalized in various network models \citep{zhang2020flexible,he2023semiparametric,tian2024efficient}. 

The projected gradient descent method (Algorithm 1 in \cite{ma2020universal}) estimates the latent vectors under a non-convex optimization framework. This approach  iteratively  updates $Z$ and $\alpha$ parameters along their gradient directions and then projects the updates onto a feasible regime  like \eqref{eq:constraint_original} to enforce boundedness and the centering constraint. 
% , including the boundedness and zero-mean constraints. 
% \cite{ma2020universal} established \citet{ma2020universal} proposed a projected gradient descent procedure to solve the non-convex optimization problem. It includes updates along the gradient directions of $Z$ and $\alpha$ parameters, followed by projection of estimates onto constraint sets.  
\citet{ma2020universal} establish the  convergence of the algorithmic output  to the true latent parameters under the Frobenius norm. 
However, their analysis has two key limitations when integrated with the inferential results in \cite{li2023statistical}. 
% There are several limitations of the analysis in \cite{ma2020universal} when applied to making connections with results in \cite{li2023statistical}. 
First, \cite{ma2020universal} evaluate the   discrepancy between the algorithmic  output and the true parameters $[Z^{\star},\alpha^{\star}]$, which essentially conflates both algorithmic  convergence error and the irreducible statistical error. To apply the inferential results of \cite{li2023statistical} in practice, these two sources of errors need to be properly  decoupled, necessitating a more refined analysis.  
% Their techniques cannot be directly applied to sharply characterizes  the difference between the MLE $[\hat{Z},\alpha]$ and $[Z^{\star},\alpha^{\star}]$.  
Second,  the theoretical arguments in \cite{ma2020universal} require 
 hyperparameters to be chosen in a way  relying on  the unknown true parameters, such as the boundedness set in the projection step and the learning rate in the projected gradient descent. These requirements are often impractical and limit  the developments of    reliable and fully data-driven  inference. 

% one has to properly specify the boundedness constraint, which depends on the unknown true parameters through $M_1$. 
%See more details in Remark \ref{rm:ma_eta_restricction}.  Second, its conclusion on the projected gradient descent examines the convergence of the output of algorithm, which does not directly exhibit connection with  the theoretical maximum likelihood estimator due to instrinsic non-convexity. 
% output from the projected does not show direct connection with the theoretical maximum likelihood estimator, so that 
% As a result, it is unclear whether the asymptotic distributional results in \cite{li2023statistical} can be directly applied or not. 

In addition, the initialization procedure in Algorithm 3 of \cite{ma2020universal} is built  on the universal singular value thresholding (USVT) method introduced by \citet{chatterjee2015matrix}.
% whose theoretical guarantee adopted from \cite{chatterjee2015matrix}.
Both \cite{ma2020universal} and \cite{chatterjee2015matrix} consider  bounded entries in adjacency matrices, such as binary edges in networks. 
Moreover, \cite{ma2020universal} suggest   a projection  step to regularize the range of estimates, similarly to the projected gradient descent. 
% The original proposal in \cite{chatterjee2015matrix} requires a truncation step . 
However, when the true parameters are unknown, 
% adjacency matrix contains unbounded entries, e.g., under Gaussian or Poisson models, 
% it is unclear how a projection set should be chosen in a principled way without relying on unknown true parameters.  
it is unclear how to choose the projection set in a principled way without prior knowledge of the true parameters under general data distributions.  
% Consequently, The applicability of the USVT-based initialization and its accompanying projection strategy to unbounded-edge models remains insufficiently understood.
% whether truncation is still necessary and, if so, how it should be chosen in a principled way.

% whether and how truncation should be implemented remains underexplored. 

To address the above issues in the existing literature,  we develop a comprehensive  analysis of both the maximum likelihood estimator and   the practical estimation procedure, covering the USVT-based initialization and the projected gradient descent algorithm. 
These contributions  are detailed in Sections \ref{sec:mle_theory_1} and \ref{sec:gd_theory_1}, respectively.   

%Our results can relax 

%  

% how to properly truncate without relying on knowledge 

% In practice, it 

% However, following the theoretical analysis in \cite{ma2020universal}, one needs to 

% is often found that a projection step is unnecessary \citep{zhang2020flexible}. Instead, ... (compact description and refer to Algorithm \ref{algor:pgd}.)
% Although \citet{ma2020universal} provided a theoretical error bound for their algorithm output under F norm, the relationship between this algorithmic estimator and the constrained MLE $(\hat{Z},\hat{\alpha})$ remains unclear.
% \blue{(YH: Maybe slightly change the name of the algorithm to highlight differences with \cite{ma2020universal}.}

% \blue{(YH: May update this to a paragraph of description to be more compact and combine with USVT. Decide later.)} 

% \blue{(YH: Why not add $\lambda$? May need explanations.)}

% \begin{remark}
% gradient update in \cite{ma2020universal} is not  $\dot{L}$ exactly. They differ by entries from the diagonal. (\blue{YH: We may argue that  our proof can tackle both  later? Not important now. Just keep in mind.}) 
% \end{remark}

% \subsection{New Eigengap-Free Lagrange Adjustment}

% \section{Eigengap-Free Theory for the Maximum Likelihood Estimator} \label{sec:mle_theory_1}

\section{Maximum Likelihood Estimator 
with Flexible Eigenvalue Multiplicity} \label{sec:mle_theory_1}
% under  Unrestricted  Spectrum}  

In this section, we  examine the constrained maximum likelihood estimator and its asymptotic theory without imposing the unique-eigenvalue assumption on $Z^{\star\top}Z^{\star}/n$ used in \cite{li2023statistical}. 
Specifically, we let $(\hat Z, \hat{\alpha})$ be one of 
% which is \red{not unique?} (one of) 
the solution that maximizes the log likelihood function, or equivalently, any solution to the following optimization problem
\begin{equation}
    \begin{aligned} \label{eq:prob}
\mathop{\arg\min}_{Z \in \mathbb R^{n \times k}, \alpha \in \mathbb R^n}\ \  \quad & L(Z,\alpha):=-\sum_{1 \leqslant i < j \leqslant n} \ell(\Theta_{ij}; A_{ij})  \\
\text { subject to } \quad\quad   & \|Z\|_{2 \to \infty } \leqslant M,\ \|\alpha\|_{\infty} \leqslant M, \ 1_n^\mytrans Z = 0.
\end{aligned}
\end{equation}
Different from  the constraints in \eqref{eq:constraint_original}  imposed by \cite{li2023statistical}, \eqref{eq:prob} does not require $Z^{\top}Z$ to be diagonal and thus is more flexible in practice. Nevertheless, their proposed Lagrange multiplier  \eqref{eq:mutlipler_li} cannot be used,   necessitating new  developments  in theoretical analysis. 

% restrict $\hat{Z}^{\mytrans}\hat{Z}$ to be diagonal. 
% Therefore, $\hat{Z}$ can only be identified up to $\mathcal{O}(k)$ transformation, as discussed in Section \ref{sec:model_homo}. 

% Therefore $\hat{Z}$ actually denotes an equivalence class of estimators up to   $\mathcal{O}(k)$ transformation, as discussed in Section \ref{sec:model_homo}.  
% \subsection{Conditions}
% \noindent \textbf{Inferential goal: } 
% Throughout the sequel, we use $(Z^\star, \alpha^\star)$ to denote the true values of $(Z,\alpha)$ that generate the observed data. 
% Notably, true values of latent embeddings should be an equivalent class. (Although in practice, preferred rotations may be used   for unique definition and interpretation (could cite varimax)) Given an estimate of latent vectors $\hat{Z}$, 

% \blue{(YH: Add somewhere to emphasize non-convexity of $L$ and high-dimensionality)} Remark~\ref{rm:xi_null} and \ref{rm:compareli}. 

\subsection{Implicit regularization from orthogonal Procrustes problem}\label{sec:implicit_regular}

% We explain a fundamental limitation of the idea in \cite{li2023statistical},  
We begin by explaining a fundamental limitation in the existing analysis, which in turn motivates our new strategy.  
As reviewed in Section~\ref{sec:gap_constrained_mle},
\cite{li2023statistical} address the non-identifiability of $Z$ by 
% To address that, \cite{li2023statistical} propose 
restricting $1_n^{\top}Z=0$ and $Z^{\top}Z$ to be diagonal, which are  intended  to resolve the ambiguities from arbitrary mean shift and $\mathcal{O}(k)$ transformation, respectively. % discussed in \eqref{eq:reparameterization}. 
% The first centering constraint removes the ambiguity from arbitrary mean shift in \eqref{eq:reparameterization}, and the second tries to  
% It turns out that 
Although  
the centering constraint effectively removes  the mean-shift ambiguity,  the diagonality restriction does not fully resolve the  non-identifiability from $\mathcal{O}(k)$ transformations. 

To see this, consider two candidate matrices $[Z^{\star},\alpha^{\star}]$ and $[\tilde{Z},\tilde{\alpha}] \in \mathbb{R}^{n\times (k+1)}$ satisfying $1_n^{\top}Z^{\star}=1_n^{\top}\tilde{Z}=0$ and both yield the same likelihood.
By Proposition~\ref{prop:identifiability}, we know \eqref{eq:reparameterization} holds up to a zero measure set. 
% By $1_n^{\top}Z=1_n^{\top}\tilde{Z}=0$ ,
% Since both $Z^{\star}$ and $\tilde{Z}$ are centered, 
% we must have
The centering constraint forces 
$c=0$, so there exists $Q\in \mathcal{O}(k)$ such that $\tilde{Z}Q=Z^{\star}$. 
% i.e., $\alpha 1_n^{\top}+1_n\alpha^{\top}+ZZ^{\top}=\tilde{\alpha} 1_n^{\top}+1_n\tilde{\alpha}^{\top}+\tilde{Z}\tilde{Z}^{\top}$ (diagonal elements added for simplicity of illustraion like discussed after Proposition~\ref{prop:identifiability}).   
Now further assume both $Z^{\star\top}Z^{\star}$ and $\tilde{Z}^{\top}\tilde{Z}$ are diagonal with diagonal entries in non-increasing order. 
% If further assume both $Z^{\top}Z$ and $\tilde{Z}^{\top}\tilde{Z}$ are diagonal and 
If both have unique eigenvalues, $Q$ must equal  a   signature matrix $\mathrm{diag}(q_1,\ldots, q_k)$ for $q_i\in \{-1,+1\}$, implying that   columns of $Z$ and $\tilde{Z}$  can be matched up to sign flips.  
In contrast, if $Z^{\star\top}Z^{\star}$ and $\tilde{Z}^{\top}\tilde{Z}$ are proportional to the identity matrix, then their eigenvalues are repeated, 
% diagonal with repeated eigenvalues,  
and $Q\in \mathcal{O}(k)$ is not necessarily diagonal. As a result,   there is no clear one-to-one correspondence between the columns of $Z^{\star}$ and $\tilde{Z}$. % can be any matrix in $\mathcal{O}(k)$. 
More generally, when only a subset of eigenvalues of $Z^{\star\top}Z^{\star}$ (or $\tilde{Z}^{\top}\tilde{Z}$) are repeated, only the columns associated with distinct eigenvalues are identifiable up to sign, whereas the remaining columns are not. 
% the columns of $Z^{\star}$ and $\tilde{Z}$ corresponding  to distinct eigenvalues can be equal up to sign flip whereas the others cannot.  

% We illustrate the issue discussing three cases. 
% First,  $Z^{\star\top}Z^{\star}/n$ has distinct  (limit) eigenvalues. 
% Then for any two candidate matrices $Z_1^{\star}$ and $Z_2^{\star}$ such that $Z_1^{\star}Z^{\star}_1=$  
% When $Z^{\top}Z/n$ has repeated (limit) eigenvalues,  

The preceding discussion shows that directly examining the algebraic difference between two latent-vector   matrices  would 
 inevitably require additional restrictions and cumbersome case-by-case discussions, thereby obscuring the practical implications. 
% , say MLE $\hat{Z}$ and true parameters $Z^{\star}\in \mathbb{R}^{n\times k}$, may not be a reasonable target. 
% It would inevitably require extra restrictions and messy case-specific discussions to match two matrices, causing the practical implications hard to track. 
% It requires 
Instead, to obtain a unified characterization   without unnecessary restrictions on the spectrum of $Z^{\star}$, 
% To quantify the discrepancy between $\hat{Z}$ and $Z^{\star}$, instead of directly taking their difference, 
it is more natural to measure discrepancy
% compare the two latent-vector matrices only 
% study the difference 
% between $\hat{Z}$ and $Z^{\star}$ 
only up to the best orthogonal transformation in $\mathcal{O}(k)$.  
% For the solution $(\hat{Z},\hat{\alpha})$ of problem ..., we define $\hat{Z}_q = \hat{Z} \hat{Q}^\mytrans$ with $\hat{Q} = \arg\min_{Q \in \mathcal{O}(k)} \|\hat{Z} - Z^\star Q\|_{\F}$.  
In particular, 
given any minimizer $(\hat{Z},\hat{\alpha})$ of the constrained optimization problem \eqref{eq:prob},  define 
\begin{align}\label{eq:rotation_z_def}
    \hat{Z}_q = \hat{Z} \hat{Q}^\mytrans, \quad \text{ with }\, \hat{Q} =  \underset{Q\in \mathcal{O}(k)}{\arg\min}  \|\hat{Z} - Z^\star Q\|_{\F}. 
\end{align} 
% \red{or???}
% \begin{align}\label{eq:rotation_z_def}
%     \hat{Z}_q = \hat{Z} \hat{Q}, \quad \text{ with }\, \hat{Q} = \arg\min_{Q \in \mathcal{O}(k)} \|\hat{Z}Q - Z^\star \|_{\F}. 
% \end{align}
% By the problem formulation, $[\hat{Z},\hat{\alpha}]$ and $[\hat{Z}_q,\hat{\alpha}]$ yield the same likelihood $L(\hat{Z},\hat{\alpha})=L(\hat{Z}_q,\hat{\alpha})$ and $\|\hat{Z}\|_{2\to \infty}=\|\hat{Z}_q\|_{2\to\infty}$. 
We propose to examine   $\hat{Z}_q-Z^{\star}=\hat{Z} \hat{Q}^\mytrans-Z^{\star}$ and establish  asymptotic distributions of its entries. 
Since $ \hat{Z}_q $ differs from  $ \hat{Z} $ only by  $\hat{Q}\in \mathcal{O}(k)$, they  yield  the identical  likelihood by  the arguments in Section~\ref{sec:model_homo} and hence give the same statistical interpretation. 
Meanwhile, the aligned difference $\hat{Z}_q-Z^{\star}$ 
removes ambiguities from sign flips or repeated eigenvalues, so it resolves the non-identifiability from $\mathcal{O}(k)$ transformations fundamentally. 
% since the ambiguity from sign flip or repeated eigenvalues is removed. Therefore  
Conceptually, this formulation eliminates the need to assume unique limiting eigenvalues 
% allows us to remove the requirement of unique limiting eigenvalues 
of   $Z^{\star\mytrans}Z^{\star}/n$ in theory.  
Further technical discussions are provided in Remark~\ref{rm:q_discussion} below. 
% \red{(JS: it seems that this is the last place where we mentioned the no-eigengap thing in this section. do you think we still want to mention it in Remark 2, or after Condition 2?)}

% without requiring   $Z^{\star\mytrans}Z^{\star}/n$ to have unique limiting eigenvalues. 
% This is consistent with the discrepancy measure in other models such as the random dot product graphs  \citep{xie2023efficient}. 
% The ambiguity from sign flip or common eigenvalues is removed, which intuitively, allows asymptotic theory established without an explicit restriction on the spectrum  of $Z$ in \eqref{eq:prob}. 

% By our problem formulation,  $[\hat{Z}_q,\hat{\alpha}]$ optimizes  the constrained  loss in \eqref{eq:prob} if and only if $[\hat{Z},\hat{\alpha}]$ achieves that, because they yield the same likelihood $L(\hat{Z},\hat{\alpha})=L(\hat{Z}_q,\hat{\alpha})$, and satisfy the constraints in \eqref{eq:prob} simultaneously. 
% $1_n^{\top}\hat{Z}_q=1_n^{\top}\hat{Z}\hat{Q}^{\top}=0$ and $\|\hat{Z}_q\|_{2\to\infty}\leqslant \|\hat{Z}\|_{2\to \infty}\leqslant M$  are also satisfied.  
% \red{(YH: is if and only if correct?)} 

Technically, without restricting the diagonality of $Z^{\top}Z$, it may appear that the non-convexity issue discussed in Section~\ref{sec:gap_constrained_mle} still persists. 
But interestingly, the best $\mathcal{O}(k)$ alignment induces an implicit regularization in the sense that for any two matrices $\hat{Z}$ and $ Z^{\star}\in \mathbb{R}^{n\times k}$, $\hat{Z}_q$  defined   in \eqref{eq:rotation_z_def} always satisfies 
\begin{align} \label{eq:equal_symmetricity}
    \hat{Z}_q^{\top}\, Z^{\star}=Z^{\star\top}\, \hat{Z}_q. 
\end{align}
% where $\hat{Z}_q$ is defined as in \eqref{eq:rotation_z_def} and it  rotates    $\hat{Z}$ for   maximal agreement with   $Z^{\star}$. 
% and  is defined as in \eqref{eq:rotation_z_def}. 
% \begin{align} 
%     (\hat{Z}\hat{Q}^{\top})^{\top} Z^{\star}=Z^{\star\top} (\hat{Z}\hat{Q}^{\top}), %\hspace{2em}\quad \hat{Q}=\argmin_{Q\in \mathcal{O}(k)} \|Z-Z^{\star}Q\|_{\F}. 
% \end{align}
% where $\hat{Q}$ rotates $Z^{\star}$ for  maximal agreement with   $\hat{Z}$ and  is defined as in \eqref{eq:rotation_z_def}. 
This identity follows from classical results  on the orthogonal Procrustes problem and is originally shown in \cite{ten1977orthogonal}. 
% inducing an equivalent likelihood and satisfies \eqref{eq:equal_symmetricity}. 
Thus, when studying latent vectors up to the best $\mathcal{O}(k)$ alignment, one can impose an extra constraint \eqref{eq:equal_symmetricity} without loss of generality, which addresses the non-identifiability in a more natural way.  
% To exploit the above constraint technically, 
In the proof, 
we  introduce a  new Lagrange multiplier penalty 
$ P^{\star}(Z) := \|\operatorname{vech}(Z^{\mytrans}Z^{\star}-Z^{\star\top}Z)\|_{2}^2 +  \|Z^{\top} 1_n\|_{2}^2$, 
corresponding to two sets of constraints   $Z^{\mytrans}Z^{\star}=Z^{\star\top}Z$ and $Z^{\top} 1_n=0$ to address the ambiguity from $\mathcal{O}(k)$ transformation and mean shift, respectively.  
Our proof 
% utilizes $P^{\star}(Z)=0$ for   $Z$ satisfying those constraints and 
shows that the 
% (original null space in the leading term of the Hessian matrix of $L$ can be ) 
Hessian matrix of $L(Z,\alpha)+P^{\star}(Z)$ is non-singular with high probability, thereby overcoming the non-convexity challenge explained in Remark~\ref{rm:score_expansion_idea}. 
We defer further technical explanations to  Remark~\ref{rm:score_augment} below. 

The proposed implicit constraints and  the augmenting penalty   have  several notable features.    
% number for $Z^{\top}Z$ being diagonal. 
% and consider the Lagrange-augmented problem  
% \begin{align}\label{eq:prob_aug}
%     \argmin_{Z\in \mathbb{R}^{n\times k},\alpha\in \mathbb{R}^{n}}L(Z,\alpha)+P^{\star}( Z,\alpha)
% \end{align}
% Like argued in \cite{li2023statistical}, 
First, the penalty $P^{\star}(Z)$ relies on unknown true parameter $Z^{\star}$ and is neither computable nor usable in practice. 
% We emphasize that the penalty term 
It only  serves as an intermediate technical device, rather than a practical regularization term such as the lasso or ridge penalty.  
% is different from penalties in shrinkage regressions, such as lasso or ridge. 
% Instead, it only serves as an intermediate technical tool. % to help examine a candidate $Z\in \mathbb{R}^{n\times k}$ after best alignment with $Z^{\star}$. 
% Specifically, 
% our proof utilizes that the gradient of $P^{\star}(Z)$ is zero, while its Hessian's column space covers the null space of the leading of $H_L$ ... .
Second, by the construction,  the constraints $Z^{\mytrans}Z^{\star}=Z^{\star\top}Z$ and $Z^{\top}1_n = 0$ contribute a total of  $k(k+1)/2$ number of constraints, which coincides with the  dimension of the intrinsic null space of the Hessian matrix of $L(Z,\alpha)$ mentioned in Remark~\ref{rm:score_expansion_idea}. 
Third, 
even when $Z^{\star\top}Z^{\star}/n$ has repeated limiting eigenvalues, 
we show that the leading term in the Hessian matrix of $L(Z,\alpha)+P^{\star}(Z)$ is still non-singular with high probability. In contrast,  using the augmenting penalty  in \cite{li2023statistical}, i.e., replacing $P^{\star}(Z)$ by \eqref{eq:mutlipler_li}, 
% $L(Z,\alpha)+P(Z)$, \
the corresponding leading term  has exactly zero eigenvalues. 
% our conclusion of non-singularity always holds, whereas 
% our proof does not require $Z^{\star}Z^{\star}/n$ to have  distinct limiting eigenvalues, whereas 
% we can show that using the penalty in \cite{li2023statistical} the leading term in the Hessian matrix of augmented loss function $L(Z,\alpha)+P(Z)$ still has zero eigenvalues. 
Details are provided in Remark B.5 of the Supplementary Material.

\vspace{2pt}

\begin{remark}[Uniqueness of $\hat{Z}_q-Z^{\star}$]  \label{rm:q_discussion}
% Even when $\hat{Z}$ and $Z^{\star}$ are only defined up to $\mathcal{O}(k)$ transformations, the aligned difference  $\hat{Z}_q-Z^{\star}$ is unique. 
First, given two fixed matrices $\hat{Z}$ and $Z^{\star}\in \mathbb{R}^{n\times k}$, 
$\hat{Q}=UV^{\top}$ with $U\Sigma V^{\top}$ denoting the singular value decomposition of $Z^{\star \top}\hat{Z}$   is always one of the solutions to $\min_{Q\in \mathcal{O}(k)}\|\hat{Z}-Z^{\star}Q\|_{\F}$ \citep{gower2004procrustes}. 
This minimizer is unique 
% $\hat{Q}$ defined in \eqref{eq:rotation_z_def} 
% is unique 
if $Z^{\star \top}\hat{Z}$ is non-singular. We next focus on this regime, because   
we show $Z^{\star \top}\hat{Z}$ is non-singular with high probability in the proof of Theorem~\ref{thm:homoMLE} (see Remark D.1 in the Supplementary Material). 
% \red{(YH: Is it possible to add a remark we can refer to? Or some equation that is already in the proof. Or we just add $\hat{Q}$ is unique in Theorem 1 along with other uniqueness.)} 
% Otherwise, $\hat{Q}=UV^{\top}$ with $U\Sigma V^{\top}$ denoting the singular value decomposition of $Z^{\star \top}\hat{Z}$   is always one of the minimizer of $\|\hat{Z}-Z^{\star}Q\|_{\F}$. 
% Thus, we focus on this unique $\hat{Q}$    
% It admits an analytical expression based on the solution to the orthogonal Procrustes problem  . 
% It is unique if $\hat{Z}^{\top}Z^{\star}$ is non-singular,  which will be shown to hold with high probability in our proof.  
% By \eqref{eq:rotation_z_def},
Then we can interpret $\hat{Q}$ as the  orthogonal transformation that best aligns     $Z^{\star}$ with $\hat{Z}$, and symmetrically,  $\hat{Q}^{\top}$ best aligns $\hat{Z}$ with $Z^{\star}$ in the sense that   $\hat{Q}^{\top}=  \argmin_{Q\in \mathcal{O}(k)} \| \hat{Z}Q - Z^\star \|_{\F}$. 
Second,  the aligned difference  $\hat{Z}_q-Z^{\star}$ remains unchanged if replacing $\hat{Z}$ with $\hat{Z}\tilde{Q}$ for any $\tilde{Q}\in \mathcal{O}(k)$ when $Z^{\star \top}\hat{Z}$ is non-singular. 
% Second, when $\hat{Z}$ is only identified up to $\mathcal{O}(k)$ transformations, the aligned difference  $\hat{Z}_q-Z^{\star}$ remains unchanged. 
% Specifically, if replacing $\hat{Z}$ with $\hat{Z}\tilde{Q}$ for any $\tilde{Q}\in \mathcal{O}(k)$, the updated $\hat{Z}_q-Z^{\star}$ based on \eqref{eq:rotation_z_def} remains same as before \citep{gower2004procrustes}. 
% Consider the equivalence classes represented by $\hat{Z}$ and $Z^{\star}$ up to $\mathcal{O}(k)$ transformations. 
% , since it is invariant to replacing $\hat{Z}$ with $\hat{Z}Q$ for any $Q\in \mathcal{O}(k)$.  
% Specifically, the identifiability up to the $\mathcal{O}(k)$ transformation  
% For example, let $\hat{Z}_qQ_1$ be another equivalent e..
% since two matrices are equal up to $\mathcal{O}(k)$ defines a equivalence relation,  
 % each $\hat{Z}$  has an associated equivalence class $\{Z\in \mathbb{R}^{n\times k}: Z=\hat{Z} Q,\, Q\in \mathcal{O}(k)\}$. % represent 
Therefore, in such a non-singular regime relevant to our analysis, $\hat{Z}_q-Z^{\star}$ is a well-defined quantity to examine the difference between a fixed $Z^{\star}$ and the equivalence class of $\hat{Z}$ with respect to $\mathcal{O}(k)$ transformations.   
% If singularity occurs in practice, 
% uniquely defined and characterizes the intrinsic difference between the two equivalence classes represented by $\hat{Z}$ and $Z^{\star}$ with respect to $\mathcal{O}(k)$ transformations.    
% $\hat{Z}$ can represent a whole equivalence class or any representative element in the equivalence class. Regardless to that, the difference.      
\end{remark} 

\begin{remark}[Details for score augmentation] \label{rm:score_augment}
To shed light on theoretical derivations, 
we provide a sketch on how we address the   non-convexity issue highlighted in Remark~\ref{rm:score_expansion_idea} via score augmentation.  
Let $y$ denote the vectorization of $[Z,\alpha]$ like  in Remark~\ref{rm:score_expansion_idea}, and define  $P^{\star}(y)=P^{\star}(Z)$ with its gradient denoted as $S_{P^{\star}}(y)$.  Our argument consists of three parts.  
% proceeds in three steps. 
First, we show $S_L(\hat{y})=0$ if and only if $S_L(\hat{y}_q)=0$ where $\hat{y}_q$ denotes the vectorization of $[\hat{Z}_q,\hat{\alpha}]$ defined in \eqref{eq:rotation_z_def}. 
Second, we argue $S_{P^{\star}}(\hat{y}_q)=0$  based on the implicit regularization \eqref{eq:equal_symmetricity}. Third,    the first-order condition $S_L(\hat{y})=0$ in Remark~\ref{rm:score_expansion_idea} can therefore be equivalently reformulated as the augmented score equation $ S_L(\hat{y}_q ) + S_{P^{\star}}(\hat{y}_q)=0$. We will then show that the augmented score has a non-singular Hessian matrix $H(y^{\star})$ with high probability,  without imposing restrictions on the spectrum of $Z^{\star\top}Z^{\star}$. Heuristically, this enables an expansion $\hat{y}_q-y^{\star}=-\{H(y^{\star})\}^{-1}S_L(y^{\star})+\text{small-order residuals}$ to examine the asymptotic distribution of $\hat{Z}_q-Z^{\star}$ and $\hat{\alpha}-\alpha^{\star}$. 
% Our proof utilizes that the gradient of $P^{\star}(y)$ is zero, i.e., $S_{P^{\star}}(y)=0$, when $Z$ satisfies the constraints $Z^{\mytrans}Z^{\star}=Z^{\star\top}Z$ and $Z^{\top} 1_n=0$. 
% Therefore, we can augment the first-order condition $0=S_L(\hat{y})$ in Remark~\ref{rm:score_expansion_idea} to 
% \begin{align}\label{eq:first_score_0}
%     0=S_L(\hat{y})=S_L(\hat{y}_q ) = S_L(\hat{y}_q ) + S_{P^{\star}}(\hat{y}_q)
% \end{align}
% where the second equation follows by the equivalence of likelihood function, $\hat{y}_q$ ... and the third equation follows by $S_{P^{\star}}(\hat{y}_q)=0$ for any $\hat{Z}_q$ satisfying \eqref{eq:equal_symmetricity}. 
% With the extra term, the induced penalized score has non-singular Hessian with high probability. 
% Further discussions on the column spaces of Hessian matrices are provided in Remark  \ref{rm:compareli} in the Supplementary Material.  
% Intuitively, we bypass two-step procedure 
Although this score expansion is motivated from classical maximum likelihood theory, our problem has  specific properties that  pose unique technical difficulties. 
 For example, the first-order condition is not directly available due to the constraints in \eqref{eq:prob}, the Hessian matrix has a special structure, and the residual terms are  high-dimensional. Our proof carefully   overcomes all the technical issues arising from our unique problem properties. 
 Notably, our analysis applies to any solution of \eqref{eq:prob} without  the need to diagonalize $\hat{Z}^{\top}\hat{Z}$, and thus  differs  significantly from \cite{li2023statistical}.  
 % and thus differs from \cite{li2023statistical} requiring a diagonalization of $\hat{Z}$. 
 % since our problem formulation differs from that  \cite{li2023statistical},  
\end{remark}

% \vspace{1pt}  

\subsection{Asymptotic Theory for Constrained Maximum Likelihood Estimator} 
We next develop asymptotic theory for the constraint-relaxed  estimator in \eqref{eq:prob}.  
 For the  formal theoretical developments, we impose the following regularity condition on the  edgewise distribution.

% \begin{condition}\label{cond:truevalue}
% Assume the true parameters $(Z^\star, \alpha^\star)$ satisfy: 
% \begin{itemize}\setlength{\itemsep}{0pt}
%     \item[(i)] There exists a positive constant $M_{1}$ such that $\left\|Z^\star\right\|_{2 \to \infty} \leqslant  M_{1}$ and
%  $\|\alpha^\star\|_{\infty} \leqslant M_{1}$.

% \item[(ii)] There exists a positive constant $M_{2}$ such that $\sigma_{\min}[Z^{\star \mytrans} Z^\star/n]
% \geqslant M_{2}$. \blue{(define $\sigma_{\min}$)}
 
% \item[(iii)] $1_n^{\mytrans}Z^{\star}=0$ and $Z^{\star \mytrans} Z^\star$ is diagonal. 
% \end{itemize}
% \end{condition}
% \smallskip 

\begin{condition}[Edgewise distributions] \label{cond:parfunction}
Let $\mathcal{X} = \{x \in \mathbb R: p(x \mid \theta) > 0\}$  denote the support of $p(x \mid \theta) $. Assume $\ell(\theta;x)=\log p(x \mid \theta)$ in \eqref{eq:model} satisfies the following conditions.
\begin{itemize}\setlength{\itemsep}{0pt}
\item[(i)]  For any fixed  $x \in \mathcal{X}$, $\ell(\theta;x)$ is three times differentiable with respect to $\theta$, with its first to third derivatives with respect to $\theta$ denoted by   $\ell^{\prime}(\theta;x)$, $\ell^{\prime \prime}(\theta;x)$, and $\ell^{\prime \prime \prime}(\theta;x)$, respectively.
Moreover, for any constant $b > 0$, there exist constants $\kappa_1(b), \kappa_2(b), \kappa_3(b) > 0$ such that
\[
\kappa_1(b) \leqslant -\ell''(\theta;x) \leqslant \kappa_2(b) \quad \text{ and } \quad |\ell'''(\theta;x)| \leqslant \kappa_3(b)
\quad\text{for all } \theta \in [-b,b] \text{ and } x \in \mathcal{X}.
\] 
% , and $\ell(\theta;x)$ is a fixed function that does not vary with respect to $n$. 
% Moreover, there exist  positive sequences $1 \ll \kappa_{1,n} \ll \kappa_{2,n} \ll \kappa_{3,n} \ll \kappa_{4,n} \ll \kappa_{5,n} \lesssim \log(n)$ such that for each $s = 1,2,3,4$, the bounds $1/\kappa_{s+1,n} \leqslant -\ell^{\prime\prime}(\theta;x) \leqslant \kappa_{s+1,n}$ and $|\ell^{\prime\prime\prime}(\theta;x)| \leqslant \kappa_{s+1,n}$ hold for any $x \in \mathcal{X}$ and $|\theta| \leqslant \kappa_{s,n}$. 
 
\item[(ii)] There exist constants $K, s > 0$ such that for any $t \geqslant 0$, $\Pr(|\ell^{\prime}(\Theta_{ij}^\star;A_{ij})| > t) \leqslant 2 \exp(-(t/K)^{s})$. 
\end{itemize}
\end{condition}

% Condition \ref{cond:truevalue} is similar to Assumptions II--IV in \citet{li2023statistical}, while 

% \sout{Condition \ref{cond:parfunction} assumes that the log likelihood is concave over a bounded parameter regime.  It is similar to  Assumptions V--VI in} \citet{li2023statistical} \sout{with  slightly different formation.} \blue{(YH: This discussion  does not apply. )} 

Condition \ref{cond:parfunction}   is similar to  Assumptions V and VI in \citet{li2023statistical}. First, Condition  \ref{cond:parfunction} (i) implies that $\ell(\theta,x)$  is concave and its  derivatives   remain bounded when $\theta$ belongs to a bounded regime. This can be satisfied by common exponential family distributions, such as standard normal,  Bernoulli, and Poisson distribution with canonical link functions. 
Second,  Condition \ref{cond:parfunction} (ii) implies that the edgewise score  function follows a sub-Weibull tail distribution, which generalizes the common sub-gaussian and sub-exponential tail properties \citep{vladimirova2020sub,kuchibhotla2022moving}. In addition, a common edgewise distribution family $\ell(\theta; x)$  is considered just for notational simplicity, and all the conclusions can be similarly generalized by allowing the  distributional types of $A_{ij}$ to vary across  $1\leqslant i<j\leqslant n$. 

\begin{remark}[Existence  of the constrained MLE] \label{rm:mle_exist}
The constrained MLE defined in \eqref{eq:prob} always exists under our assumptions. 
In particular, 
  the objective function $L(Z,\alpha)$ in \eqref{eq:prob} is continuously differentiable with respect to $[Z,\alpha]$ by Condition \ref{cond:parfunction} (i),  
    % , \red{(YH: With respect to what? $Y=[Z,\alpha]$?)} 
    and the feasible set in \eqref{eq:prob} is compact by its formulation.  
    % \red{(YH: Compact is from constraints, not  Condition \ref{cond:parfunction}?)} 
    Hence,  by the Weierstrass extreme value theorem, 
$L(Z,\alpha)$ attains its minimum over the compact feasible set, so the defined minimizer $\hat{Y}$ exists. 
\end{remark}

To facilitate the presentation,  we concatenate $z_i$ and $\alpha_i$ for each node $1 \leqslant i \leqslant n$ as 
\begin{align}
\label{eq:Sigma_i_of_Y}
    y_i = \begin{bmatrix}
        z_i\\
        \alpha_i
    \end{bmatrix} \in \mathbb R^{k+1}, \quad \text{and define}\quad \Sigma_i(Y) = -\sum_{j = 1}^n \ell^{\prime \prime}(\Theta_{ij}; A_{ij}) w_jw_j^\mytrans \in \mathbb R^{(k+1)\times(k+1)}, 
    % \quad  w_i =\begin{bmatrix}
    %     z_i\\
    %     1
    % \end{bmatrix}\in \mathbb R^{k+1}, 
\end{align}
  where $w_j={\partial \Theta_{ij}}/{\partial y_i} = [z_j^\mytrans, 1]^{\top}\in \mathbb{R}^{k+1}$, and $Y=[y_1,\ldots, y_n]^{\top}=[Z,\alpha]\in \mathbb{R}^{n\times (k+1)}$. 
% For each node $1 \leqslant i \leqslant n$, we define $y_i = (z_i^\mytrans, \alpha_i)^\mytrans \in \mathbb R^{k+1}$, $w_i = (z_i^\mytrans, 1)^\mytrans \in \mathbb R^{k+1}$, and $\Sigma_i = -\sum_{j \neq i} \ell^{\prime \prime}(\Theta_{ij}; A_{ij}) w_jw_j^\mytrans \in \mathbb R^{(k+1)\times(k+1)}$. 
Given a   set of $m$ node indices $\mathcal{I} = (i_1, \ldots ,i_m)$, we denote the concatenated vector
\begin{align} \label{eq:blockvecmatdef}
    y_{\mathcal{I}} =\begin{bmatrix}
        y_{i_1}\\
        \vdots\\
        y_{i_m}
    \end{bmatrix} \in \mathbb R^{m(k+1)}\quad\text{ and }\quad \Sigma_{\mathcal{I}}(Y) = \begin{bmatrix}
        \Sigma_{i_1} & \cdots & 0\\
        \vdots & \ddots & \vdots\\
        0 & \cdots & \Sigma_{i_m} 
    \end{bmatrix}\in \mathbb R^{m(k+1) \times m(k+1)}
\end{align}
being a block-diagonal matrix. 
% $y_{\mathcal{I}} = [y_{i_1}^{\top}, \ldots, y_{i_m}^{\top}]^{\top} \in \mathbb R^{m(k+1)}$ and define the block-diagonal matrix $\Sigma_{\mathcal{I}}(Y) = \operatorname{blkdiag}(\Sigma_{i_1}, \ldots, \Sigma_{i_m}) \in \mathbb R^{m(k+1) \times m(k+1)}$.
For the maximum likelihood estimator, we let $\hat{Y}_q=[\hat{Z}_q, \hat{\alpha}]$ with $\hat{Z}_q$ in \eqref{eq:rotation_z_def} and define $\hat{y}_{q,\mathcal{I}}$ and  $\Sigma_{\mathcal{I}}(\hat{Y}_q)$ similarly to \eqref{eq:blockvecmatdef} with $Y$ replaced by $\hat{Y}_q$. 
We next state the main asymptotic result.

\begin{theorem} \label{thm:homoMLE}
    Assume Conditions \ref{cond:truevalue}--\ref{cond:parfunction}, and the constant $M$ in \eqref{eq:prob} satisfies %that  there exists a constant $c>0$ such that $ M-M_1 > c $  
\begin{align}\label{eq:Mrestrction}
    M \geqslant 2M_1\  \text{ with  }M_1\text{ in Condition~\ref{cond:truevalue}.}
\end{align} 
% for a constant $c>0$ and 
% $M_1$ in Condition~\ref{cond:truevalue}. \red{(YH: TBA: move this above)}
    % \gray{the true parameters $(Z^{\star},\alpha^{\star})$ belong to the feasible regime.} \red{(need $M - M_1 \geqslant c$)} \red{@Yuang (YH: This needs to be  finalized or $M\geqslant 2M_1$ in supp?)} 
\begin{enumerate}
    \item[(i)] For any constant $\varepsilon>0$, there exist constants $C_\varepsilon > 0$ and  $N_{\varepsilon}\in \mathbb{N}$ such that when $n\geqslant N_{\varepsilon}$, with probability $1 - O(n^{-\varepsilon})$, $\hat{\alpha}$ is unique and $\hat{Z}$ is unique up to the orthogonal group $\mathcal{O}(k)$, and
    \begin{align*}
        \big\|\hat{Y}_q - Y^\star\big\|_{2 \to \infty} \leqslant C_\varepsilon n^{-1/2} \log(n).
    \end{align*}
   
    \item[(ii)]
% In addition, 
For any fixed index set $\mathcal{I}$ of $m$ nodes, as $n$ goes to infinity, %the sequence 
    \begin{align*}
        \big\{{\Sigma}_{\mathcal{I}}(\hat{Y}_q)\big\}^{1/2}\, \big(\hat{y}_{q,\mathcal{I}} - y_{\mathcal{I}}^\star\big) \stackrel{d}{\to} \mathcal N\big(0,\mathrm{I}_{m(k+1)}\big).
    \end{align*}
\end{enumerate}
% \blue{(YH: Covariance also depends on other parameters?
%     \begin{align*}
%         \big\{{\Sigma}_{\mathcal{I}}(\hat{Y}_q)\big\}^{1/2}\big(\hat{y}_{q,\mathcal{I}} - y_{\mathcal{I}}^\star\big) \stackrel{d}{\to} \mathcal N\big(0,\mathrm{I}_{|\mathcal{I}|}\big).
%     \end{align*} 
% in the subscript, $\mathcal{I}$ better to be added after $q$?)} 
% \blue{\begin{align*}
%    \text{update to } \hat{Y}_{q,\mathcal{I}}
% \end{align*} }
\end{theorem}

Theorem~\ref{thm:homoMLE} first establishes the uniqueness of    the constrained maximum likelihood estimator $\hat{Y}=[\hat{Z},\hat{\alpha}]$ up to $\mathcal{O}(k)$ transformations. Although $\hat{Z}$ is not strictly unique as a matrix, it
can be naturally interpreted as an equivalent class up to $\mathcal{O}(k)$. The aligned estimator $\hat{Y}_q=[\hat{Z}_q,\hat{\alpha}]$ serves as a convenient canonical representative for inference. 
The  unknown orthogonal matrix $\hat{Q}$ merely determines the coordinate system for representing the latent space and does not hinder   statistical interpretation for the intrinsic discrepancy. 
% This does not hinder the statistical interpretation  since any estimator $\hat{Z}$ can only be interpreted up to $\mathcal{O}(k)$ without additional assumptions. 
% we can view  it as  a unique equivalence class or  any representative element in that unique equivalence class throughout this paper. 
% Without additional assumptions,  $\hat{Z}$ can only be interpreted up to $\mathcal{O}(k)$ transformations and thus we can view any estimator $\hat{Z}$.  
Second, the two-to-infinity error bound suggests  the maximum likelihood estimator achieves uniform consistency across all parameters as $n\to \infty$.  
% latent vectors and baseline parameters.  
% \red{(YH: 2toinfy discussion TBA.)}
Third, Theorem~\ref{thm:homoMLE} further establishes entrywise asymptotic normality of the aligned estimator $\hat{Y}_q$. 
% the asymptotic normality in 
% Theorem \ref{thm:homoMLE} also holds with $\Sigma_{\mathcal{I}}(\hat{Y}_q)$ replaced by $\Sigma_{\mathcal{I}}(Y^{\star})$. Our  proof actually first shows that the asymptotic inverse covariance matrix 
% % \red{(inverse?)} 
% of $\hat{y}_{q,\mathcal{I}}-y^{\star}_{\mathcal{I}}$ is  $\Sigma_{\mathcal{I}}(Y^{\star})$ and then argues it can be consistently estimated with $\Sigma_{\mathcal{I}}(\hat{Y}_q)$. 
When $p( \cdot \mid \theta)$ follows a natural exponential family distribution with canonical link function, we actually  have that 
$\Sigma_i$ equals the Fisher information matrix of $y_i=[z_i^{\mytrans},\alpha_i]^{\mytrans}$, that is,  $\Sigma_i=\EXPT\big(\frac{\partial L}{\partial y_i}\frac{\partial L}{\partial y_i}^\mytrans\big)$.  

Although the conclusions in  Theorem~\ref{thm:homoMLE}  resemble  that in \cite{li2023statistical}, the scope here is broader without enforcing $\hat{Z}^{\top}\hat{Z}$ to be diagonal or   $Z^{\star\top}Z^{\star}$ to have distinct limiting eigenvalues. Therefore, Theorem~\ref{thm:homoMLE} provides a unified inferential result showing that the same   asymptotic structure remains valid regardless of eigenvalue  multiplicity.  This broadens the applicability of the inferential theory and paves the way for connecting  to practical  algorithms under more flexible and realistic constraints, which will    be discussed in Section~\ref{sec:gd_theory_1}. 

% motivated by Section~\ref{sec:implicit_regular}. 
% under the more flexible parameter regime motivated by Section~\ref{sec:implicit_regular}. 
% \smallskip 

\begin{remark}[Invariance with respect to $M$]\label{rm:m_value_vary}
Interestingly, our proof reveals that the constrained MLE defined  in \eqref{eq:prob} is invariant to the choice of $M$, provided $M$ is  sufficiently large as required in  \eqref{eq:Mrestrction}. 
In other words, any two values of $M$ satisfying  \eqref{eq:Mrestrction} yield identical solutions in  \eqref{eq:prob} with high probability. 
% and Theorem \ref{thm:homoPGD}  \red{(TBA: Theorem 1 also holds.)}
This is formally justified in Remark D.2 of the Supplementary Material. 
Therefore, $\hat{Y}$ is well-defined without requiring a uniquely specified $M$.  
This invariance is also reflected in the 
 practical algorithm in Section \ref{sec:gd_theory_1}, where  the empirical estimator approximates $\hat{Y}$ in \eqref{eq:prob} without explicitly specifying  $M$. 
\end{remark}

The asymptotic results in  Theorem \ref{thm:homoMLE}  provide a foundation for various downstream inferential tasks, where the target of interest is a  transformation  of $Y^{\star}=[Z^{\star},\alpha^{\star}]$. For example, for any pair of nodes $(i,j)$, their corresponding edgewise mean  is 
\begin{align}\label{eq:mulink}
    \EXPT(A_{ij}\mid \Theta^{\star})=\mu(\Theta_{ij}^{\star}) = \mu\big(\langle z_i^{\star}, z_j^{\star}\rangle + \alpha_i^\star + \alpha_j^\star \big), 
\end{align}
where $\mu(\cdot)$ denotes the link function between the expectation and parameters under the distribution $p(\cdot \mid \theta)$ in \eqref{eq:model}. More generally, let $g(\cdot)$ be a fixed function of $y_{\mathcal{I}}\in \mathbb{R}^{m(k+1)}$, where $\mathcal{I}$ is a fixed index set as in Theorem~\ref{thm:homoMLE}. 
% we consider a generic function $g(\cdot )$ that is a function of $y_{\mathcal{I}}\in \mathbb{R}^{m(k+1)}$ for a fixed index set $\mathcal{I}$ as in Theorem~\ref{thm:homoMLE}. 
By Theorem~\ref{thm:homoMLE}, the induced maximum likelihood estimator of $g(y^{\star}_{\mathcal{I}})$ is $g(\hat{y}_{q,\mathcal{I}})$, and we next establish its asymptotic distribution. %similarly to   and is presented below.  
\begin{corollary} \label{cor:gyqIhat_AN}
For a given fixed index set $\mathcal{I}$  of $m$ nodes, define $y_{\mathcal{I}}^{\star}$ and  $\hat{y}_{q,\mathcal{I}}$  as in Theorem~\ref{thm:homoMLE}. 
Let $g:\mathbb{R}^{m(k+1)}\to \mathbb{R}$ be a fixed   function that is twice continuously differentiable, and denote    its gradient by  $\nabla g(\cdot)$.   
Assume the conditions of Theorem \ref{thm:homoMLE}, and the regularity condition on $g(\cdot)$ stated in Condition E.1 in the Supplementary Material. Then, as $n\to \infty$, $ \{ g(\hat{y}_{q,\mathcal{I}}) - g(y_{\mathcal{I}}^\star) \}/ \widehat{se}(\hat{g}_{\mathcal{I}}) \stackrel{d}{\to} \mathcal N(0,1),$ where we denote  
\begin{align*}
    % \big\{ g(\hat{y}_{q,\mathcal{I}}) - g(y_{\mathcal{I}}^\star) \big\}\big/ \hat{se}(\hat{y}_{q,\mathcal{I}}) \stackrel{d}{\to} \mathcal N(0,1),\\
\widehat{se}(\hat{g}_{\mathcal{I}}):= \big\{ \nabla g(\hat{y}_{q,\mathcal{I}})^\top \Sigma_{\mathcal{I}}(\hat{Y}_q)^{-1} \nabla g(\hat{y}_{q,\mathcal{I}}) \big\}^{1/2}. 
\end{align*} 
\end{corollary}

% Corollary~\ref{cor:gyqIhat_AN} enables inference for various transformed statistics. As an example, f
As an illustration, consider inference for the edgewise mean in \eqref{eq:mulink} by Corollary~\ref{cor:gyqIhat_AN}. 
% we can apply Corollary~\ref{cor:gyqIhat_AN}  to 
  % the targeted mean in \eqref{eq:mulink}. 
  Construct $\hat{\Theta}_{ij}:=\langle \hat{z}_{q,i}, \hat{z}_{q,j}\rangle+\hat{\alpha}_i+\hat{\alpha}_j= \langle \hat{z}_{i}, \hat{z}_{j}\rangle+\hat{\alpha}_i+\hat{\alpha}_j$. 
By Corollary~\ref{cor:gyqIhat_AN} and the chain rule, 
 we have $\{\mu(\hat\Theta_{ij}) - \mu(\Theta_{ij}^\star)\}/\widehat{se}(\hat{\mu}_{ij})\xrightarrow{d} \mathcal N(0,1)$ under suitable conditions,  where  
\begin{align}\label{eq:se_mu}
    \widehat{se}(\hat{\mu}_{ij}):=\big|\mu^{\prime}(\hat\Theta_{ij})\big|\ \left\{ \begin{bmatrix}
    \hat{w}_j\\
    \hat{w}_i
\end{bmatrix}^{\top} 
\begin{bmatrix}
   \big[ \Sigma_{i}(\hat{Y}) \big]^{-1}& 0\\
    0 &    \big[\Sigma_{j}(\hat{Y})\big]^{-1}
\end{bmatrix} 
\begin{bmatrix}
    \hat{w}_j\\
    \hat{w}_i
\end{bmatrix}\right\}^{1/2},
\end{align} 
and $\mu'(\cdot)$ denotes the derivative of $\mu(\cdot)$.

\section{Adaptive Algorithms and Theoretical  Guarantees} \label{sec:gd_theory_1}

In this section, we   develop a fully data-adaptive computational framework whose output converges to the maximum likelihood estimator $[\hat{Z},\hat{\alpha}]$ with high probability. 
In particular, Sections~\ref{sec:pgd_adaptive}--\ref{sec:add_pgd} study the  adaptive projected gradient descent, covering the algorithm, its convergence theory, and additional implementation details, respectively. Sections~\ref{sec:initial1}--\ref{sec:initial_theory} then develop a range-adaptive singular value thresholding for initialization along with theoretical guarantees.

\subsection{Projected Gradient Descent Algorithm with Adaptive Learning Rate} \label{sec:pgd_adaptive}

As reviewed in Section~\ref{sec:gap_homo}, solving the MLE \eqref{eq:prob} is typically achieved by the projected gradient descent algorithm \citep{ma2020universal}. 
However, bridging the algorithm output  in practice and the idealized MLE \eqref{eq:prob} is faced with several challenges. 
First, the objective function is non-convex and the dimension of parameters is high. 
% contains high-dimensional parameters.
Second,   Theorem~\ref{thm:homoMLE}  
suggests that the solution in \eqref{eq:prob} is relevant to inferring true parameters $Y^{\star}$ if $M\geqslant 2M_1$, but in practice $M_1$ is unknown. Therefore, it can be unclear how $M$ should be chosen. 
% which can be mysterious  when $M_1$ is unknown in practice. 
Third, existing   analyses in \cite{ma2020universal}   require certain tuning parameters to be   chosen in a way relying on unknown model truth, such as  the learning rate and projection set, and therefore become impractical.  

To address all these challenges, 
we 
devise a novel backtracking line search scheme that produces data-adaptive learning rates and develop new proof techniques showing that the explicit projection onto an unknown set can be avoided. 
% allows skipping the redundant projection step. 
% The proposed algorithm with adaptive learning rates is presented in Algorithm \ref{algor:pgd}, and we next explain the details. 
% We introduce our proposed projected gradient  algorithm with adaptive learning rates, with 
A pseudo-code summary of the proposed algorithm is given in  Algorithm \ref{algor:pgd}, and we now explain the details.  
% Each iteration of gradient descent in 
Overall, Algorithm~\ref{algor:pgd} takes iterates over projected gradient descent, and each iteration consists  of three major stages. 
The first stage computes the   search direction
\begin{align} \label{eq:direction_def}
    d(Y) = -[J_n\nabla_{Z} L(Y), \nabla_{\alpha} L(Y)],
\end{align}
where $J_n=\mathrm{I}_n-1_n1_n^{\top}/n$, and $\nabla_Z L(Y)$ and $\nabla_{\alpha} L(Y)$ represent partial derivatives of $L(Y) = L(Z,\alpha)$ with respect to $Z$ and $\alpha$, respectively. Note that $d(Y)$ differs from the gradient descent direction $-\nabla_Y L(Y) $ only by a projection matrix $J_n$ multiplied to $\nabla_Z L(Y)$. As $1_n^{\top}J_n=0$, multiplying $J_n$ ensures that the updates to $Z$ part are centered, so that updated latent vectors from each iteration satisfy the identifiability constraint $1_n^{\top}Z=0$. 
% \red{(Note: d already has negative sign!)} 
The second stage runs a line search  to    adaptively determine an appropriate learning rate $\eta$.
% for algorithm convergence. 
Due to the non-convexity and high dimensionality of the  optimization problem in \eqref{eq:prob}, classical  line search methods and analyses used in convex optimization cannot directly yield convergence guarantee. 
To address that, 
we introduce a novel set of search rules. Specifically,
given  a current iterate $Y$ (for example,  $Y=Y^r$ for $Y^r$ estimate in the $r$-th iteration),  a candidate step size $\eta$ is accepted  if the following conditions are simultaneously satisfied 
\begin{align}
   L(Y + \eta d(Y)) - L(Y) - C_{\operatorname{ls}} n\eta^2 \langle \nabla_Y L(Y) , d(Y) \rangle  \leqslant 0, \label{eq:search_rule_1}\\ 
\text{ and }\ \max_{1\leqslant i\leqslant n}\big\{L_i(Y + \eta d_i(Y)) - L_i(Y) - C_{\operatorname{ls}} n\eta^2 \langle \nabla_Y L_i(Y) , d_i(Y) \rangle  \big\}\leqslant 0, \label{eq:search_rule_2}
\end{align}
where we define $d_i(Y)=-e_ie_i^{\top}\nabla_Y L(Y)$, $e_i$ is an $n$-dimensional indicator vector with only the $i$-th position being one and zero otherwise, and $L_i= - \sum_{j\in \{ 1,\ldots, n\} \backslash \{i\}} \ell(\Theta_{ij}; A_{ij} ) $, and $C_{\operatorname{ls}}$ is a hyperparameter that can be simply set to 1 in our framework. 
% \red{(YH: I add a constant to mimic  classical statement. Is it ok?)} 
Similarly to classical backtracking line search  \citep{nocedal2006numerical}, we start with an initial value of $\eta$  and  iteratively shrink it by a contraction  factor $\beta\in(0,1)$ until both  \eqref{eq:search_rule_1} and \eqref{eq:search_rule_2} are  satisfied. In   practice, one may also set a  limit for the number of shrinking iterations to prevent pathological cases. 
% in practice. 
Finally, the third stage updates the parameters along the descent direction with the selected learning rate.  
% \red{(YH: Need to add $L_i$ definition)}
% \red{(to check: is $e_i$ needed earlier?)}

% , where we define 
% In Algorithm \ref{algor:pgd}, \blue{explain notation...}  

\begin{algorithm}[!htbp]
    \SetKwBlock{IndentBlock}{}{}  
    \caption{Projected Gradient Descent with Adaptive Line Search.}
	\label{algor:pgd}
    %\setstretch{1}
    \normalsize 
	\KwIn{Data: $ A \in \mathbb{R}^{n\times n}$. \hspace{0.40em} Initial estimate: $Y^0 \in \mathbb R^{n \times (k+1)}$.\hspace{0.4em}  Initial step size: $\eta_{\operatorname{init}} > 0$. \\
      \hspace{2.65em} Backtracking rate and steps: $\beta \in (0,1)$ and $R'\in \mathbb{N}$. \hspace{0.4em} Number of iterations: $R$.}

\SetArgSty{textnormal}

\BlankLine 

\For{$r=0,\ldots, R-1$}{

\smallskip 

% $\triangleright$ {Compute descent direction}. 
\textbf{I. Compute descent direction.} 
\quad  $d(Y^r) = -[J_n \nabla_Z L(Y^r), \nabla_\alpha L(Y^r)]$.

% \medskip   
\vspace{5pt}

% $\triangleright$ Line search.  \SetAlgoVlined  %\SetAlgoNoLine 
\textbf{II. Line search.} 
\IndentBlock{ 
Initialize $\eta = \eta_{\operatorname{init}}$\; 

% \SetAlgoLined 

\If{\eqref{eq:search_rule_1} or \eqref{eq:search_rule_2} is violated  with    $Y=Y^r$ and the current $\eta$
}{ \SetAlgoNoLine 
\Repeat{\eqref{eq:search_rule_1} and  \eqref{eq:search_rule_2}   both hold, or the number of backtracking exceeds $R'$}{$\eta =\beta \eta$}
}
return final step size for the $r$-th iteration: $\eta_r = \eta$. 
}

\vspace{5pt}
% \bigskip    

\textbf{III. Update.}  \quad    $Y^{r+1}=Y^r + \eta_r d(Y^r)$. 

% \smallskip 
}

    \BlankLine
    \KwOut{$Y^R.$}
\end{algorithm}

We next discuss connections and differences between the proposed search rules \eqref{eq:search_rule_1}--\eqref{eq:search_rule_2} and existing methods in the literature.
First, \eqref{eq:search_rule_1} takes a form similar to the classical Armijo condition in convex optimization; see Remark \ref{rm:classical_line_search} for a description. 
% , with a specific description provided in Remark \ref{rm:classical_line_search}.  
However, \eqref{eq:search_rule_1} specifies a factor $n\eta^2$ before the cross product term $ \langle \nabla_Y L(Y) , d(Y) \rangle $ which is quadratic in $\eta$ rather than linear as in the classical Armijo condition. 
% and differs from classical Armijo condition. 
This modification  is crucial for addressing the theoretical challenges from intrinsic non-convexity and high-dimensionality of our problem. 
% This subtle difference is proposed to tackle with theoretical challenges arising from the intrinsic non-convexity and high-dimensionality of our problem. 
% \red{(YH: Any other reasons to mention?)} 
Second, \eqref{eq:search_rule_2} 
is equivalent to imposing $n$  conditions on $L_i(Y)$ and $d_i(Y)$ with $i\in \{1,\ldots,n\}$  simultaneously.  
For each $i=1,\ldots, n$, the inequality constraint takes a form similar to \eqref{eq:search_rule_1}. 
But $Y+\eta d_i(Y)$ differs from $Y$ only by updating its $i$-th row $y_i$ along the  direction $-\nabla_{y_i}L(Y)$.  Intuitively, \eqref{eq:search_rule_2} evaluates the ``goodness'' of update along each row of $Y$, whereas \eqref{eq:search_rule_1} evaluates the update of the whole matrix $Y$. Technically, imposing \eqref{eq:search_rule_2}  ensures  row-wise  control needed in the convergence analysis. 
% in derivation. 
This helps overcome the necessity of  projecting onto an unknown set in the  original Algorithm 1 in \cite{ma2020universal}.  
Notably, both \eqref{eq:search_rule_1} and \eqref{eq:search_rule_2}  are  data-adaptive and practically flexible.

% it is imposed to ensure the node-wise convergence property. 
% \textit{interpretation for \eqref{eq:search_rule_2}}:  \eqref{eq:search_rule_2} can be viewed as row-wise update rule. this can help relax projection onto bounded set 

\begin{remark}[Classical Armijo condition]\label{rm:classical_line_search}
% As an example, 
% For clarity and ease of comparison, 
% As an illustration,
For comparison, 
we give a simple example 
of  minimizing a convex function $f(x)$ with $x\in \mathbb{R}^k$ by the line search method under the classical Armijo condition \citep{nocedal2006numerical}. At the $r$-th iteration, let $x_r$ denote the estimate for the current iterate and let $d_r$ be a descent direction satisfying $\langle \nabla f(x_r), d_r\rangle <0$. 
% suppose we have obtained parameter estimate $x_r$ and a search direction $d_r$ 
% To construct an updated estimate $x_r+\eta_r  d_r$, 
% To determine a step size $\eta$ to construct new estimate $x_r+\eta d_r$, 
% At a current estimate, say $x_r$ in the $r$-th iteration, we find a descent direction $d_r$ satisfying $\langle \nabla f(x_r), d_r\rangle <0$. 
By iteratively shrinking from an initially large $\eta$, classical Armijo condition accepts the first value of $\eta$ that satisfies 
% For a current parameter $Y$, update direction $d(Y)$ satisfying , and a learning rate $\eta$, classical Armijo condition often specifies 
$  f(x_r+\eta d_r) - f(x_r) \leqslant c \eta  \langle  \nabla f(x_r), d_r\rangle $ for $c\in (0,1)$, and then updates the iterate to $x_r+\eta   d_r$. 
% ,  where $c\in (0,1)$ is a hyperparameter. 
% Then a final updated estimate $x_r+\eta   d_r$ is constructed. 
% The classical formula takes a form  similar to our proposed
Our proposed conditions \eqref{eq:search_rule_1}--\eqref{eq:search_rule_2} have a similar form but are fundamentally different by replacing    $c\eta$
% The classical Armijo rule is linear in $\eta$ and with 
% but  $c\eta$ is replaced 
with $C_{\mathrm{ls}}n\eta^2$. 
This distinction is a critical innovation to accommodate the high-dimensional network model, where the number of latent vectors grows with $n$.
% the high-dimension of parameters (increasing with $n$) under the network model. 
Accordingly, unlike $c\in (0,1)$ required in the Armijo rule, the constant $C_{\mathrm{ls}}$ can be chosen more flexibly and can be equal to or even greater than 1. 
% set to 1.  
% The constant $C_{\mathrm{ls}}$ in \eqref{eq:search_rule_1}--\eqref{eq:search_rule_2} can take flexible values, such as 1 or larger, which differs from $c\in (0,1)$ in Armijo condition. 
% classically suggest using a small $c$ like $10^{-4}$ but 
%typically 
% and we do not require $C_{\mathrm{ls}}$. This means that how small $c$ 
% \blue{But  our analysis does not rely on $c$ being small because... (YH: What may be added? may decide after simulations)}
\end{remark}

% First, consider update current $Y$ to $Y+\eta d(Y)$. Then  \eqref{eq:search_rule_1} takes a form similar to classical Armijo condition in convex optimization. 
% \bigskip 

\subsection{Convergence Theory for Projected Gradient Descent} \label{sec:pgdtheory}

% To obtain the convergence guarantee of Algorithm \ref{algor:pgd},  we also need appropriate conditions imposed on the initial estimator $Y^0$ and the initial step size $\eta_{\text{init}}$. 
% Technically,
This subsection presents the convergence theory for Algorithm \ref{algor:pgd}. 
Because latent vectors are identifiable only up to orthogonal transformation $\mathcal{O}(k)$, 
we measure the distance between two  matrices of latent vectors $\hat{Z}$ and $Z\in \mathbb{R}^{n\times k}$ by 
\begin{align*}
    \mathrm{dist}(\hat{Z}, Z):=\min_{Q\in \mathcal{O}(k)} \,\|\hat{Z} - Z Q\|_{\F}.  
\end{align*}
% which accounts for the identifiability of latent vectors up to the orthogonal transformation $\mathcal{O}(k)$. 
This is similarly considered in other studies \citep{ma2020universal,he2023semiparametric,tian2024efficient}, and it can be typically shown to be of the same order as 
% which argues it can be equivalent to  
$\|\hat{Z}\hat{Z}^{\top}-ZZ^{\top}\|_{\F}/\sqrt{n}$ \citep{tu2016low}. 
% up a multiplicative factor of  order $1/\sqrt{n}$. 
% due to the unidentifiability of latent vectors. 
Then for $\hat{Y}=[\hat{Z}, \hat{\alpha}]$ and $Y = [Z, \alpha]$, we define the squared overall error as $  \mathrm{dist}^2(\hat{Y},Y ) :=   \mathrm{dist}^2(\hat{Z}, Z) + \|\hat{\alpha}-{\alpha}\|_2^2$, where   the discrepancy between $\hat{\alpha}$ and ${\alpha}$ is directly measured by the vector $\ell_2$ norm. 

To establish the convergence, we impose a regularity condition on the   inputs of Algorithm \ref{algor:pgd} and discuss its  implications afterwards.

\vspace{-0.7em}
\begin{condition}[Inputs of Algorithm~\ref{algor:pgd}]\label{cond:initial1}
Assume: 
\begin{enumerate}
    \item[(i)] \ \ 
    The initial estimate $Y^0 = [Z^0, \alpha^0]$ satisfies (a)  $1_n^\mytrans Z^0 = 0$,  (b) $
  \operatorname{dist}^2(Y^0, Y^\star) \lesssim n^{1-\varsigma_0}$ for a constant  $\varsigma_0 \in (0,1/2)$, and (c)  $ \|Y^0\|_{2\to\infty} \leqslant b_0 $ for a constant $b_0>0$.      
    \item[(ii)]\ The initial step size $\eta_{\operatorname{init}}$ satisfies $1/C_{\operatorname{init}}\leqslant n\eta_{\operatorname{init}}\leqslant C_{\operatorname{init}}$   
    % $1/(C_{\operatorname{init}}n)\leqslant \eta_{\operatorname{init}}\leqslant C_{\operatorname{init}}/n$ 
    for a constant $C_{\operatorname{init}}>1$. 
    \item[(iii)] The backtracking iteration limit $R'$ satisfies $R'\to \infty$ as $n\to \infty$. % $1 \ll R' \ll \log(n)$.
\end{enumerate}
\end{condition}

% Condition \ref{cond:initial1} (i) imposes three assumptions on the initial estimate $Y^0$. 
To keep the convergence theory broadly applicable, Condition~\ref{cond:initial1}~(i) gives an abstract initialization condition on $Y^0$. In   Section~\ref{sec:initial1}, we 
will provide  a concrete example based on the singular value thresholding  and verify   Condition \ref{cond:initial1} (i) holds with   high probability.
% To show the general applicability of Algorithm \ref{algor:pgd}, we impose Condition \ref{cond:initial1} (i) on the initial estimator $Y^0$ to establish the convergence guarantee in Theorem \ref{thm:homoPGD}, while we will also provide  a specific example based on the singular value thresholding  in Section \ref{sec:initial1} and show it satisfies Condition \ref{cond:initial1} (i) with a high probability. 
We next explain the implications of the three requirements (a)--(c) in Condition \ref{cond:initial1} (i). 
First, the constraint  $1_n^{\top}Z^0=0$ ensures that the columns of $Z^0$ are centered.  This requirement is mild since any estimate can be centered 
% as we can always center an estimate 
without changing the induced likelihood value as discussed after Condition \ref{cond:truevalue}. Moreover, with this initialization, the next iterate also satisfies the centering constraint as $1_n^{\top}Z^1=1_n^{\top}Z^0-\eta_0 1_n^{\top} J_n \nabla_Z L(Y^0) =0$ by $1_n^{\top}J_n=0$. Therefore,  
throughout all the iterations in Algorithm \ref{algor:pgd}, latent vectors remain centered, and the translational non-identifiability is avoided. 
% have zero mean  and there is no unidentifiability due to mean shift. 
Second, part (b) imposes an overall error bound on $Y^0$.
% the Frobenius norm of $Y^0$ after accounted for best orthogonal transformation of latent vectors. 
% \blue{TBA...} 
It is similar to Assumption 8 on the overall error of initialization in \cite{ma2020universal}.
In our problem setting, the required rate in \cite{ma2020universal} is slightly weaker as it reduces to $  \operatorname{dist}^2(Y^0, Y^\star)\leqslant cn$ for a sufficiently small $c$.  
% which, though, only requires $  \operatorname{dist}^2(Y^0, Y^\star)\leqslant cn$ for a sufficiently small $c$ and thus  is slightly weaker than Condition \ref{cond:initial1} (i).  
% However, we point out  
Nevertheless,  the sharper rate in Condition \ref{cond:initial1} (i) does not lead to a significant difference and can be shown to be satisfied by the initial estimator presented in Section \ref{sec:initial1}. 
% which has been shown to be satisfied by \blue{...} 
Third,  part (c) implies that row-wise $\ell_2$ norms of $Y^0$ are uniformly bounded.
Notably,  $b_0$ only needs to be a universal constant and does not rely on any  unknown true model parameters. 
This differs from \cite{ma2020universal}  that enforces boundedness by explicitly projecting each iterate onto a bounded set relying on the unknown truth. 

Condition \ref{cond:initial1} (ii) 
implies that  
% required rate for 
the initial learning rate  should be  of the order of $1/n$. 
% specifies the required rate for the initial learning rate  in the gradient descent. 
% This order is consistent with prior requirement in \cite{ma2020universal} under our considered settings. 
It is a weak requirement because no precise constant needs to be specified, thanks to 
% our construction of 
the adaptive search rules \eqref{eq:search_rule_1}--\eqref{eq:search_rule_2}.  
% The  final outcome of Algorithm \ref{algor:pgd} would not be sensitive to a particular value of  $C_{\operatorname{init}}$. 
This overcomes the implicit dependence on the unknown ground truth of step sizes in previous analysis; more technical details are discussed in Remark \ref{rm:ma_eta_restricction}. 

% Condition \ref{cond:initial1} (iii) gives  a budget on the number of backtracking steps 
% In the line search,
In Condition~\ref{cond:initial1}~(iii), $R'$   represents  the budget on the number of backtracking steps and is required only to diverge with $n$.   
% limit of iterations employed in the line search stage.  
% Setting a slow rate such as 
This is simply to rule out pathological scenarios where the backtracking procedure would continue indefinitely. 
% line search without a cap may run endlessly  in practice. 
% ensures that the line search will not run endlessly in practice. 
% To establish Theorem \ref{thm:homoPGD}, 
% % we require
% % it suffices that $R'$  diverges with $n$, 
% % increases as $n$ increases, 
% Condition~\ref{cond:initial1}~(iii) requires $R'$ to diverge with $n$, 
% This is a mild requirement as 
%   the divergence rate can be slow, 
% does not have to be fast, 
Setting a slow rate, 
such as $R' \asymp \log n$, would be practically innocuous. 
% Actually, Corollary \ref{cor:line_iteration} below 
Actually, Section \ref{sec:add_pgd} will show that the 
backtracking can end in finite steps with   high probability.

We are now ready to present the convergence guarantee for the output of Algorithm~\ref{algor:pgd}. 

% \vspace{-0.8em} 

\begin{theorem} \label{thm:homoPGD}
    Assume Conditions \ref{cond:truevalue}--\ref{cond:initial1}. For any constant $\varepsilon > 0$, there exist constants $C > 0$, $c_0 \in (0,1)$, and $ N_\varepsilon \in \mathbb N$ such that when $n \geqslant N_\varepsilon$,
    \begin{align*}
         \mathrm{dist}^2(Y^r, \hat{Y}) \leqslant C (1- c_0)^r   \mathrm{dist}^2(Y^0, \hat{Y})
    \end{align*}
    holds for any $0 \leqslant r \leqslant R$ with probability $1 - O(n^{-\varepsilon})$, where $\hat{Y}=[\hat{Z},\hat{\alpha}]$ denotes any  representative element in the equivalence class of the constrained maximum likelihood estimator defined in \eqref{eq:prob}. 
\end{theorem}

% \blue{TBD: $c$ depends on $\beta$ and $b_1$; specific form given in the proof \red{TBA}} 

% Theorem \ref{thm:homoPGD} is similar to Theorem \blue{TBA...} in \cite{ma2020universal} 
% when $\kappa_1(\delta_n)$ is a constant. \blue{(YH: can we give more explicit relationship? like $\kappa$ reduces to $e^{M_1}\kappa$ in \cite{ma2020universal})}  
% \blue{(Initialization: can we change $\eta_{\text{initial}}$ to a parameter that can be user-specified? Like $Y^0$, we give some conditions. We can give an example that $\eta_{\text{initial}}=1/(4\|\Sigma_L(Y^0)\|_{\text{op}})$ where $\Sigma_L $  is ...? We may also use $Y^r$ after updating in each iteration.)}

Theorem \ref{thm:homoPGD} establishes an R-linear convergence guarantee for  Algorithm \ref{algor:pgd}. 
It implies that the distance between $Y^r$ and  $\hat{Y}$ is bounded by $(1-c_0)^r$  which converges to 0 as $r\to \infty$. Therefore, as the number of iterations $R$ increases, the output of Algorithm \ref{algor:pgd} converges to the constrained maximum likelihood estimator $\hat{Y}$, up to the equivalence class induced by $\mathcal{O}(k)$.    
Specifically, when $R\gg \log n$,   the algorithmic  error 
% degenerates to zero  and 
is negligible compared to the statistical error in Theorem \ref{thm:homoMLE}. 
% Therefore, we expect 
This suggests that the asymptotic distribution in Theorem \ref{thm:homoMLE} can also hold for the algorithm output $Y^R$ approximately, thereby justifying applying inferential procedures to $Y^R$ in practice. 
% Theorems  \ref{thm:homoMLE} and \ref{thm:homoPGD} separate the statistical error and algorithm error.   
% Theorem \ref{thm:homoPGD} examines the distance between algorithm output $Y^R$ and $\hat{Y}$, and  and thus differs from 
% As a comparison, 
This conclusion is substantially different from the convergence guarantee  
(Theorem 9) in \cite{ma2020universal} which essentially controls the distance between algorithm output $Y^R$ and true model parameters $Y^{\star}$.  
% Due to different targets, 
Because their target is different, 
the corresponding error bound in \cite{ma2020universal} mixes the statistical error with the algorithmic convergence error, 
and therefore does not 
% which cannot 
directly justify the use of the inferential results in Theorem \ref{thm:homoMLE}.  
% Moreover, the F-norm-based distance $\mathrm{dist}^2(Y^R,Y^{\star})\lesssim \mathrm{dist}^2(Y^R,\hat{Y}) + \mathrm{dist}^2(\hat{Y},Y^{\star}) = O(n...)$, \blue{(may add 2toinfty norm bound above to discuss here)} which   are consistent with the rates in \cite{ma2020universal} 
% (which didn't separate two errors explicitly and cannot justify the use of inference). 

Interestingly, 
% Algorithm~\ref{algor:pgd} and 
Theorem~\ref{thm:homoPGD} does not explicitly depend on  the constrained constant $M$ in \eqref{eq:prob}.
This is non-trivial because Theorem~\ref{thm:homoMLE}  requires $M$ to be sufficiently large to study $\hat{Y}$, and we remove the explicit projection onto an unknown bounded set compared to Algorithm~1 in \cite{ma2020universal}.  
This is made possible by two key theoretical results.   
First, as discussed in Remark \ref{rm:m_value_vary}, 
we show  the definition of $\hat{Y}$ is not sensitive to a particular choice of $M$, provided $M$ is large enough. 
Second, we develop new theoretical techniques showing  that   $\|Y^r\|_{2\to \infty}$ remains  in a uniformly bounded region throughout the iterations with high probability.
% , whenever the initial estimator $Y^0$ satisfies   Condition~\ref{cond:initial1}. 
% This bounded-iterate argument allows us to establish convergence to $\hat{Y}$ without carrying $M$ into either the algorithm or the statement of Theorem~\ref{thm:homoPGD}. 

The boundedness phenomenon is related to the   implicit regularization noted by  \citep{ma2019implicit} under  different  models, but the analysis here is  fundamentally different. 
% yet requires a fundamentally different analysis here. 
% Related phenomena have been noted as implicit regularization   \citep{ma2019implicit} under  different  models.  
Specifically,  our framework must tackle  the intrinsic non-linearity with respect to $Y$, e.g., the link function in \eqref{eq:mulink} is non-linear. 
Also, it needs 
  to accommodate our fully data-adaptive scheme that does not depend on unknown true parameters. These challenges cannot be directly handled by existing techniques. As a by-product, our analysis bypasses the construction of the leave-one-out sequences in \cite{ma2019implicit}, which could be of independent theoretical interest. 
These  efforts not only make technical contributions but also further support the practical relevance of Algorithm~\ref{algor:pgd} under more flexible settings.   
% Nevertheless, our analysis is substantially different because it needs  to address the non-linearity of the link function in \eqref{eq:mulink},  
% based on a substantially different framework that 
% bypasses the construction of 
% % We   do  not need to construct 
% the leave-one-out sequences in \cite{ma2019implicit}, and supports adaptive learning rates that do not depend on unknown true parameters.  
% This property of independence with unknown parameters not only highlights the theoretical novelty of our result, but  also    demonstrates the flexibility of Algorithm \ref{algor:pgd} and supports its practical use. 

% and validity of using Algorithm \ref{algor:pgd} in practice. 

% we  can consider a relatively large value of $M$
% such that  both the true parameters $Y^{\star}$ and the initial estimator $Y^0$ specified by Condition \ref{cond:initial1} satisfy the constraints in 

% show that as long as the initialization starts in a ``good''  region  defined in Condition \ref{cond:initial1}, the iterative  updates remain in a bounded region and equivalent to MLE (not sensitive to bounded). 

% \red{(YH: $c$ may be too universal. Maybe change a notation. The proof can still use $c$ primarily. Just add one more sentence to define a specified notation for Theorem 2 inequality. Check Remark~\ref{rm:gd_rate}. @Yuang.)}

The constant $c_0 \in (0,1)$ in Theorem \ref{thm:homoPGD} quantifies the contraction rate of the error upper bound. 
% specifies how fast the error bound shrinks. 
% with a specific form given in \eqref{eq:fnorm_c_form} of the Supplementary Material. 
Its explicit expression is technically involved, as it reflects the adaptiveness of our algorithm and generality of our distributional framework, so we defer it to the Supplementary Material; see  Remark F.3.  
% Due to the adaptiveness of our algorithm and universality of our distributional framework, it involves a lot of  technical complications to explicitly describe $c_0$, which   therefore, is only presented in the Supplementary Material.  
% Note Theorem \ref{thm:homoPGD} only gives an upper bound of the convergence error, it may not uniformly characterize  
% we note that the specification of $c$ 
Nevertheless, we mention that the derived formula shows that the convergence rate depends on the same key model characteristics as in \cite{ma2020universal}. In particular, $c_0$ decreases as the condition number of $Z^{\star\top}Z^{\star}/n$ or $M_1$ in Condition~\ref{cond:truevalue} increases.
% ; see  in the Supplementary Material for more details. 
Thus, the obtained  convergence rate  is consistent with that in \cite{ma2020universal} even under the more adaptive and  general  scenarios considered in this paper.  
% The derived form reveals  how the convergence rate depends on quantities such as the condition number of $Z^{\star\top}Z^{\star}/n$ and $M_1$ in Condition~\ref{cond:truevalue}, which is consistent with those in \cite{ma2020universal}. 
% Specifically,  $c_0$ increases, that is, the error upper bound shrinks at a faster rate, if the condition number of $Z^{\star\top}Z^{\star}/n$ or $M_1$ increases. 

% The derived form reveals that $c_0$ increases, that is, the error upper bound shrinks at a faster rate, if the condition number of $Z^{\star\top}Z^{\star}/n$ or $M_1$ in Condition~\ref{cond:truevalue} grows. Details are discussed in Remark~\ref{rm:gd_rate} in the Supplementary Material. 
% These quantitative insights  are consistent with those in \cite{ma2020universal}. 

\vspace{-0.6em} 
\begin{remark}\label{rm:thm_novel}
Although the algorithmic structure of Algorithm \ref{algor:pgd} is similar to that in \cite{ma2020universal}, we emphasize the underlying  proofs of Theorem~\ref{thm:homoPGD} follow substantially different arguments. First, we directly deal with the non-convex loss over $Y=[Z,\alpha]$ and explicitly disentangle the statistical and algorithmic errors. Second, our  proof of Theorem \ref{thm:homoPGD} not only establishes the F-norm convergence but also shows uniform boundedness of $y_i^r$ across all iterations. The latter cannot be shown by existing techniques. 
% in \cite{ma2020universal}. 
% This result has not been established under related studies of latent space models to our knowledge. 
% This is achieved by showing that each $y_i$ converges to $y_i^{\star}$ for $1\leqslant i \leqslant n$ up to  a small error.
% ; see Section \ref{pf:homoPGD3} of the Supplementary Material. 
Third, all of our convergence conclusions are fully data-adaptive and do not require unknown    information about true parameters. 
Fourth, to fulfill the requirements of the initialization, we develop new data-adaptive  singular value thresholding estimator  in Section \ref{sec:initial1} below.
Taken together, these features show that our analysis is 
% not a minor modification of existing arguments, but 
a genuinely different framework with both theoretical and practical advantages.
\end{remark}

% $M_1$ or $1/M_2$  increases for $M_1$ and $M_2$ given in Condition \ref{cond:truevalue}. 
% This implies that the error bound shrinks at a faster rate if the the condition number of $Z^{\star\top}Z^{\star}/n$ 
% As the condition number of $Z^{\star\top}Z^{\star}/n$ can be bounded by $M_1^2/M_2$, we know that the error bound shrinks at a faster rate if the condition number of $Z^{\star\top}Z^{\star}/n$ is larger. 
% This implies that the shrinkage of the error bound is faster if 
% But conceptually, 
% our derivation reveals that $c$ depends on the curvature the log likelihood $\ell(\theta;x)$ and step sizes $\eta_r$ across iterations. 
% %$b_0$ increases step size $\eta$ decreases ....  
% % Although Theorem \ref{thm:homoPGD} only gives upper bound of convergence error, the 
% % we didn't establish ``if and only if'' conclusions, and convergence 

\subsection{Implementation Details of Projected Gradient Descent} \label{sec:add_pgd}
This subsection discusses several implementation details and variants of Algorithm~\ref{algor:pgd}. 
First, Corollary~\ref{cor:line_iteration} will show that setting a backtracking   budget $R'$ would not impact reaching an acceptable step size with high probability.

\begin{corollary} \label{cor:line_iteration}
The line search step in Algorithm \ref{algor:pgd} has the following  guarantees: 
\begin{enumerate}[label=(\roman*)]
    \item 
     Assume Conditions \ref{cond:truevalue}--\ref{cond:initial1}.  
    % Assume the same conditions as those of Theorem \ref{thm:homoPGD}. 
    For any $\varepsilon>0$, there exists $N_{\varepsilon}\in \mathbb{N}$, such that when $n\geqslant N_{\varepsilon}$,  the line search stage  terminates within $R'$ steps across all iterations in Algorithm \ref{algor:pgd}, with probability $1-O(n^{-\varepsilon})$. 
    \item 
         Assume Conditions \ref{cond:truevalue}--\ref{cond:parfunction} and Condition  \ref{cond:initial1} (i).  
    % Assume the same conditions as those of Theorem \ref{thm:homoPGD}, and additionally 
    Additionally, assume that the initial estimate $Y^0 = [Z^0,\alpha^0]$ satisfies  $\|Z^0- Z^\star Q_z\|_{2\to \infty} + \|\alpha^0 - \alpha^\star\|_{\infty} \lesssim \zeta_{0,n}$ for some sequence $\zeta_{0,n} \ll 1$ and some $Q_z \in \mathcal{O}(k)$. Choose the initial step size as 
        \begin{align} \label{eq:eta_init_choice}
\eta_{\operatorname{init}}= \frac{1}{6\max_{1 \leqslant i< j\leqslant n} \{1 -\ell^{\prime \prime}(\Theta_{ij}^0)\}}\cdot \min \left\{ \frac{1}{ \|Z^0\|_{\op}^2},\, \frac{1}{n}\right\}, 
    \end{align}  
where $\Theta_{ij}^0=\alpha_i^0+\alpha_j^0+\langle z_i^0,z_j^0\rangle$. 
%     \begin{align} \label{eq:eta_init_choice}
% \eta_{\operatorname{init}}= \frac{1}{6\max_{1 \leqslant i< j\leqslant n} \{1 -\ell^{\prime \prime}(\Theta_{ij}^0)\}}\cdot \frac{1}{ \|W^0\|_{\op}^2}. 
%     \end{align} 
    % where $W^0=[Z^0, 1_n]\in \mathbb{R}^{n\times (k+1)}$. 
    %which (YH: It is hard to directly see this holds for general readers. Should be claimed as a result we prove or )
Then $\eta_{\operatorname{init}}$ satisfies  Condition~\ref{cond:initial1} (ii), and for any $\varepsilon>0$,  there exists $N_{\varepsilon}\in \mathbb{N}$, such that when $n\geqslant N_{\varepsilon}$, the line search  conditions \eqref{eq:search_rule_1}--\eqref{eq:search_rule_2} are satisfied at $\eta = \eta_{\operatorname{init}}$ for all iterations  with probability $1-O(n^{-\varepsilon})$.  
% $Y^r$ when $n \geqslant N_\varepsilon$, requiring no backtracking.
% \gray{which satisfies Condition~\ref{cond:initial1} (ii) with probability $1-O(n^{-\varepsilon})$. Moreover, with probability $1-O(n^{-\varepsilon})$, the backtracking conditions \eqref{eq:search_rule_1}--\eqref{eq:search_rule_2} are satisfied at $\eta = \eta_{\operatorname{init}}$ for all iterates $Y^r$ when $n \geqslant N_\varepsilon$, requiring no backtracking.} \red{(@Yuang YH: I updated order.)}
\end{enumerate}
\end{corollary}

Case (i) in  Corollary~\ref{cor:line_iteration}  indicates  that when the backtracking budget  $R'$ diverges with $n$, the imposed cap would not interfere with reaching a desired step size with  high probability. 
 % It suggests a useful balance between  asymptotic analysis and practical implementation. 
Case (ii) considers a scenario with a stronger initialization. Specifically, 
  if a sharper two-to-infinity error bound  for $Y^0$ is available, we can construct a  closed-form choice of $\eta_{\operatorname{init}}$  in \eqref{eq:eta_init_choice}  such that the line search conditions  \eqref{eq:search_rule_1}--\eqref{eq:search_rule_2}  hold with high probability. Thus, $\eta_r=\eta_{\operatorname{init}}$ throughout, and the backtracking loop needs not to be triggered. This two-to-infinity-error assumption is much stronger than the row-wise boundedness in Condition~\ref{cond:initial1}. It  could hold under a stronger structural assumption on $\mathbb{E}(A)$, such as low rank \citep{cape2019two}. 
 Notably, $\eta_{\operatorname{init}}$ in \eqref{eq:eta_init_choice} takes a form similar to the choice in \cite{ma2020universal};  see   Remark \ref{rm:ma_eta_restricction} for more discussions. This helps explain why step sizes of a similar form have worked well in prior studies   \citep{ma2020universal,he2023semiparametric}.

\vspace{-4pt}
\begin{remark}[Comparison with step sizes  in \cite{ma2020universal}] \label{rm:ma_eta_restricction}
% We compare with the choice of step size  in \cite{ma2020universal}.  
\citet{ma2020universal} set  step sizes for $Z$ and $\alpha$ as  $\eta_Z=\eta /\|Z^0\|_{\mathrm{op}}^2$ and $ \eta_\alpha=\eta /(2 n)$, respectively, where $\eta$ is a hyperparameter.
To establish  theoretical guarantees,  $\eta$ is required to  satisfy   inequalities that involve unknown model parameters; see,  e.g., the last three inequalities in their proof of Lemma 25. %\red{(YH: Its Theorem 7 only says $\eta\leqslant c$ for a universal constant $c$???)} 
In contrast, our proposed line search rules \eqref{eq:search_rule_1} and \eqref{eq:search_rule_2} determine the step size adaptively from the data and avoid the need to specify such unknown constants.
% avoid choosing unknown constants  and therefore offers a useful  strategy that safely bridges  practice and theory.  
Moreover, for the choice in \eqref{eq:eta_init_choice}, we have $\eta_{\operatorname{init}}\asymp \min\{1/ \|Z^0\|_{\operatorname{op}}^2 , 1/n\}$,
% with   high probability \red{(YT: why WHP?)}, 
matching the choices in \cite{ma2020universal} up to a multiplicative constant. 
% linking to the original choices in \cite{ma2020universal} up to a multiplicative constant. % depending on the curvature of the log likelihood $\ell(\theta;x)$. 
\end{remark}

% \vspace{-5pt}
For clarity of exposition, we present the projected gradient descent Algorithm~\ref{algor:pgd}  in a streamlined  form above. We emphasize that our analytical framework readily accommodates  different  variations in practice and next discuss several  examples below. 

\vspace{-6pt}
\begin{remark}[Component-wise step sizes] \label{rm:comp_step_size}
The above discussions consider a common step size for updating both $Z$ and $\alpha$ in each iteration. 
But our results can be straightforwardly generalized when considering component-wise step sizes. 
In particular, we can replace $Y+\eta d(Y)$ and $Y+\eta d_i(Y)$ in \eqref{eq:search_rule_1}--\eqref{eq:search_rule_2} with  $Y+ d(Y) \eta_D$ and  $Y+ d_i(Y) \eta_D$, respectively, where $\eta_D=\mathrm{blkdiag(\eta_Z,\ldots, \eta_Z,\eta_{\alpha})}\in \mathbb{R}^{(k+1)\times (k+1)}$. Then we can conduct the line search with $\eta_Z=\beta\eta_{Z}$ and $\eta_{\alpha}=\beta\eta_{\alpha}$ similarly. 
All of our developments and conclusions extend with only routine notational modifications.
% remain similar just with extra  tedious  details. 
\end{remark}

\vspace{-1.8em} 
\begin{remark}[Iteration-varying initial step sizes] \label{rm:stepsize_vary}
The above discussions use a common initial step size $\eta_{\operatorname{init}}$ across all iterations $1\leqslant r\leqslant R$. More generally, our framework   also allows $\eta_{\operatorname{init}}$ to vary with $r$, which remains  theoretically valid and may  provide  additional adaptiveness.  
For example, 
Case (ii) in Corollary~\ref{cor:line_iteration} motivates constructing  $\eta_{\operatorname{init},r}$ in the same form as \eqref{eq:eta_init_choice}, 
% similarly to  $\eta_{\operatorname{init}}$, 
but  with $Y^0=[Z^0,\alpha^0]$  replaced by $Y^r=[Z^r,\alpha^r]$ at the $r$-th iteration. The same argument then shows that   $\eta_{\operatorname{init},r}$  satisfies Condition~\ref{cond:initial1} (i) with high probability, and the convergence guarantee continues to hold, that is, 
% Moreover, since 
% Moreover, it implies that  
% Theorems \ref{thm:homoMLE}  and \ref{thm:homoPGD}  imply that 
$Y^r$ approaches $Y^{\star}$ up to the identifiability constraint as $r$ increases. As a result, 
we expect $\eta_r=\eta_{\mathrm{init},r}$ when $r$ is large, thereby reducing  the computational cost from triggering the backtracking loop.  
\end{remark} 
% Moreover,  Theorems \ref{thm:homoMLE}  and \ref{thm:homoPGD}    imply that as $r$ increases, $Y^r$ approaches $Y^{\star}$ (up to  identifiability constraints).  
% Similarly to Case (ii) of Corollary~\ref{cor:line_iteration}, 
% we expect the line search to be skipped more frequently 
% in later iterations, improving  computational efficiency.  

\vspace{-1.8em} 
\begin{remark}[Stopping criteria]
\label{rmk:stopping}
Besides using a targeted $R$ for termination, another common practice  in gradient-descent-based algorithms is to set stopping criteria. Our theoretical analysis is flexible and can be   extended to justify that scenario too. As an example, one widely used criterion  is  to stop at the $r$-th iteration if $\|\nabla_Y L(Y^r)\|_{\mathrm{F}} \leqslant \epsilon$   for  a prespecified $\epsilon$. 
Our proof of Theorem~\ref{thm:homoPGD} shows that with high probability, the iterates  $Y^r$ in Algorithm \ref{algor:pgd} satisfy $\mathrm{dist}(Y^r, \hat{Y}) \asymp n^{-1} \|\nabla_Y L(Y^r)\|_{\mathrm{F}}$, so  $\|\nabla_Y L(Y^r)\|_{\mathrm{F}} \leqslant \epsilon$ guarantees $\mathrm{dist}(Y^r, \hat{Y}) =O_p( \epsilon/n)$. 
Hence when $\epsilon$ is chosen  sufficiently small, we expect   the difference between $Y^r$ and $\hat{Y}$   to be  ignorable for practical use. 
This provides theoretical justification for using the gradient norm as a practical stopping criterion.
\end{remark}

\subsection{Initialization: Range-Adaptive Singular Value Thresholding } \label{sec:initial1}

As mentioned in Section \ref{sec:pgdreview}, 
 the spectral-based universal singular value thresholding \citep{chatterjee2015matrix} is one of the most widely  used  approaches for initializing the projected gradient descent under latent space models  \citep{zhang2020flexible,he2023semiparametric,li2023statistical}. 
% For example, 
A representative example is 
the initialization method (Algorithm 3)  in \cite{ma2020universal}, which builds on the USVT under the logic transform. 
% developed based on  USVT  under the logic transform.  
% It first applies  the singular value thresholding to the observed matrix $A$ to obtain a low-rank approximation matrix, and then projects each element of the low-rank matrix elementwisely onto ... 
% \cite{ma2020universal} has developed  Algorithm 3  based on USVT for the latent space model with the logic transform. 
However, similarly to the challenges in the projected gradient descent procedure, 
this approach requires an explicit projection step relying on unknown model parameters to obtain the theoretical guarantee. 
% for the initial estimator.  

% It is easy to implement and can retain theoretical guarantee under a wide class of models. 
% However, to apply the original USVT algorithm to our problem and reach desired guarantee, one TBA ....
% Moreover, 
% to achieve the above convergence guarantee, we need appropriate bounds on two-to-infinity norm of the initial estimate $Y^0$ as shown in Condition \ref{cond:initial1}. 
% Such results were not directly provided in \cite{chatterjee2015matrix} or \cite{ma2020universal}. 

To obtain an initial estimator that is practically applicable while still retaining the theoretical guarantee, we propose a  Range-Adaptive Singular Value Thresholding (RA-SVT) method, with a pseudo-code summary in Algorithm \ref{algo:initialpgd}. 
% A detailed description is provided in  Algorithm \ref{algo:initialpgd}. 
Specifically, after applying the singular value thresholding to the observed matrix $A$ and obtaining  $\tilde{E}$,  
we construct the adaptive interval 
\begin{align}\label{eq:adaptiveint}
    \big[ \tilde E_{( \gamma_n )}, \tilde E_{( n^2 - \gamma_n )}\big] \, \cap \, \big[ \tilde E_{( l_1 )}, \tilde E_{(l_2 )} \big], %\big[ \mathring E_{( \hat{k}_1 )}, \mathring E_{(\hat{k}_2 )} \big], 
\end{align} 
where $\tilde{E}_{(1)}\leqslant \cdots \leqslant \tilde{E}_{(n^2)}$
denote the 
ordered statistics of all elements in $\{\tilde{E}_{ij}:1\leqslant i,j\leqslant n \}$, $\gamma_n$ is a hyperparameter to trim extreme values, and $1\leqslant l_1<l_2\leqslant n^2$ denote the smallest and largest indices such that $\mu^{-1}(e)$ is well-defined for $e\in [ \tilde E_{( l_1 )}, \tilde E_{(l_2 )} ]$ and the link function $\mu(\cdot)$  in \eqref{eq:mulink}. 
% with $\mu(\mathbb{R})$ representing the image of the link function $\mu(\cdot)$ over $\mathbb{R}$. 
% $ \mathring E_{(l_1 -1 )} \not\in  \mu(\mathbb{R})$, and $ \mathring E_{(l_2 +1 )} \not\in  \mu(\mathbb{R})$. 

The two intervals in \eqref{eq:adaptiveint} serve two different purposes. 
The first one $[ \tilde E_{( \gamma_n )}, \tilde E_{( n^2 - \gamma_n )}]$ trims extreme values in a data-adaptive way, while the second interval $[ \tilde E_{( l_1 )}, \tilde E_{(l_2 )} ]$ constrains \eqref{eq:adaptiveint} to stay inside the image of $\mu$, so that $\mu^{-1}(\cdot)$ can be applied to any projected value. 
% Together, \eqref{eq:adaptiveint} defines a maximal feasible range that excludes extreme values and remains compatible with the inverse link. 
Further  discussions are given after Condition \ref{cond:meanfunction} below. 
Our proposed Algorithm \ref{algo:initialpgd} is connected to the initialization method  in \cite{ma2020universal} when the adaptive interval \eqref{eq:adaptiveint} is replaced by the fixed interval $[ e^{-M_1}/2, 1/2]$ and  $\mu(x)=e^x/(1+e^x)$. 
% Our proposed Algorithm \ref{algo:initialpgd} is connected to  Algorithm 3 in \cite{ma2020universal}  if replacing the projection set \eqref{eq:adaptiveint} with  $[ e^{-M_1}/2, 1/2]$ and taking $\mu(x)=e^x/(1+e^x)$.  
The advantage of our construction is that the adaptive interval $[ \tilde E_{( \gamma_n )}, \tilde E_{( n^2 - \gamma_n )}]$  produces a  suitably bounded interval  without relying on unknown model parameters.

\begin{algorithm}[!htbp]
\caption{Range-Adaptive Singular Value Thresholding} \label{algo:initialpgd}
%\setstretch{1.5}
\normalsize 
\KwIn{Data: $ A \in \mathbb{R}^{n\times n}$. \quad Hyperparameters: $\tau_n$,  $\gamma_n$.} 

    Let $\sum_{i=1}^{n} \sigma_{i} u_{i}v_{i}^{\top}$ denote the singular value decomposition of $ A$. 
     
    Let $\tilde{E} =\sum_{\{i:\, \sigma_{i} > \tau_n\}} \sigma_{i} u_{i} v_{i}^{\top}$. Elementwisely project $\tilde{E}$ onto \eqref{eq:adaptiveint} to obtain $\mathring{E}$. 
    
    Let $\mathring \Theta \in \mathbb R^{n \times n}$ with  entries 
    $\mathring{\Theta}_{ij}=\mu^{-1}\big(\mathring{E}_{ij}\big)$ across $ 1 \leqslant i,j\leqslant n$.  
    % \begin{align}\label{eq:ring_thetaij_def}
    %     \mathring  \Theta_{ij}  = \mu^{-1}\left[\Pi_{\mathcal{S}_{\mathring{E}} }\left(\mathring{E}_{ij} \right) \right], 
    % \end{align}
 %    \gray{\begin{align*}
 %     \mathring  \Theta_{ij}  = 
 %     \begin{cases} 
 %     \mu^{-1}(\mathring E_{ij}) & \text{if } \mathring E_{ij} \in \mu(\mathbb R) \cap \big[\mathring E_{( \gamma_n )}, \mathring E_{( n^2 - \gamma_n )}\big] \\
 %     0 & \text{otherwise}, 
 %     \end{cases}
 %    \end{align*}
 % where $\mathring{E}_{(m)}$ denotes the $m$-th order statistic of the entries of $\mathring{E}$. }
 
    Let $\mathring \alpha = (n \mathrm{I}_n + 1_n 1_n^\mytrans )^{-1} \mathring \Theta 1_n$.

    Let $\mathring Z = \mathcal{S}_{k}(\mathring{\Theta} - \mathring \alpha 1_n^\mytrans - 1_n \mathring \alpha^\mytrans)$, 
    % where $\mathcal{S}_k(M)=\sum_{i=1}^k \sigma_i\cdot 1(\sigma_i>0) u_iu_i^{\top}$, takes up to  of the eigendecomposition of $M=\sum_{i=1}^n\sigma_iu_iu_i^{\top}$. 
   where $\mathcal{S}_k(\cdot)$ denotes the operator that takes the top-$k$ positively truncated square root of the input matrix, formally defined in (H.1) in the Supplementary Material.  %\red{(YH: I forgot to update this?)}
   % denotes taking the 
   %  eigenvalue decomposition components corresponding to the positive eigenvalues up to $k$ of the input matrix.
   %  \red{(YH: TBD confirm $\mathring{\Theta}$ is symmetric)}
    % top-$k$ singular value components of the input matrix. (@Yuang. Shall we define it?)} 

\BlankLine
\KwOut{$\mathring Y = [\mathring Z, \mathring \alpha].$} 
\end{algorithm} 

\subsection{Asymptotic Theory for Range-Adaptive Singular Value Thresholding}\label{sec:initial_theory}
We next establish the asymptotic guarantee for the 
proposed RA-SVT estimator 
% output 
$\mathring Y=  [\mathring Z, \mathring \alpha]$ from  Algorithm \ref{algo:initialpgd}.
The goal is to show that $\mathring Y$   satisfies the initialization requirement in  Condition \ref{cond:initial1} under suitable conditions and tuning parameters. 
To this end, we impose regularity conditions on the link function in  \eqref{eq:mulink} and the hyperparameters in Conditions \ref{cond:meanfunction} and \ref{cond:usvt} below. %, where their implications discussed afterwards. 
%\red{(YH: Maybe add a condition on $\nu$, say $\nu\in \mu(\mathbb{R})$ and is a finite value??)}
 
% Similarly to Condition \ref{cond:parfunction}, we impose a regularity condition on the link function in  \eqref{eq:mulink}. 

\begin{condition}[Link function] \label{cond:meanfunction}
Assume $\mu(\theta)$ in \eqref{eq:mulink} satisfies the following conditions. 
\begin{itemize}\setlength{\itemsep}{0pt}

\item[(i)] The function $\mu(\theta)$ is continuously differentiable with its derivative denoted by $\mu^{\prime}(\theta)$. Moreover, for any constant $b > 0$, there exist constants $\kappa_4(b), \kappa_5(b) > 0$ such that
\begin{align*}
    \kappa_4(b) \leqslant \mu^{\prime}(\theta) \leqslant \kappa_5(b) \quad \text{for all } \theta \in [-b,b].
\end{align*}
\item[(ii)] For any $t \geqslant 0$, $\Pr(|A_{ij} - \mu(\Theta_{ij}^\star)| > t) \leqslant 2 \exp(-(t/K)^{s})$ for fixed constants $K,s>0$, which,  without loss of generality, can be assumed to be same as those in   Condition \ref{cond:parfunction} (ii). 
\end{itemize}
\end{condition}
Under the natural exponential family distributions with canonical link functions, 
Condition~\ref{cond:meanfunction}   follows directly from Condition~\ref{cond:parfunction}. Specifically, in that setting, 
  $ \mu'(\theta)=-\ell''(\theta;x)$ and $A_{ij}-\mu(\Theta_{ij}^{\star}) = \ell'(\Theta_{ij}^{\star}; A_{ij})$. Thus, 
Condition~\ref{cond:meanfunction} is equivalent to  Condition~\ref{cond:parfunction} with $\kappa_4(b)=\kappa_1(b)$ and $\kappa_5(b)=\kappa_2(b)$, suggesting no extra restrictions in the canonical setting. 
%\gray{Condition~\ref{cond:meanfunction} (i) implies that  over any closed interval on $\mathbb{R}$, $\mu(\cdot)$ is strictly increasing, and thus it is invertible with a well-defined inverse function $\mu^{-1}(\cdot)$.} \\
Moreover, Condition~\ref{cond:meanfunction}~(i) guarantees the  invertibility of $\mu(\cdot)$ and  existence of indexes $l_1$ and $l_2$ in \eqref{eq:adaptiveint}. 
In particular, it implies that 
% Moreover, Condition~\ref{cond:meanfunction} (i) implies that 
$\mu(\cdot)$ is continuous and strictly increasing over any closed interval on $\mathbb{R}$. Hence $\mu(\mathbb{R})$, the image of $\mu(\cdot)$ over $\mathbb{R}$,  is an interval, 
and its inverse function $\mu^{-1}(x)$ is well-defined for all $x\in \mu(\mathbb{R})$. 
Therefore, $l_1$ and $l_2$ in \eqref{eq:adaptiveint} exist and are easy to identify.   
% are uniquely identified by comparing to the boundaries  of $\mu(\mathbb{R})$, which is computationally simple.  
For example, if $\mu(\mathbb{R})=(a,b)$, $l_1$ and $l_2$ are simply the smallest and largest indices  such that $\tilde{E}_{(l_1)}>a$ and  $\tilde{E}_{(l_2)}<b$.  
This guarantees that  applying $\mu^{-1}(\cdot)$ after the projection onto the adaptive interval \eqref{eq:adaptiveint} is valid.  
Canonical link functions, such as $\mu(x)=e^x/(1+e^x)$ for  Bernoulli distribution and $\mu(x)=e^x$ for Poisson distribution, indeed satisfy Condition~\ref{cond:meanfunction} (i).

% Condition~\ref{cond:meanfunction} (i) is indeed satisfied. 

% $\mu(\cdot)$ is indeed strictly increasing and its image  $\mu(\mathbb{R})$ over the real line $\mathbb{R}$ is an interval.  Finding indices $l_1$ and $l_2$ in \eqref{eq:adaptiveint}  is easy and can be achieved by examining two ends points of $\mu(\mathbb{R})$.  

% with the use of an inverse function
% $[-b,b]$ for any constant $b>0$. Thus $\mu(\cdot)$ is invertible with well-defined inverse function $\mu^{-1}(\cdot)$ over $[-b,b]$, consistent with the use of an inverse function applied in Algorithm \ref{algo:initialpgd}.  

%\red{(YH: Can we further justify and say $\mu^{-1}$ applicable for whole real line?)}

% Define a \gray{curvature} \blue{(YH: need a different word)} maximum gradient mapping for $\mu(\cdot)$ as 
% \begin{align*}
%     \kappa_2(B) := \sup_{\substack{|\theta| \leqslant B}} \max \left\{ \mu'(\theta), \  1/\mu'(\theta) \right\}. 
% \end{align*}
% Under the natural exponential family distributions with canonical link functions, by   $ \mu'(\theta) = -\ell''(\theta;x) $, 
% we have $ \kappa_2(B)= \kappa_1(B) $ defined in \eqref{eq:kappa1r}, suggesting $\kappa_2(\cdot)$ also reflects the curvature of the log likelihood function. 

\begin{condition}[Tuning parameters] \label{cond:usvt}
    % Assume 
    The tuning parameters of Algorithm \ref{algo:initialpgd} satisfy:
    \begin{enumerate}
       \item[(i)]\ Threshold: $n^{\frac{1}{2}} \ll \tau_n \ll n^{\frac{1}{2} + \frac{1}{k+3}}$;
        \item[(ii)] Trimming quantile: $\tau_n n^{\frac{3}{2} - \frac{1}{k+3}} \ll \gamma_n \ll n^2$.
    \end{enumerate} 
\end{condition}

% In Algorithm \ref{algo:initialpgd}, 
The threshold $\tau_n$ controls how much truncation is applied to $A$. 
% First, $\tau_n$ is used to threshold the singular values of $A$. 
% Our lower bound requirement $\tau_n\gg n^{1/2}$  is slightly larger than the 
Original USVT \citep{chatterjee2015matrix} recommends a threshold of order $n^{1/2}$, 
% suggests choosing $\tau_n$ equal to $n^{1/2}$ up to a constant, 
which is slightly smaller than our lower bound $\tau_n\gg n^{1/2}$ in Condition~\ref{cond:usvt}. 
This difference arises because \cite{chatterjee2015matrix} assumes that entries of the observed data matrix are  bounded with probability one,  whereas our framework also accommodates unbounded data with a unified analysis.  
% whereas we do not impose such assumptions to accommodate for  unbounded data under a unified framework.  
In practice, however, the gap is insignificant. 
% But this difference is not very significant. 
For example, we can choose $\tau_n=n^{1/2}\log n$, which differs from $n^{1/2}$ only by a slow $\log n$ rate and  already satisfies Condition~\ref{cond:usvt} (i).  Condition~\ref{cond:usvt} (i) is stated as a range instead of a fixed choice  for flexibility and generalizability. 
% While  such a fixed choice is acceptable,    Condition~\ref{cond:usvt} (i) further provides a range of feasible $\tau_n$ for flexibility and generalizability. 
Intuitively,  $\tau_n$ should  be sufficiently small to preserve useful signals in $\mathbb{E}(A)$, but within a suitable range, a slightly larger $\tau_n$ may help remove more noises. The optimal choice of $\tau_n$ could be an interesting question. But this is not pursued in this paper as Algorithm \ref{algo:initialpgd}  is only used as an initialization for Algorithm \ref{algor:pgd}.

The trimming quantile $\gamma_n$ controls the length of the adaptive projection interval \eqref{eq:adaptiveint}. 
For a fixed $\tilde{E}$, increasing $\gamma_n$   shortens   $[ \tilde E_{( \gamma_n )}, \tilde E_{( n^2 - \gamma_n )}]$, and then  typically shortens \eqref{eq:adaptiveint} too. 
% and then  \eqref{eq:adaptiveint} is expected to be shorter in general. 
% An appropriate $\gamma_n$ allows the adaptive interval \eqref{eq:adaptiveint} to cover the true parameters in the matrix $\mu^{\star}:= \mu(\Theta^{\star})$  asymptotically without knowing $\Theta^{\star}$ (or $M_1$ in Condition~\ref{cond:truevalue}). 
The admissible range of $\gamma_n$ in Condition~\ref{cond:usvt}~(ii) reflects a tradeoff. 
On the one hand, $\gamma_n$ needs to be sufficiently large to screen out extreme entries in  $\tilde{E}$ and produce a  bounded interval, which is essential for     effective  regularization. 
% This motivates the lower bound in Condition~\ref{cond:usvt}, 
Technically, the lower bound in Condition~\ref{cond:usvt}
 arises from our derived  error rate $\|\tilde{E}-\mu(\Theta^{\star})\|_{\F}^2=O_p(\tau_n n^{\frac{3}{2} - \frac{1}{k+3}})$. 
% Intuitively, since $\tau_n n^{\frac{3}{2} - \frac{1}{k+3}}$ delineates the overall error order in $\tilde{E}$, it could serve as a useful indication for the least proportion of noises for truncation in \eqref{eq:adaptiveint}. 
Our intuition of setting $\tau_n n^{\frac{3}{2} - \frac{1}{k+3}}$ as the  lower bound is that it  delineates the  order of the overall estimation error in $\tilde{E}$, 
% error order in $\tilde{E}$ 
and thus serves as a useful benchmark for the  minimum trimming level. 
% indication for the smallest  proportion of noises for truncation. 
% our argument requires a lower bound $\gamma_n\gg \|\tilde{E}-\mu^{\star}\|_{\F}^2$, which is unknown and thus approximated by utilizing $\|\tilde{E}-\mu^{\star}\|_{\F}^2=O_P(\tau_n n^{\frac{3}{2} - \frac{1}{k+3}})$. 
% unknown. The specific $\tau_n n^{\frac{3}{2} - \frac{1}{k+3}}$ in the lower bound arises from an upper bound of the F-norm error $\|\tilde{E}-\mu^{\star}\|_{\F}^2$, which naturally depends on the threshold $\tau_n$ and latent dimension $k$. 
% The lower bound $\gamma_n\gg \tau_n n^{\frac{3}{2} - \frac{1}{k+3}}$ arises from the   Frobenius norm 
On the other hand, if $\gamma_n$ is too large,   the entries in $\mu(\Theta^{\star})$ are increasingly likely to fall outside the interval in  \eqref{eq:adaptiveint}, and the induced bound of $\operatorname{dist}^2(\mathring{Y}, Y^{\star})$  in \eqref{eq:error_usvt} below will deteriorate. A suitably large $\gamma_n$ allows \eqref{eq:adaptiveint} to cover the entries in $\mu(\Theta^{\star})$  asymptotically without knowing $\Theta^{\star}$ (or $M_1$ in Condition~\ref{cond:truevalue}).  
The upper bound  $\gamma_n\ll n^2$ is natural, because  the induced error bound in  \eqref{eq:error_usvt} is consistent with  existing results; see more discussions below.    

\begin{theorem} \label{thm:homousvt}
    Assume Conditions \ref{cond:truevalue}, \ref{cond:parfunction}, \ref{cond:meanfunction}, and \ref{cond:usvt}. For any constant $\varepsilon > 0$, there exist constants  $C > 0$ and $ N_\varepsilon \in \mathbb N$ such that when $n \geqslant N_\varepsilon$,
    \begin{align}\label{eq:error_usvt}
 \operatorname{dist}^2(\mathring Y, Y^\star)
  \leqslant  \frac{C\gamma_n}{n} \quad \text{ and }\quad  \|\mathring Y\|_{2\to \infty} \leqslant C
\end{align}
with probability $1 - O(n^{-\varepsilon})$.
\end{theorem}

Theorem~\ref{thm:homousvt} establishes two key bounds for $\mathring Y= [\mathring Z, \mathring \alpha]$. The first bound implies that the squared Frobenius norm error for estimating $Y^{\star}$ is $o(n)$, as $\gamma_n=o(n^2)$ by Condition~\ref{cond:usvt}. This is consistent with conclusions in \cite{chatterjee2015matrix} and  \cite{ma2020universal}. 
The second bound 
shows that $\|\mathring Y\|_{2\to \infty}$ remains    bounded.
Such a result is essential for fulfilling the requirement in Condition~\ref{cond:initial1} but has not been established in \cite{chatterjee2015matrix} or \cite{ma2020universal}. 
% Such a result is essential to justify the use of $\mathring{Y}$ as an initialization of Algorithm \ref{algor:pgd} by Condition~\ref{cond:initial1}, but it has  not been established in \cite{chatterjee2015matrix} or \cite{ma2020universal}. 
In summary, Theorem~\ref{thm:homousvt}  suggests that there exist choices of $\gamma_n$ and $\tau_n$ such that the resulting estimator $\mathring{Y}$ satisfies the requirements of initial estimates in Condition~\ref{cond:initial1}. For example, we can set 
\begin{align}\label{eq:gamma_n_tau_n}
   \gamma_n \asymp  n^{2 - \varsigma_0} \text{ for any  } \varsigma_0 \in (0,1/(k+3)), \ \text{ and any }\  n^{\frac{1}{2}}\ll \tau_n\ll n^{\frac{1}{2}+\frac{1}{k+3} - \varsigma_0}.  
\end{align}
% $\gamma_n \asymp  n^{2 - \varsigma_0}$ for any $\varsigma_0 \in (0,1/(k+3))$ and any $n^{\frac{1}{2}}\ll \tau_n\ll n^{\frac{1}{2}+\frac{1}{k+3} - \varsigma_0}$. 
As discussed after Corollary \ref{cor:line_iteration}, 
stronger conclusions, such as consistent two-to-infinity error  control, may be available under additional structural assumptions on the expected adjacency matrix $\mathbb{E}(A)$. However, such assumptions are generally unavailable or hard to justify given a  nonlinear link in \eqref{eq:mulink}. 
Instead, Theorem~\ref{thm:homousvt} is tailored to the  error control actually needed for the initialization under relatively weak  structural assumptions. 
% there exist stronger conclusions on consistent two-to-infinity error  control, but those may require more structural assumptions 
% % explicit low-rankness 
% imposed on   $\mathbb{E}(A)$, which cannot be guaranteed with a nonlinear link in \eqref{eq:mulink}. 
% In summary, 
% Theorem~\ref{thm:homousvt} shows that 

\section{Simulations}\label{sec:simulations}

This section provides   simulation studies to assess  both  the  algorithmic performance and asymptotic theory developed in this  paper.  
We begin with describing the simulation setup used throughout the section.  
% to demonstrate the partical algorithmic 
% evaluates the practical behavior of the proposed computational pipeline, consisting of Range-Adaptive Singular Value Thresholding  initialization (RA-SVT, Algorithm~\ref{algo:initialpgd}) and Projected Gradient Descent with Backtracking Line Search (line-search GD, Algorithm~\ref{algor:pgd}). We organize the presentation as follows. 
% In particular, 
Section~\ref{sec:simu-fs-vs-ls} studies the convergence behavior, and Section~\ref{sec:simu-an} examines the distribution  of the obtained estimator.  

% gradient descent  across a range of initial step sizes, providing empirical evidence related to Theorem~\ref{thm:homoPGD} and Corollary~\ref{cor:line_iteration}. We also compare the line-search GD with fixed-step GD to show that the backtracking scheme helps to find appropriate step sizes. Section~\ref{sec:simu-an} examines the asymptotic normality pattern of the line-search estimator, in connection with Theorem~\ref{thm:homoMLE}. 
% Section~\ref{sec:simu:rasvt-vs-usvt} \red{...}

% \subsection{Data generation} 
% \label{sec:simu-generation}
% \paragraph*{Set-up} 
\vspace{1pt} 
\textit{Simulation setup.}\, Following \cite{ma2020universal}, we generate  $\alpha^{\star}$ by first sampling $\tilde\alpha \in \mathbb{R}^n$ with independent entries     from the uniform distribution on $[1,3]$ 
% Uniform$[1,3]$ 
and then normalizing it to  $\alpha^\star=\tilde\alpha/(1_n^{\top}\tilde{\alpha})$. 
In addition, % we generate the model parameters $(Z^\star, \alpha^\star)$ as follows.  
% Let $k=2$. 
 we let $k=2$ and generate random $\tilde Z \in \mathbb{R}^{n\times k}$ with independent entries following the  standard normal truncated to the interval $[-2,2]$. 
% $\mathcal N_{[-2,2]}(0,1) $ over .
% where $\mathcal N_{[-2,2]}(0,1)$ is the standard normal distribution restricted onto the interval $[-2,2]$. 
Then we construct  $Z^\star = 0.5\sqrt{n}\,\tilde U$, where $\tilde U$ denotes the left singular vectors of the centered matrix $J_n\tilde Z$ with $J_n$  defined as in \eqref{eq:direction_def}. 
% $\tilde Z - 1_n 1_n^\top \tilde Z / n$. 
This construction yields two equal nonzero eigenvalues of $Z^{\star\top}Z^{\star}$, thereby placing the simulations in the challenging repeated-eigenvalue regime studied in this paper. 
% has two equal nonzero singular values, which is the challenging scenario that this paper addresses. 
% Denote by $z_i^{\star\top}$ the $i$-th row of $Z^\star$.

% $\alpha_i^\star=-\tilde\alpha_i/\sum_{j=1}^n\tilde\alpha_j$ for $i=1,\ldots,n$.

Given $(Z^{\star},\alpha^{\star})$, 
we generate network data $A$ from the model \eqref{eq:model} under three common distributions, Poisson, Bernoulli, and Gaussian with their respective canonical links. Under each distribution, we take  $n\in\{500,1000,2000,4000\}$. The main text  presents results under the Poisson distribution,  
while the  results under the Bernoulli and Gaussian distributions are similar  and are deferred to {Section J} of the Supplementary Material. 

% In this section, we consider data generated from the Poisson model
% \[
% \Theta_{ij}^\star=\alpha_i^\star+\alpha_j^\star+z_i^{\star\top}z_j^\star,\qquad
% A_{ij}\mid \Theta_{ij}^\star\sim\mathrm{Poisson}\!\left(\exp(\Theta_{ij}^\star)\right)\ , \qquad 1\leqslant i<j\leqslant n.
% \]
% with $A$ symmetrized and diagonal set to zero. 

% \vspace{1pt} 
% \textit{Implementation details.}\, 
% In the implementation of the algorithms,
% for RA-SVT, following Condition~\ref{cond:usvt}, we use
% 
% $\tau_n=\sqrt{\log n\cdot \sum_{i,j}A_{ij}/n}$ and $\gamma_n=0.1\,n^{2 -1/(k+4)}$.
The proposed projected gradient descent (Algorithm~\ref{algor:pgd}) is implemented with  
$C_{\mathrm{ls}}=1, \beta=0.5$, $R=2000$, and $R'$
equal to the ceiling of $\log n$. 
% $R'=\lceil\log n\rceil$, where $\lceil x \rceil$ denotes the smallest integer greater than or equal to $x$. 
% For comparison, the fixed-step version of gradient descent (fixed-step GD) is implemented as Algorithm~\ref{algor:pgd} without lines 5-9.
% For both line-search GD and fixed-step GD, 
The initial estimator $Y^0$ is constructed by Algorithm~\ref{algo:initialpgd}. Following \eqref{eq:gamma_n_tau_n}, we choose $\gamma_n=0.1\,n^{2 -\varsigma_0}$ with $\varsigma_0=1/(k+4)$ and  $\tau_n=( \upsilon_n n\log n )^{1/2}$ with $\upsilon_n= \sum_{1\leqslant i, j\leqslant n} A_{ij}/n^2$. Here $\tau_n$ satisfies \eqref{eq:gamma_n_tau_n} with high probability, as $\upsilon_n$ concentrates around a constant. 
% with the RA-SVT output, and the GD iterations 
Following Remark~\ref{rmk:stopping}, we adopt an early stopping criterion  that 
$S_{\max}^r:= \max_{1\leqslant i\leqslant n,\,  1\leqslant j\leqslant k+1}|\nabla_{y_{ij}} L(Y^r)|\leqslant 0.01$, where $\nabla_{y_{ij}}L(Y)$ denotes the partial derivative of $L(Y)$ with respect to $y_{ij}$, i.e., the $(i,j)$-th entry of $\nabla_Y L(Y)\in \mathbb{R}^{n\times (k+1)}$. 

% $S_{\max}^r:=\|\nabla_{Y} L(Y^r)\|_{1 \to \infty}\leqslant 0.01$.  
% terminate when 
% \begin{align}\label{eq:earlystop}
%     % \|\nabla_Y L(Y^r)\|_{\max}
%     S_{\max}^r:=\max_{1\leqslant i\leqslant n,\,  1\leqslant j\leqslant k+1}|\nabla_{y_{ij}} L(Y^r)|\leqslant 0.01.
% \end{align}

% $\|\nabla_Z L(Y^r)\|_{\max}\leqslant 0.01$ or the number of iterations reaches $R=2000$; these choices are consistent with the discussions in Remark~\ref{rmk:stopping}.  
% For any latent position matrices $Z_1, Z_2$, we compute their distance by $\operatorname{dist}^2(Z_1, Z_2) = \|{Z_1}-Z_2 VU^{\top}\|_{\mathrm{F}}^2$ where $U\Sigma V^{\top}$ is the singular value decomposition of ${Z}_1^{\top}Z_2$; see  \cite{schonemann1966generalized}.

\subsection{Algorithmic Performance of the Projected Gradient Descent} \label{sec:simu-fs-vs-ls}

% This subsection examines the convergence of the projected gradient descent algorithm and demonstrates data adaptiveness of the proposed line search conditions.  
% backtracking behavior of line-search GD. 
% We first give a representative comparison between line-search GD and fixed-step GD to illustrate why adaptive step-size selection is needed in practice, and then summarize the behavior of line-search GD across a range of initial step sizes. 
We examine the adaptivity and computational trade-off of the proposed algorithm. 
To this end, we vary the initial step size and compare the adaptive procedure with a fixed-step counterpart. 
% demonstrate the adaptiveness of the proposed line search conditions across different values of the initial step size $\eta_{\mathrm{init}}$. 
Specifically, we take $\eta_0$ to be six times  \eqref{eq:eta_init_choice},
% as in \eqref{eq:eta_init_choice}  multiplied by six, 
which serves as a baseline step size determined by the initial estimator, and consider four scales $\eta_{\mathrm{init}}/ \eta_0\in \{10, 5, 1, 1/5\}$. % to vary  scales.     
The fixed-step counterpart of  Algorithm~\ref{algor:pgd} 
% We  compare Algorithm~\ref{algor:pgd} with a fixed-step counterpart that 
sets  $\eta_r=\eta_{\mathrm{init}}$ throughout the iterations, equivalent to skipping  lines 5--9 in Algorithm~\ref{algor:pgd}. 
Each configuration is evaluated over 100 Monte Carlo replications. 
% In  each setting, the experiments are conducted 100 Monte Carlo replications. 
% In each replication, 
To measure convergence, we use $S_{\max}^r$  as a quantitative measure, because 
by Remark \ref{rmk:stopping},  $S_{\max}^r$ 
% $\|\nabla_Y L(Y^r)\|_{\max}$ 
converges to 0 if and only if $Y^r$ approximates the maximum likelihood estimator $\hat{Y}$ up to the identifiability constraints.
% , making  $S_{\max}^r$ a natural quantitative measure of convergence.
 % so a small $S_{\max}^r$ is a quantitative indication on convergence. 
Accordingly, we let $R_{\mathrm{conv}}$ denote the smallest iteration $r$  such that $S_{\max}^r\leqslant 0.01$,
% the stopping rule $S_{\max}^r\leqslant 0.01$ is reached, 
and declare convergence if $R_{\mathrm{conv}}\leqslant R=2000$. 

% We say the projected gradient descent converges if \eqref{eq:earlystop} is satisfied for $r\leqslant R$. 

 Table~\ref{tab:conv_poisson_n_eta} presents the empirical convergence proportions under both adaptive  and fixed step sizes. 
 % , respectively. of replications where the algorithm converges under adaptive  and fixed step sizes, respectively. 
 % The adaptive method converges in all settings
Using the adaptive step size, the projected gradient descent algorithm converges reliably under all values of $\eta_{\mathrm{init}}$. 
 In contrast,  the fixed-step  method  fails to converge when $\eta_{\mathrm{init}}\in \{10,5\}$, indicating high sensitivity to the initial step size. 
 
 % becomes large, the algorithm does not converge across 100 replications. 
 
 % under the Poisson distribution.  
% We study line-search GD under four initial step sizes values $\eta_{\mathrm{init}}/ \eta_0\in \{10, 5, 1, 1/5\}$,  
% $
% \eta_{\mathrm{init}}\in\{10\eta_0,\ 5\eta_0,\ \eta_0,\ \eta_0/5\},
% $ 
% where $\eta_0=6\times \eqref{eq:eta_init_choice} $ represents a baseline step size determined by the initial estimator $Z^0$,  
% \[
% \eta_0=
% \frac{1}{\max_{1\leqslant i<j\leqslant n}\{1-\ell''(\Theta_{ij}^0)\}}
% \cdot
% \min\!\left\{\frac{1}{\|Z^0\|_{\op}^2},\frac{1}{n}\right\}
% \]
% is the step size in \eqref{eq:eta_init_choice} of Corollary~\ref{cor:line_iteration} with the theoretical constant removed, 
% and $Z^0$ is constructed from the RA-SVT initializer. 

\begin{table}[!htbp]
\centering
\setlength{\tabcolsep}{9pt}   % default is 6pt
\caption{Empirical convergence proportions over 100 Monte Carlo replications under the Poisson model.} 
\label{tab:conv_poisson_n_eta}
% \scriptsize
\begin{tabular}{r c cccc c cccc}
\toprule
\multicolumn{2}{c}{} &
\multicolumn{4}{c}{Adaptive step size} & \multicolumn{1}{c}{} & \multicolumn{4}{c}{Fixed step size} \\[2pt]
\cline{3-6}\cline{8-11}
\noalign{\vskip 0.35ex}
\multicolumn{2}{c}{} &
\multicolumn{4}{c}{$\eta_{\mathrm{init}}/\eta_0$} & \multicolumn{1}{c}{} & \multicolumn{4}{c}{$\eta_{\mathrm{init}}/\eta_0$} \\[2pt]
\noalign{\vskip 0.35ex}
$n$ &  & $10$ & $5$ & $1$ & $1/5$ &  & $10$ & $5$ & $1$ & $1/5$ \\
\cline{1-1}\cline{3-6}\cline{8-11}
\noalign{\vskip 0.45ex}
$500$  &  & $1.00$ & $1.00$ & $1.00$ & $1.00$ &  & $0.00$ & $0.00$ & $1.00$ & $1.00$ \\
$1000$ &  & $1.00$ & $1.00$ & $1.00$ & $1.00$ &  & $0.00$ & $0.00$ & $1.00$ & $1.00$ \\
$2000$ &  & $1.00$ & $1.00$ & $1.00$ & $1.00$ &  & $0.00$ & $0.00$ & $1.00$ & $1.00$ \\
$4000$ &  & $1.00$ & $1.00$ & $1.00$ & $1.00$ &  & $0.00$ & $0.00$ & $1.00$ & $1.00$ \\
\bottomrule
\end{tabular}
\end{table}

% \red{(YH: Why not add $n=4000$?)}

% To motivate the line-search scheme, 

% the solving  may not be reliable if 
% is practically important because a suitable step size is typically unknown a priori and may vary across datasets.
% More generally, one may also use the proposed line search conditions to help  choose   initial step sizes too.  

\begin{figure}[!htbp]
\centering
\begin{subfigure}[t]{0.39\linewidth}
\centering
\includegraphics[width=\linewidth]{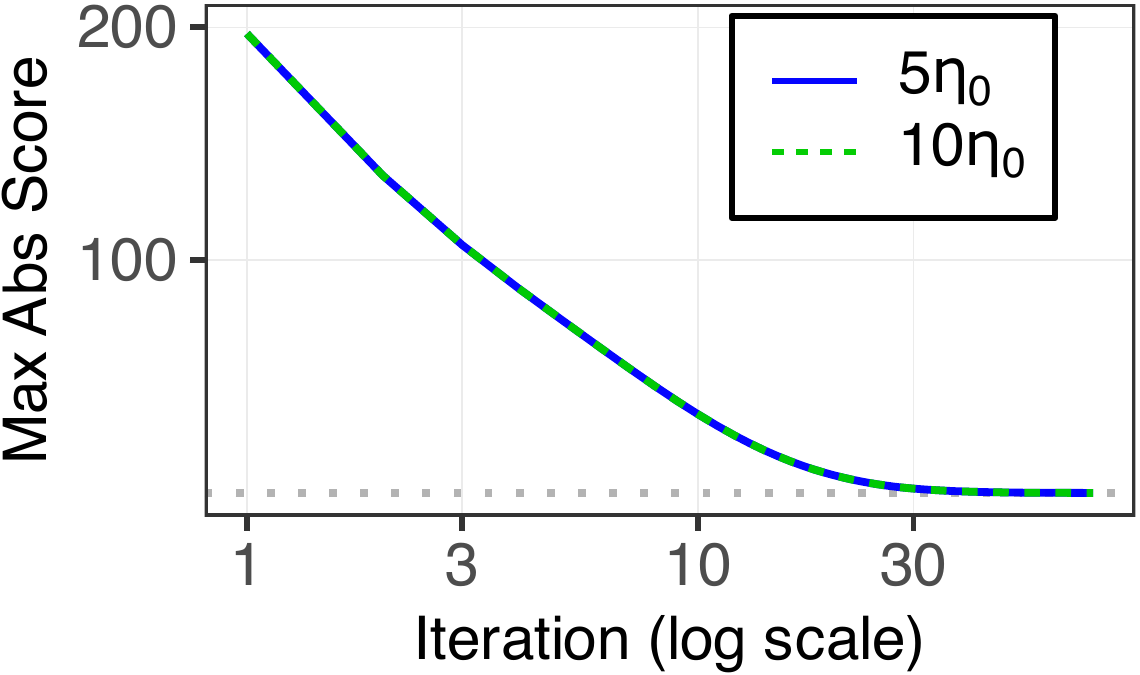}
\caption{Adaptive step size}\label{fig:adpt_score}
\end{subfigure}
\hspace{0.1\linewidth}
\begin{subfigure}[t]{0.39\linewidth}
\centering
\includegraphics[width=\linewidth]{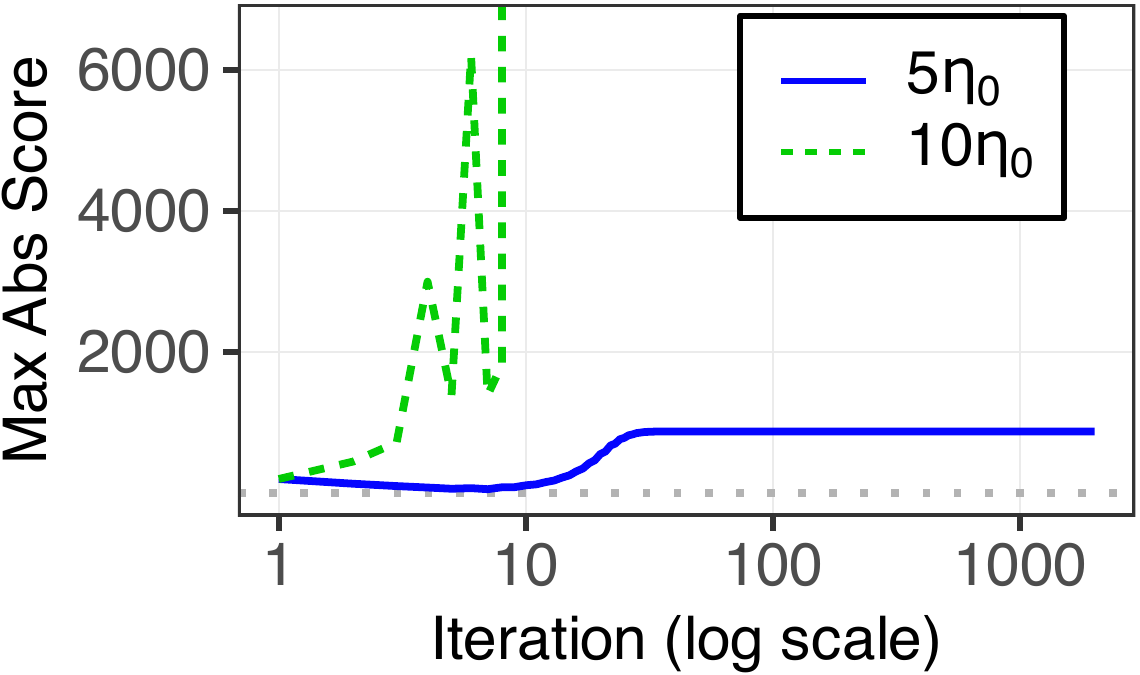}
\caption{Fixed step size} \label{fig:fixed_score}
\end{subfigure}
% \vspace{-0.8em}
\caption{Maximum absolute score $S_{\max}^r$ versus  iteration $r$ under the Poisson model with $n=1000$ in  one Monte Carlo replication. Panels (a) and (b) correspond to the adaptive and fixed step sizes, respectively. The gray dotted line marks the stopping threshold $0.01$. }
\label{fig:ll_score_seed3123_poisson_n1000}
% \vspace{0.25em}
% \parbox{1\linewidth}{\footnotesize 
% \textit{Notes.} 
% }
\end{figure}

% Figure~\ref{fig:ll_score_seed3123_poisson_n1000} compares the max score trajectories with $n=1000$ for fixed-step GD and line-search GD with $\eta_{\mathrm{init}}\in\{5\eta_0,10\eta_0\}$. 
% For fixed-step GD, the score either explodes after a few iterations, or remains far from zero throughout. By contrast, the corresponding line-search runs converge stably for both choices of $\eta_{\mathrm{init}}$. The same pattern is observed across 100 replications; we display only one run for clarity. These results show that fixed-step GD is more sensitive to the choice of step size than line-search GD. This distinction is practically important because a suitable step size is typically unknown a priori and may vary across datasets. As a result, fixed-step GD may require more careful tuning, whereas the proposed line-search procedure selects the step size adaptively and yields more stable performance. 
To  explain the patterns in  Table~\ref{tab:conv_poisson_n_eta}, Figure~\ref{fig:ll_score_seed3123_poisson_n1000} presents $S_{\max}^r$ against the iteration index $r$ for a randomly chosen replication with $n=1000$ and $\eta_{\mathrm{init}}/\eta_0\in \{10, 5\}$. 
% distinct  convergence patterns  under the adaptive and fixed step sizes, 
% we consider $\eta_{\mathrm{init}}/\eta_0\in \{10, 5\}$ and present the  maximum absolute score $S_{\max}^r$
% % $\|\nabla_Y L(Y^r)\|_{\max}$ 
% versus  $r$ in one Monte Carlo  replication  under  $n=1000$. 
% Recall Remark \ref{rmk:stopping} suggests that $S_{\max}^r$ 
% % $\|\nabla_Y L(Y^r)\|_{\max}$ 
% converges to 0 if and only if $Y^r$ approximates the maximum likelihood estimator $\hat{Y}$, up to the identifiability constraints. 
% under the adaptive and fixed step sizes, respectively. 
% Figure~\ref{fig:ll_score_seed3123_poisson_n1000} presents  $\|\nabla_Y L(Y^r)\|_{\max}$ over  $r$. 
Figure~\ref{fig:ll_score_seed3123_poisson_n1000}\subref{fig:adpt_score}  shows that the maximum absolute score steadily decreases to 0  when using the adaptive step size. 
By contrast,  Figure~\ref{fig:ll_score_seed3123_poisson_n1000}\subref{fig:fixed_score} shows that for the fixed-step method, $S_{\max}^r$ stagnates at a non-zero level 
% is trapped at a limiting non-zero value 
when $\eta_{\mathrm{init}}=5\eta_0$, and even diverges  when $\eta_{\mathrm{init}}=10\eta_0$. 
This suggests that the produced estimator can be far from $\hat{Y}$. 
% These trajectories are consistent with the consis observed patterns are consistent across all replications. 
These different trajectories demonstrate the effectiveness and adaptivity of the proposed line search method.
% without prior knowledge in practice.
% , This distinction suggests that the proposed line search conditions can guarantee approximating the maximum likelihood estimator, demonstrating a practical advantage. 

\begin{figure}[!htbp]
\centering
\begin{subfigure}[t]{0.46\linewidth}
\centering
\includegraphics[width=0.83\linewidth]{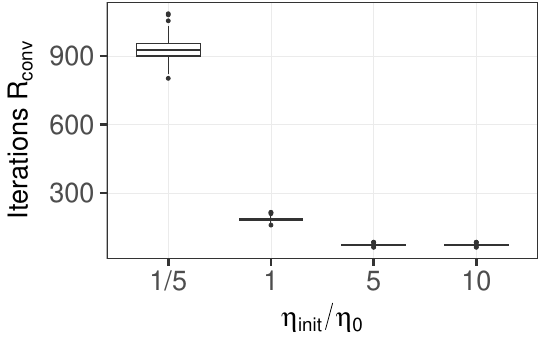}
\caption{Iterations to  convergence}\label{fig:gd_iter_num} 
\end{subfigure}
\hspace{0.05\linewidth}
\begin{subfigure}[t]{0.46\linewidth}
\centering
\includegraphics[width=0.83\linewidth]
{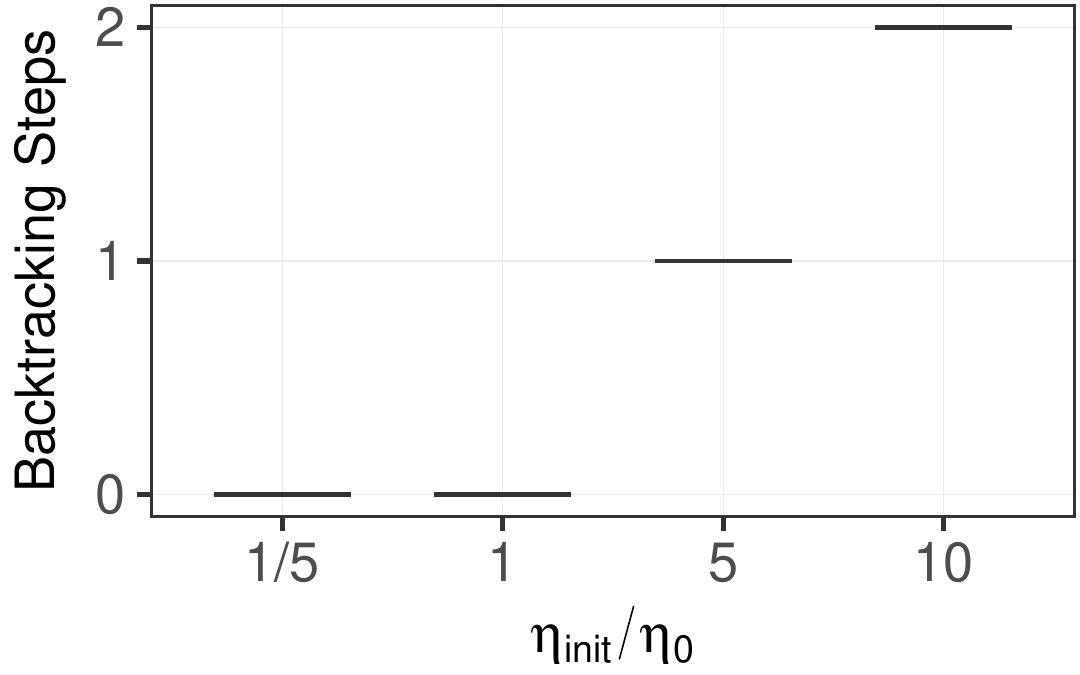}
% \includegraphics[width=0.83\linewidth]{figures/bubblecount_avg_backtracks_vs_eta_ratio_n1000_linesearch_poisson.pdf}
% \caption{Backtracking step averaged over iterations}\label{fig:count_back}
\caption{Backtracking steps}\label{fig:count_back}
\end{subfigure}
% \vspace{-0.3cm}
\caption{Computational trade-off under the adaptive step size for the Poisson model with $n=1000$ over 100 Monte Carlo replications. Panel (a) shows the boxplot of  $R_{\mathrm{conv}}$   versus $\eta_{\mathrm{init}}/\eta_0$. Panel (b) shows the boxplot of the backtracking steps per iteration  versus $\eta_{\mathrm{init}}/\eta_0$.}
% bubble plot of the average backtracking steps  versus $\eta_{\mathrm{init}}/\eta_0$; bubble labels show the numbers  of   replications attaining each value.}
\label{fig:poisson_n1000_ls_boxplots}
\end{figure}

Although the adaptive method converges across different $\eta_{\mathrm{init}}$, the choice still  affects computational efficiency. 
% Although Table~\ref{tab:conv_poisson_n_eta} shows that under the fixed step size,  the projected gradient descent algorithm can still converge  when $\eta_{\mathrm{init}}/\eta_0$ is relatively small. 
% However, we point out that using a smaller step size pays a price of more iterations of the projected gradient descent for convergence. 
% Figure~\ref{fig:poisson_n1000_ls_boxplots} illustrates the computational trade-off induced by the initial step size. 
% To illustrate that trade-off, 
Figure~\ref{fig:poisson_n1000_ls_boxplots}\subref{fig:gd_iter_num} shows the distribution of 
the number of gradient descent iterations to convergence, i.e., $R_{\mathrm{conv}}$,  
versus different $\eta_{\mathrm{init}}/\eta_0$. 
When  $\eta_{\mathrm{init}}$ is small,  the updates are more conservative, and  $R_{\mathrm{conv}}$  is correspondingly larger. 
% It shows that smaller $\eta_{\mathrm{init}}$ leads to higher $R_{\mathrm{cov}}$, i.e., more iterations of gradient descent to  converge. This is  because the resulting updates are conservative. 
As $\eta_{\mathrm{init}}$ increases, $R_{\mathrm{conv}}$   decreases  significantly.  
Moreover, Figure~\ref{fig:poisson_n1000_ls_boxplots}\subref{fig:count_back} presents the empirical distribution of the backtracking steps per iteration. 
% over 100 replications of implementing the projected gradient descent.  
% for each run of Algorithm~\ref{algor:pgd}, that is,  $\sum_{r=1}^{R_{\mathrm{conv}}}b_r/R_{\mathrm{conv}}$ where  $b_r$ represents the backtracking steps in the $r$-th iterations. 
% It shows that the backtracking steps increase with respect to $\eta_{\mathrm{init}}$. 
For each $\eta_{\mathrm{init}}/\eta_0$,  the   corresponding boxplot collapses to a line segment, because the backtracking steps remain identical across all gradient descent iterations and Monte Carlo replications. 
% Nevertheless, 
The results suggest  that more backtracking steps are consistently needed when $\eta_{\mathrm{init}}/\eta_0$ is large. 
% causing the degeneracy of the boxplots to lines in Figure~\ref{fig:poisson_n1000_ls_boxplots}\subref{fig:count_back}. 
% It shows that keep increasing $\eta_{\mathrm{init}}$  could result in more backtracking steps needed for convergence. 
% to reduce an overly large initial step size. 
% shows that    using the backtracking line search achieves smaller gradient descent iterations with  . 
In summary, Figure~\ref{fig:poisson_n1000_ls_boxplots} reveals the computational trade-off that larger initial step sizes reduce the number of iterations to convergence but could incur more backtracking steps. 
A moderately large initial step size offers a favorable balance between rapid descent  and limited backtracking. 
% Since the backtracking steps can increase mildly, we recommend choose a slightly large $\eta_{\mathrm{init}}$ to guarantee sufficiently quick decay and allows for moderate  backtracking. 

% \vspace{-2.5em}

\subsection{Asymptotic Distribution}
\label{sec:simu-an}
% As Section~\ref{sec:simu-fs-vs-ls} already demonstrates algorithmic convergence, 
We next present  the empirical distribution of the obtained estimator and examine whether  it aligns with the asymptotic theory in Theorem~\ref{thm:homoMLE}. 
%, showing its alignment to the asymptotic theory in Theorem~\ref{thm:homoPGD}. 
% We demonstrate that the empirical distribution of the obtained estimator aligns with the asymptotic distribution of the maximum likelihood estimator in Theorem~\ref{thm:homoPGD}. 
% In view of the convergence guarantee in Theorem~\ref{thm:homoPGD} and the numerical evidence in Section~\ref{sec:simu-fs-vs-ls}, the estimator computed by line-search GD can be regarded as a practical counterpart to the theoretical constrained MLE studied in Section~\ref{sec:mle_theory_1}. This subsection therefore examines the asymptotic normality pattern of the line-search estimator as an empirical illustration of Theorem~\ref{thm:homoMLE}. 
To this end, it suffices to consider the adaptive step size with $\eta_{\mathrm{init}}=\eta_0$ since the convergence behavior is stable as shown in Section~\ref{sec:simu-fs-vs-ls}. 
% since the adaptive method yields stable convergence across different. 
% $\eta_{\mathrm{init}}$ in Section~\ref{sec:simu-fs-vs-ls}. 
% and run line-search GD for $n\in\{500,1000,2000,4000\}$, with 200 replications for each value of $n$. 
Throughout this subsection, each setting is evaluated over 200 Monte Carlo replications.

% For simplicity, in this subsection we write the line-search estimator as $\hat Y=[\hat Z,\hat\alpha]$ and denote its aligned version, in the sense of \eqref{eq:rotation_z_def}, by $\hat Y_q=[\hat Z_q,\hat\alpha]$. 
% For the illustration of 
% We first examine the  asymptotic normality of an entry in the latent-vector estimator.  
% Specifically, 
We consider two representative estimands. The first is an entrywise latent-vector parameter. Specifically, consider  ${z}_{11}^{\star}$, the first 
coordinate of the node $1$ latent vector. 
% row and first column of $\hat{Z}_q$, denoted by $\hat{z}_{q,11}$. 
This corresponds to taking $g(y_1)= \langle e_1^{(k+1)}, y_1 \rangle$ in Corollary~\ref{cor:gyqIhat_AN}, where $e_1^{(k+1)}\in \mathbb{R}^{k+1}$ is the standard basis vector for the first coordinate in the $(k+1)$-dimensional space. It follows that the  standardized statistic  ${t}(\hat{z}_{q,11}):= 
(\hat z_{q,11}-z_{11}^\star)/\widehat{{se}}(\hat z_{q,11}) \xrightarrow{d} \mathcal N(0,1)$ with $\widehat{{se}}(\hat z_{q,11}):=\{e_1^{(k+1) \top}\Sigma_1(\hat{Y}_q)^{-1}e_1^{(k+1)}\}^{1/2} $. Second, we consider a nonlinear transformation corresponding to the edge mean $\mu(\Theta_{12}^{\star})$. By Corollary~\ref{cor:gyqIhat_AN}, we construct the standardized statistic $ t(\hat{\mu}_{12}) := \{\mu(\hat\Theta_{12}) - \mu(\Theta_{12}^\star)\}/\widehat{se}(\hat{\mu}_{12})$, where $\hat{\mu}_{12}:=\mu(\hat{\Theta}_{12})$ and  $\widehat{se}(\hat{\mu}_{12})$ is defined in \eqref{eq:se_mu}.

% To illustrate the entrywise asymptotic normality, 
% we consider the first row and first column of $\hat{Z}_q$, denoted by $\hat{z}_{q,11}$.   
% By Corollary~\ref{cor:gyqIhat_AN}, let ${t}(\hat{z}_{q,11}):= 
% (\hat z_{q,11}-z_{11}^\star)/\widehat{{se}}(\hat z_{q,11})$ denote the standardized statistic for $\hat{z}_{q,11}$, where $\widehat{{se}}(\hat z_{q,11})$ is obtained  by  $z_{11}^{\star}=e_1^{\top}y_1^{\star}$, and $e_1\in \mathbb{R}^{k+1}$ is $(k+1)$-dimensional vector with only the first entry being one. 
% % based on Corollary~\ref{cor:gyqIhat_AN}
% % Theorem~\ref{thm:homoMLE} and is given in {Section J} of the Supplementary Material. 
% % Figure~\ref{fig:qq_poisson_main} displays the QQ plots of ${t}(\hat{z}_{q,11})$ against the standard normal $\mathcal{N}(0,1)$. 
% In addition, we demonstrate the asymptotic distribution of the estimator after non-linear transformation. 
% As an example, we consider the estimator for  the edge mean $\mu(\Theta_{12}^{\star})$. Let $ t(\hat{\mu}_{12}) := \{\mu(\hat\Theta_{12}) - \mu(\Theta_{12}^\star)\}/\widehat{se}(\hat{\mu}_{12})$ denote the standardized statistic for $\hat{\mu}_{12}:=\mu(\hat\Theta_{12})$, where $\widehat{se}(\hat{\mu}_{12})$ is defined as in \eqref{eq:se_mu}.
% the precise form of $\widehat{se}(\hat{\mu}_{12})$ is obtained based on Corollary~\ref{cor:gyqIhat_AN} and is given in {Section J} of the Supplementary Material. 

\begin{figure}[!htbp]
\centering
\begin{subfigure}[t]{0.24\linewidth}
\centering
\includegraphics[width=\linewidth]{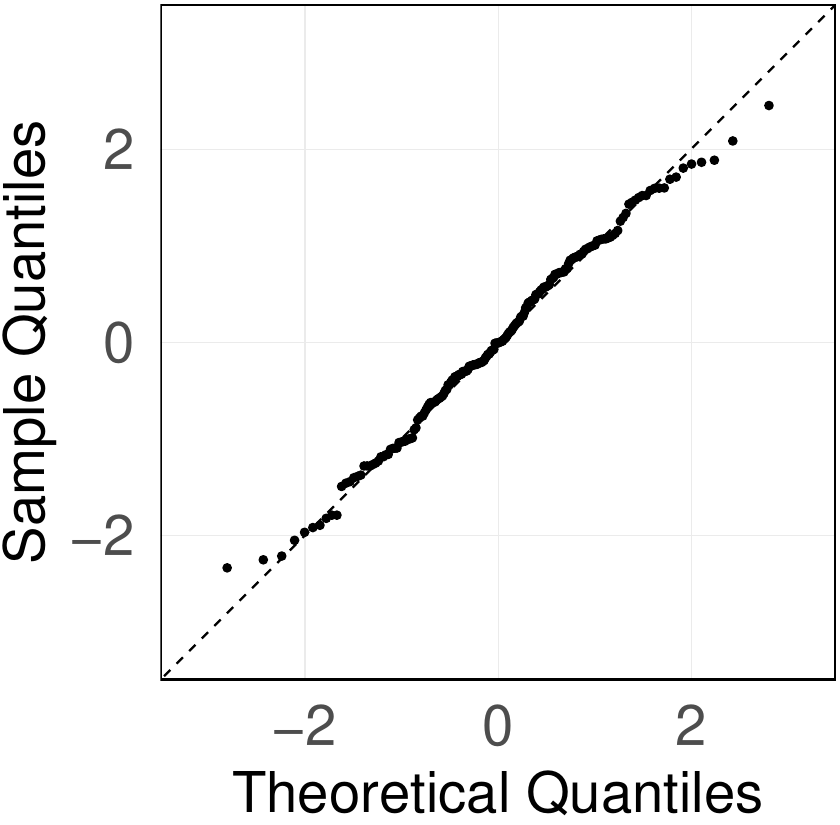}
\caption{$n=500$}
\end{subfigure}
\begin{subfigure}[t]{0.24\linewidth}
\centering
\includegraphics[width=\linewidth]{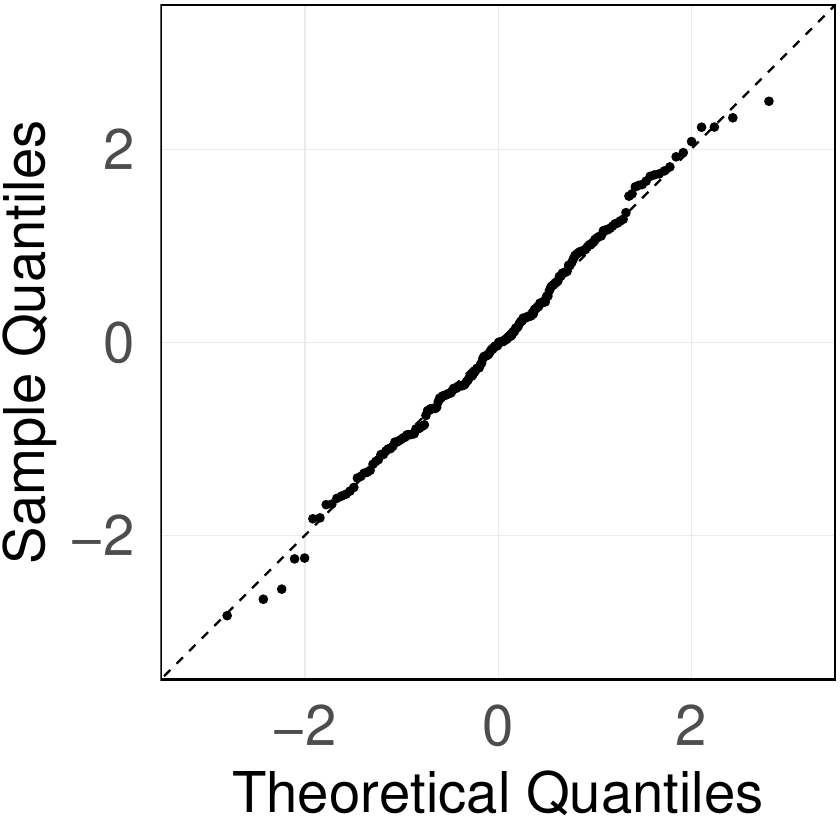}
\caption{$n=1000$}
\end{subfigure}
\begin{subfigure}[t]{0.24\linewidth}
\centering
\includegraphics[width=\linewidth]{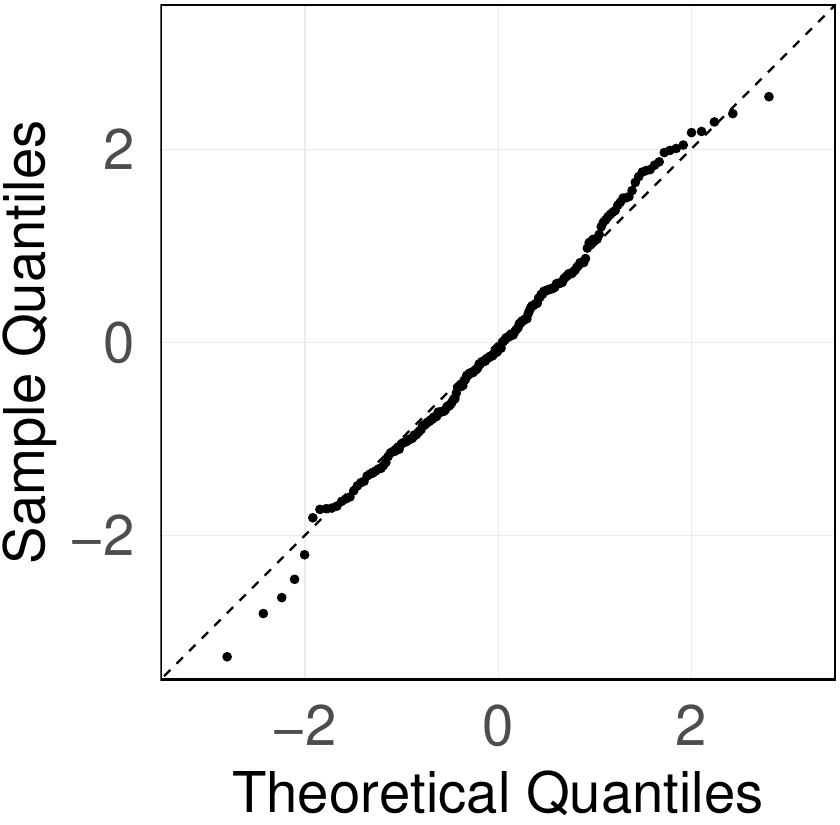}
\caption{$n=2000$}
\end{subfigure}
\begin{subfigure}[t]{0.24\linewidth}
\centering
\includegraphics[width=\linewidth]{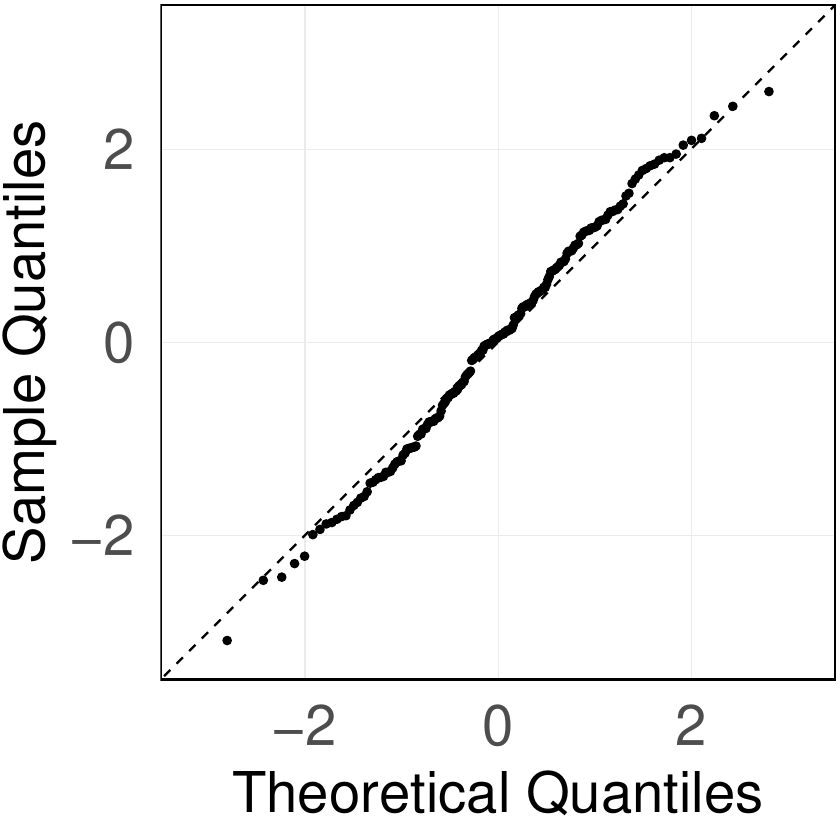}
\caption{$n=4000$}
\end{subfigure}
% \vspace{-0.3cm}
\caption{QQ plots of $t(\hat z_{q,11})$ against $\mathcal N(0,1)$ under the Poisson model.
% with $n\in\{500,1000,2000,4000\}$.
}
\label{fig:qq_poisson_main}
\end{figure}

\begin{figure}[!htbp]
\centering
\begin{subfigure}[t]{0.24\linewidth}
\centering
\includegraphics[width=\linewidth]{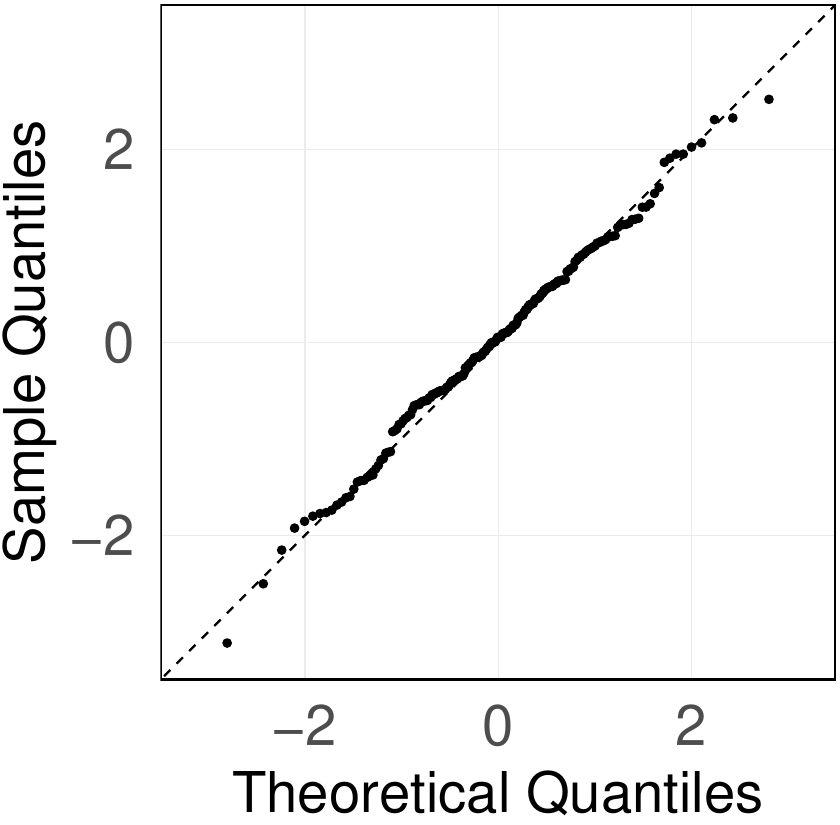}
\caption{$n=500$}
\end{subfigure}
\begin{subfigure}[t]{0.24\linewidth}
\centering
\includegraphics[width=\linewidth]{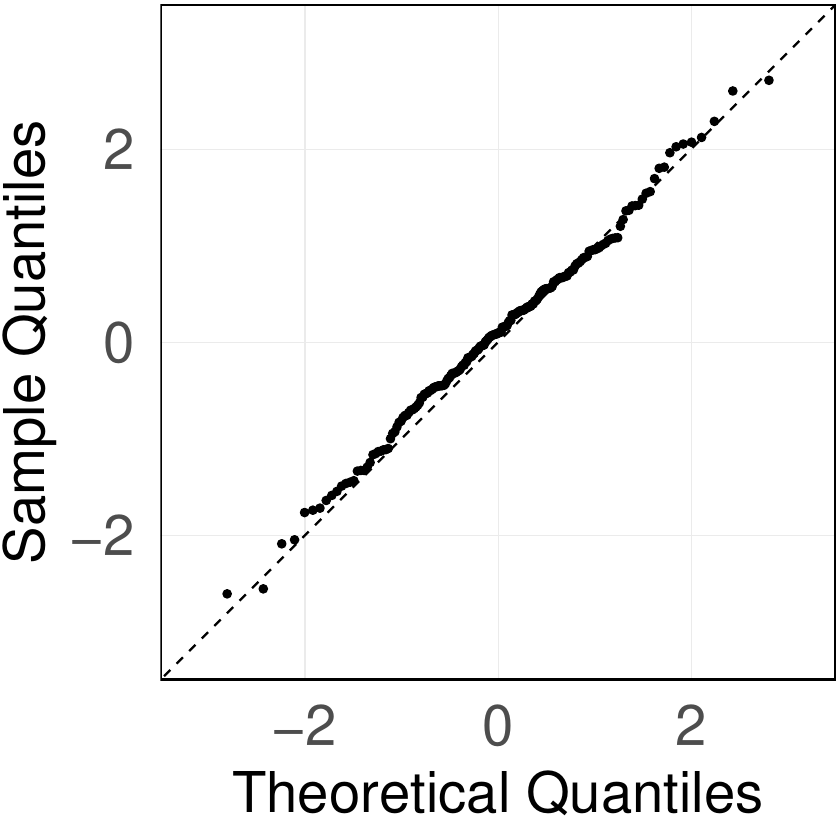}
\caption{$n=1000$}
\end{subfigure}
\begin{subfigure}[t]{0.24\linewidth}
\centering
\includegraphics[width=\linewidth]{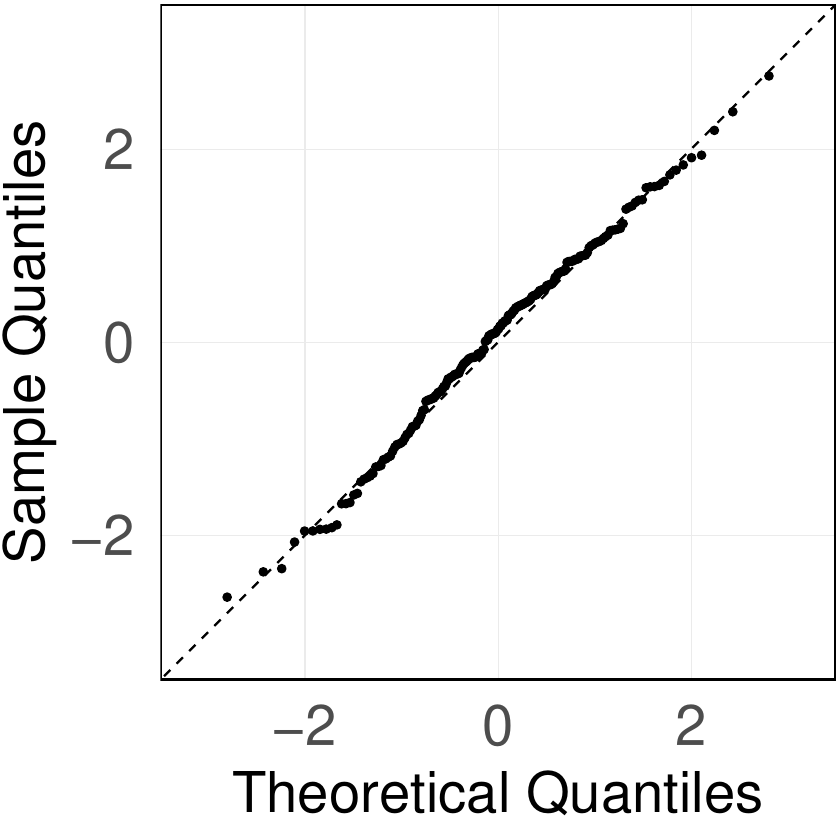}
\caption{$n=2000$}
\end{subfigure}
\begin{subfigure}[t]{0.24\linewidth}
\centering
\includegraphics[width=\linewidth]{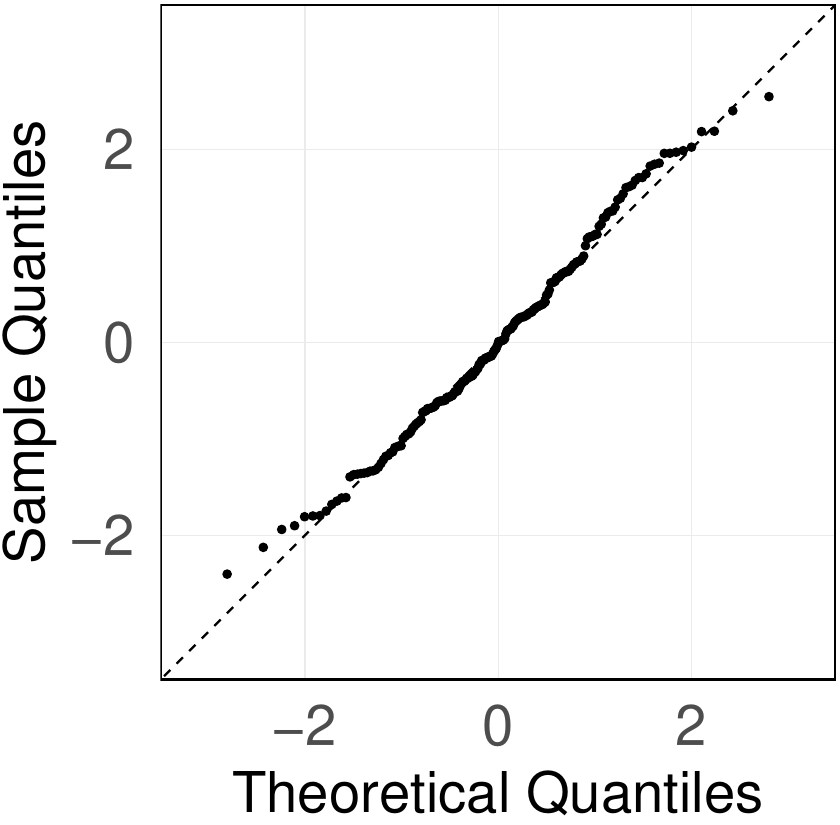}
\caption{$n=4000$}
\end{subfigure}
% \vspace{-0.3cm}
\caption{QQ plots of $t(\hat \mu_{12})$ against $\mathcal N(0,1)$ under the Poisson model.
% with $n\in\{500,1000,2000,4000\}$.
}
\label{fig:qq_poisson_mu12}
\end{figure}

 Figure~\ref{fig:qq_poisson_main} displays the QQ plots of ${t}(\hat{z}_{q,11})$ against the standard normal $\mathcal{N}(0,1)$, and 
Figure~\ref{fig:qq_poisson_mu12} presents the corresponding  QQ plots for $t(\hat{\mu}_{12}) $. % against $\mathcal{N}(0,1)$. 
In both cases, the empirical distributions of the standardized statistics closely align with the standard normal distribution. This supports the normal approximations by our asymptotic theory.  
Notably, the results are obtained when  $Z^{\star\top}Z^{\star}$ has repeated eigenvalues, demonstrating that a distinct-eigenvalue assumption is not essential for the asymptotic distributions considered here.

% \begin{figure}[htbp]
% \centering
% \includegraphics[width=1\linewidth]{figures/qq_direct_merged_poisson_rep200.pdf}
% \caption{QQ-plots of $\hat t(z_{11})$ for Poisson model with $n\in\{500,1000,2000,4000,8000\}$.}
% \label{fig:qq_poisson_main}
% \end{figure}

% \begin{figure}[htbp]
% \centering
% \includegraphics[width=1\linewidth]{figures/qq_muTheta12_direct_merged_poisson.pdf}
% \caption{QQ-plots of $\hat t(\mu_{12})$ for Poisson model with $n\in\{500,1000,2000,4000,8000\}$.}
% \label{fig:qq_poisson_mu12}
% \end{figure}

\section{Data Analysis}\label{sec:data}

To illustrate the proposed inference procedures, we analyze  the New York Citi Bike dataset \citep{bikedata}.  
The raw data record rides between   bike stations in New York City. 
Each ride  contains two stations and a start time. 
% two stations and a trip starting time.  
We aim to compare the travel patterns  during two commuting peak periods 8:00-9:00 and 18:00-19:00 on the weekday August 1st, 2019. 
We preprocess the data by keeping the rides that last between one minute and 3 hours and focusing on stations whose total degrees are  at least one in both  hours. 
% with total degrees that are smaller and equal than 1 in either of the two hours. 
The processed data contains 703 bike stations and 7,426  and  9,264 rides over the two  peak hours, respectively.

\begin{figure}[!htbp] 
    \centering
   \begin{subfigure}{0.32\textwidth}
		\centering
		\includegraphics[width=1\linewidth]{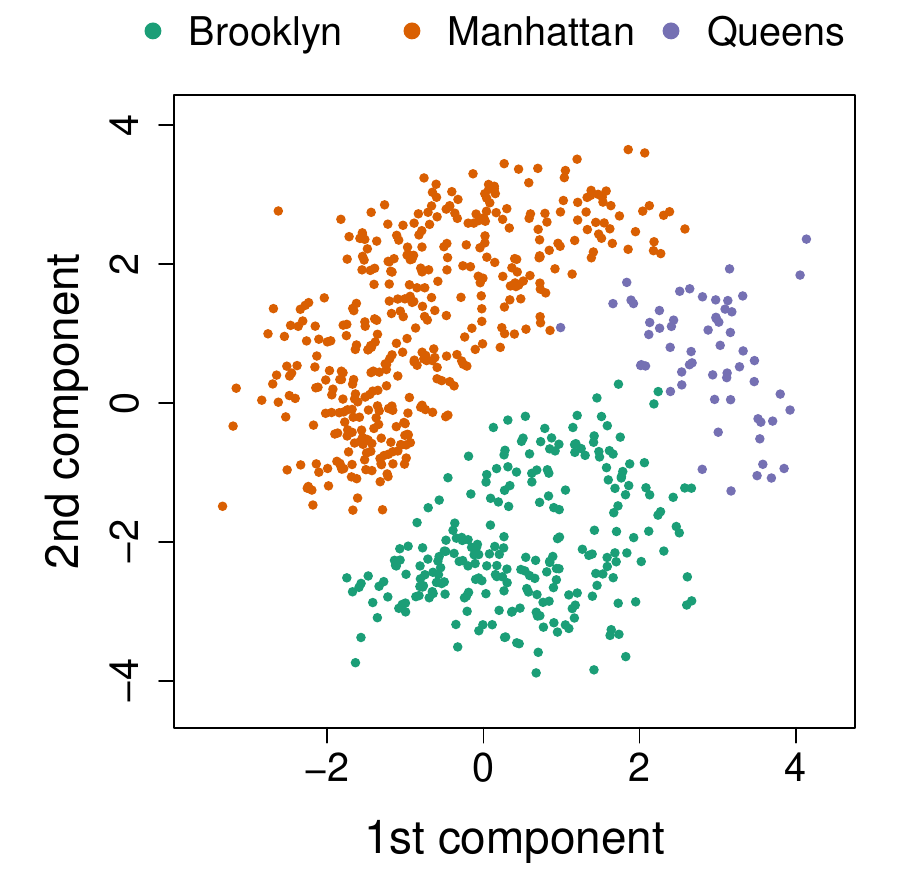}
  \caption{Morning latent vectors}  \label{fig:morning_latent}
\end{subfigure}
\begin{subfigure}{0.32\textwidth}
		\centering
		\includegraphics[width=1\linewidth]{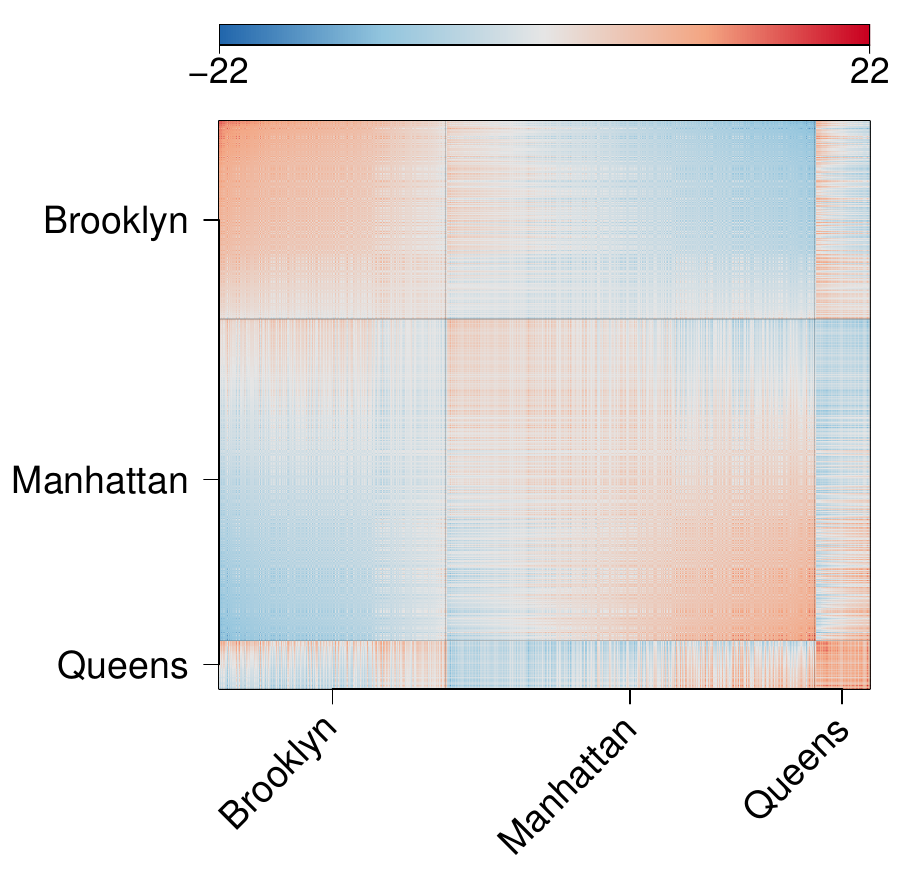}
  \caption{Morning latent similarities}\label{fig:heatmap1}
  \end{subfigure}
\begin{subfigure}{0.32\textwidth}
		\centering
		\includegraphics[width=1\linewidth]{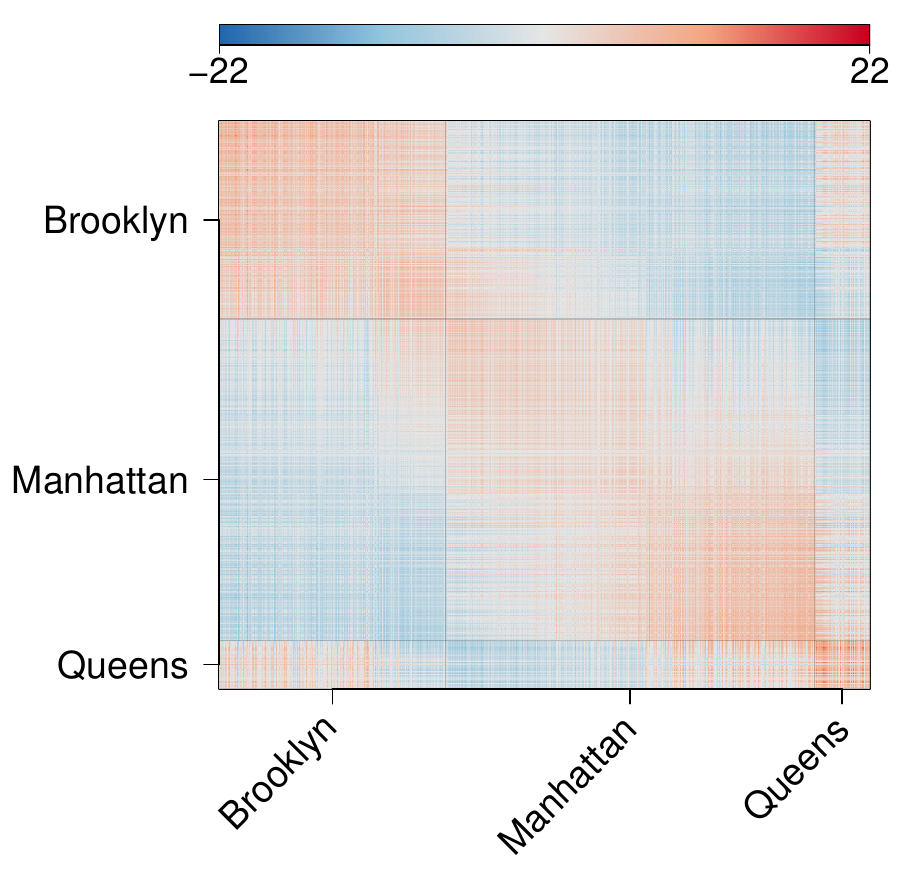}
  \caption{Evening latent similarities}\label{fig:heatmap2}
  \end{subfigure}
% \begin{subfigure}{0.32\textwidth}
% 		\centering
% 		\includegraphics[width=1\linewidth]{data/nybike_innerprod_difference.pdf}
%   \caption{Differences}
%   \end{subfigure}
    \caption{Estimation results.  Panel (a) visualizes the estimated latent positions  for the morning network.  Each point is colored based on which borough the corresponding bike station is located. Panels (b) and (c) show  two heatmaps over estimated latent similarity matrices, whose $(i,j)$-th entry is the inner product of the latent vectors for stations $i$ and $j$.} 
    \label{fig:heatmap}
\end{figure}

For each hour, we construct a weighted network,  where the network nodes represent stations, and the edge weight between two stations is the number of rides between them during that hour. 
We  fit the proposed model separately to the two hourly networks using two-dimensional latent vectors under a Poisson model. This model choice is natural because the edge weights are count-valued.  Figure~\ref{fig:heatmap} summarizes the estimation results. 
% We apply the processed estimation method with two-dimensional latent vectors to each hour's weighted networks under Poisson distribution.
% The edge weights are count-valued and therefore it is natural to model through Poisson distribution. 

% They further demonstrate the  pronounced block patterns that distinguish within-borough from between-borough travel.  
% This echoes the borough-level clusters shown in the panel (a). 

% Panel (c) further the  reveal pronounced block patterns that distinguish within-borough from between-borough travel. 

 % 2. We visualize estimated canonical parameters $\hat{\Theta}_{1,ij}$ and $\hat{\Theta}_{2,ij}$ as heatmaps. 

Figure~\ref{fig:heatmap}\subref{fig:morning_latent} visualizes the estimated latent positions for the morning network. 
It clearly shows three clusters that align well with the geographic  borough structure of the stations. 
This alignment suggests that the fitted model captures meaningful and interpretable latent structure from the observed network.
  % We also provide a map of true station locations with borough names in Section~\ref{sec:supp_data} of the Supplementary Material for comparison.  
% that in the morning peak, the estimated latent positions reveal true station boroughs shown in Figure~\ref{fig:true_bike} in the Supplementary Material. 
% They are consistent with true station locations, which are visualized in Figure ... 
The evening network exhibits a similar pattern. Due to space limitation, the detailed results are  deferred to the Supplementary Material,  along with a reference map of true station locations for comparison. 
% Similar patterns are observed for the evening network 
% This demonstrates that fitted models give reasonable and interpretable model parameters. The estimated latent vectors for evening data are similar and are presented in Figure~\ref{fig:true_latent_positions} of the Supplementary Material. 

Figures~\ref{fig:heatmap}\subref{fig:heatmap1} and \ref{fig:heatmap}\subref{fig:heatmap2} 
present heatmaps of inner products between estimated latent vectors for  the two networks. Specifically, the $(i, j)$-th entry of each heatmap corresponds to $\langle \hat{z}_{1,i},\hat{z}_{1,j}\rangle$ or $\langle \hat{z}_{2,i},\hat{z}_{2,j}\rangle$, where   $\hat{z}_{1,i}$ and $\hat{z}_{2,i}$ denote  the estimated latent vectors for station $i$ in the morning and evening  networks, respectively. 
 These inner products measure the latent affinity between stations after adjusting for node-specific baseline effects and  are informative for understanding fundamental  travel patterns. 
 Both heatmaps reveal noticeable block structure that  is consistent with the three-borough clustering pattern in  Figure~\ref{fig:heatmap}\subref{fig:morning_latent}. 

\begin{figure}[!htbp]
    \centering
\begin{subfigure}{0.32\textwidth}
		\centering
		\includegraphics[width=1\linewidth]{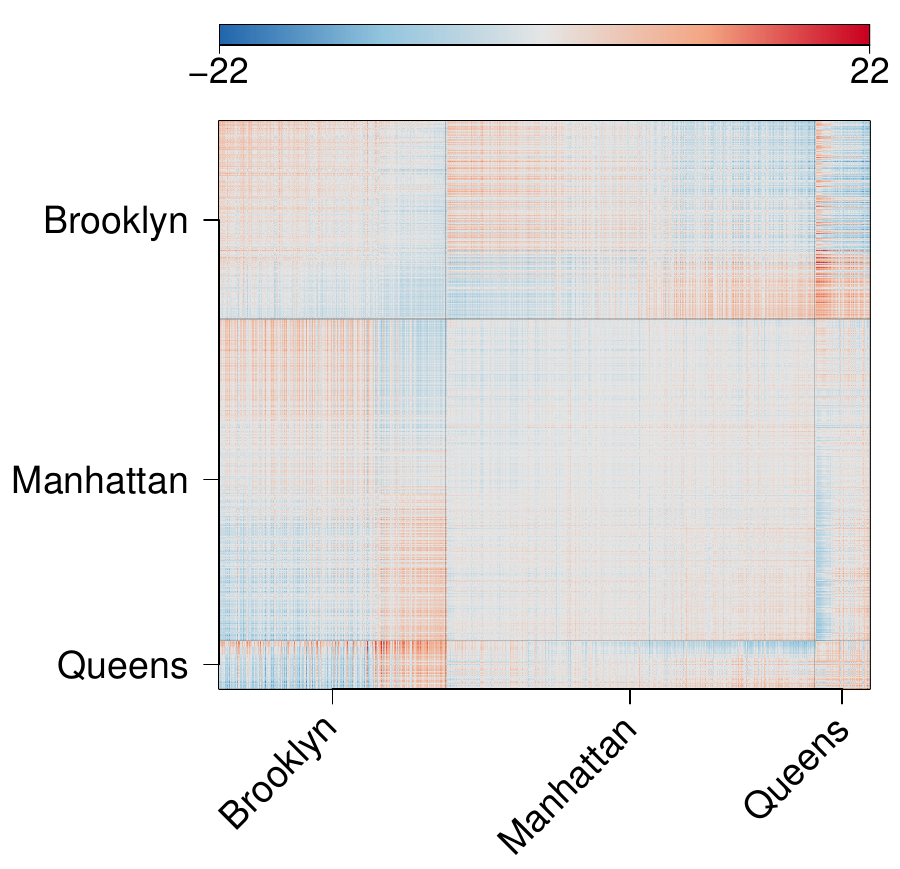}
  \caption{Estimated  differences}\label{fig:rawdiff}
  \end{subfigure}
 \begin{subfigure}{0.325\textwidth}
		\centering
		\includegraphics[width=1\linewidth]{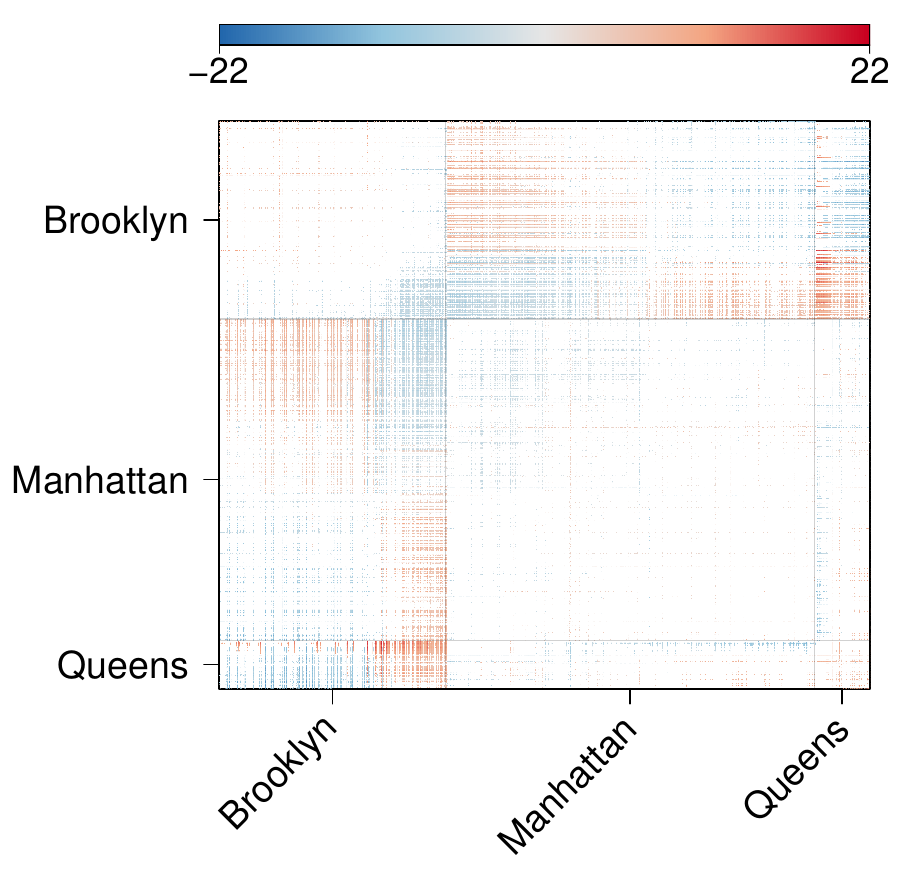}
  \caption{Significant differences}\label{fig:sigdiff}
\end{subfigure} 
  \begin{subfigure}{0.325\textwidth}
		\centering
		\includegraphics[width=1\linewidth]{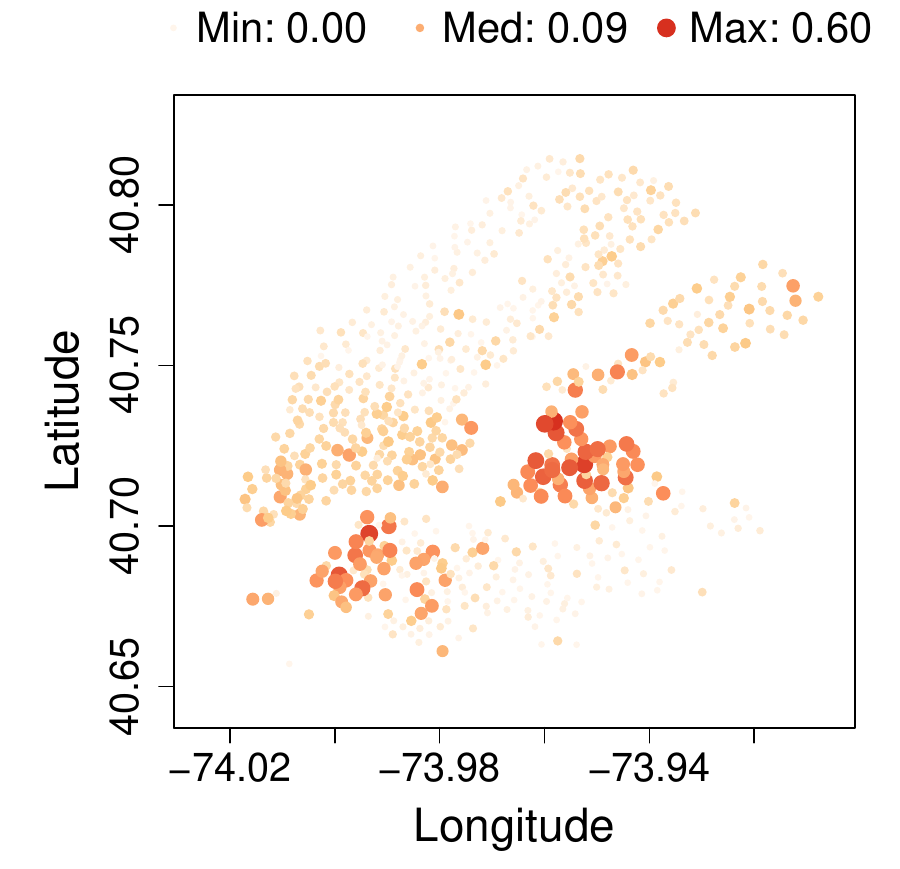}
  \caption{Station-wise rejection rates}\label{fig:rej_rates}
\end{subfigure}
    \caption{Two-sample comparison.  Panel (a) displays the difference  of  the latent similarity matrices in the morning and evening hours. Panel (b) displays the same matrix after masking the entries that are insignificant after the Benjamini–Hochberg correction in white.  Panel (c) shows the geographic locations of the stations, where the size and color of station $i$ encode the proportion of the rejected two-sample tests involving station $i$.} 
    \label{fig:net_compare}
\end{figure}

In addition to estimating each network individually, comparing the morning and evening networks can further provide  insights into how urban mobility patterns change across commuting periods. As a first step, we examine the raw difference 
% Investigating different connectivity patterns in the morning and evening transit networks could give insights into human mobility and city transits. 
% We first take direct differences 
% $\langle \hat{z}_{1,i},\hat{z}_{1,j}\rangle-\langle \hat{z}_{2,i},\hat{z}_{2,j}\rangle$ 
between the two latent similarity matrices,  presented as the heatmap  in Figure \ref{fig:net_compare}\subref{fig:rawdiff}. 
% The raw difference matrix shows diffuse changes across the whole network. There still exists certain block structure, but it could be hard to interpret. 
This raw difference matrix shows more pronounced patterns in the between-borough block compared to the single-network heatmaps in Figures~\ref{fig:heatmap}. But the changes span across all station pairs, hindering clear interpretation. 
% Although this raw difference matrix still exhibits some block structure, the overall pattern is diffuse and difficult to interpret directly.

To identify significant differences, we perform pairwise two-sample tests for differences in latent similarities between the two periods. 
% construct two-sample tests  to detect differences in latent similarities between the two hours.  
 Specifically, for each pair of stations $(i,j)$, we test whether $\langle {z}_{1,i},{z}_{1,j}\rangle$ and $\langle {z}_{2,i},{z}_{2,j}\rangle$ differ based on the estimator $\langle \hat{z}_{1,i},\hat{z}_{1,j}\rangle-\langle \hat{z}_{2,i},\hat{z}_{2,j}\rangle$.
 % , where $\hat\Theta_{1,ij}=\hat z_{1,i}^\top \hat z_{1,j}$ and $\hat\Theta_{2,ij}=\hat z_{2,i}^\top \hat z_{2,j}$. 
 Each test is constructed from the asymptotic   results in Theorem~\ref{thm:homoMLE} and Corollary~\ref{cor:gyqIhat_AN}.  
Given the $n(n-1)/2$ pairwise comparisons, we apply the Benjamini–Hochberg procedure at the level 0.05 to filter  insignificant pairs. 
% control the false discovery rate at level 0.05 using the Benjamini–Hochberg procedure. 
% Figures~\ref{fig:true_bike}~\subref{fig:sigdiff} and \ref{fig:true_bike}~\subref{fig:rej_rates} summarize the testing results. 

Figure~\ref{fig:net_compare}\subref{fig:sigdiff} displays the matrix of significant pairwise differences, with insignificant entries masked out by white color. 
In contrast to the raw difference heatmap, the statistically significant changes are concentrated primarily in off-diagonal blocks corresponding to between-borough station pairs.
This observation indicates that morning–evening travel  differences are driven more by inter-borough connectivity than by within-borough rides. 
 
 Figure \ref{fig:net_compare}\subref{fig:rej_rates}  provides a station-level summary. For each station, the station-wise rejection rate is defined as the proportion of pairwise tests involving that station that are rejected. Larger values  therefore identify stations whose latent connectivity profiles differ significantly, indicating broader reconfiguration of their connectivity patterns. When mapped geographically, the stations with high rejection rates, encoded by larger and darker points, cluster near the boundary between Brooklyn and Manhattan, separated by the East River. This spatial concentration suggests that the transition from morning to evening travel is especially pronounced for stations near waterfront areas or major inter-borough commuting corridors.
 % in the areas that near the boundary between Brooklyn and Manhattan, where the middle blank space represents the East River. 
 % This suggests that the transition from morning to evening travel may be related to waterfront activities or nearby commuting corridors. 
 % concentrated along that particular commuting corridors or waterfront  regions, rather than reflecting a uniform network-wide change. 

In summary, this example demonstrates the importance of statistically principled inference in network analysis. The raw estimated values may not clearly distinguish genuine structural change from ordinary estimation noise. By contrast, the inferential analysis produces   much  more interpretable patterns and a clearer picture of structural change. 
The developments in this paper could help identify credible structures and translate high-dimensional network estimates into scientifically meaningful conclusions on network connectivity patterns.

\section{Discussion}

In this work, we establish a unified framework that bridges the maximum likelihood estimator theory and practical algorithms under the latent space network models.  
First, we develop new theoretical results for the constrained maximum likelihood estimator that overcome the restricted eigen-gap assumption in prior analysis.  Second, we address the impractical dependence on unknown ground truth in existing algorithms by introducing  adaptive procedures. Specifically, for the projected gradient descent, we construct new  line search conditions that can adaptively select the  learning rate and prove that the explicit projection onto the unknown constraint set is unnecessary. 
For the universal singular value thresholding, we develop  an adaptive interval for elementwise projection.  
Third, we prove convergence of the adaptive algorithms to the ideal  constrained maximum likelihood estimator.

There are several natural generalizations based on the current results. 
First, the asymptotic theory in Theorem~\ref{thm:homoMLE} focuses on a fixed set of nodes. % or edge sets.  
In networks with $n$ total nodes and $n(n-1)/2$ edges, it is of interest to establish theoretical guarantees for simultaneous inference, such as false discovery rate control when comparing two networks, as demonstrated in the data analysis. 
Notably, our current asymptotic theory enables separate inference for baseline degrees and latent vectors, which  
  form two sets of high-dimensional parameters with different interpretations.
  Developing unified simultaneous inference across these parameters remains an important question. 
Second, the framework can be extended to other model formulations. For example, when additional edgewise or nodewise covariates are available, they can be added into  the network model \cite{ma2020universal} or handled via joint modeling \citep{huang2018pairwise,li2025high}. More general kernels beyond Euclidean inner products may also be considered \cite{rubin2022statistical}. The current arguments and  asymptotic results can be generalized accordingly by updating the likelihood formulations. 
Third, the estimated latent vectors can be used in downstream tasks, such as regression \citep{lunde2023conformal} or causal inference  \citep{hayes2022estimating}. The derived uncertainty and algorithmic approximation errors need to be properly accounted for when plugging in estimated parameters.

% direction is to consider Third, sparse model? 

% $U\Lambda U^{\top}$ framework  \citep{rubin2022statistical}

%%%%%%%%%%%%%%%%%%%%%%%%%%%%%%%%%%%%%%%%%%%%%%
%% Example with single Appendix:            %%
%%%%%%%%%%%%%%%%%%%%%%%%%%%%%%%%%%%%%%%%%%%%%%
% \begin{appendix}
% \section*{Title}\label{appn} %% if no title is needed, leave empty \section*{}.
% Appendices should be provided in \verb|{appendix}| environment,
% before Acknowledgements.

% If there is only one appendix,
% then please refer to it in text as \ldots\ in the \hyperref[appn]{Appendix}.
% \end{appendix}
% %%%%%%%%%%%%%%%%%%%%%%%%%%%%%%%%%%%%%%%%%%%%%%
% %% Example with multiple Appendixes:        %%
% %%%%%%%%%%%%%%%%%%%%%%%%%%%%%%%%%%%%%%%%%%%%%%
% \begin{appendix}
% \section{Title of the first appendix}\label{appA}
% If there are more than one appendix, then please refer to it
% as \ldots\ in Appendix \ref{appA}, Appendix \ref{appB}, etc.

% \section{Title of the second appendix}\label{appB}
% \subsection{First subsection of Appendix \protect\ref{appB}}

% Use the standard \LaTeX\ commands for headings in \verb|{appendix}|.
% Headings and other objects will be numbered automatically.
% \begin{equation}
% \mathcal{P}=(j_{k,1},j_{k,2},\dots,j_{k,m(k)}). \label{path}
% \end{equation}

% Sample of cross-reference to the formula (\ref{path}) in Appendix \ref{appB}.
% \end{appendix}

%%%%%%%%%%%%%%%%%%%%%%%%%%%%%%%%%%%%%%%%%%%%%%
%% Support information, if any,             %%
%% should be provided in the                %%
%% Acknowledgements section.                %%
%%%%%%%%%%%%%%%%%%%%%%%%%%%%%%%%%%%%%%%%%%%%%%
\begin{acks}[Acknowledgments]
% The authors would like to thank the anonymous referees, an Associate
% Editor and the Editor for their constructive comments that improved the
% quality of this paper.
% YH was partially supported by NSF Grant DMS-2515523 and Wisconsin Alumni Research Foundation. YH is the corresponding author. 
\end{acks}

%%%%%%%%%%%%%%%%%%%%%%%%%%%%%%%%%%%%%%%%%%%%%%
%% Funding information, if any,             %%
%% should be provided in the                %%
%% funding section.                         %%
%%%%%%%%%%%%%%%%%%%%%%%%%%%%%%%%%%%%%%%%%%%%%%
% \begin{funding}
% The first author was supported by NSF Grant DMS-??-??????.

% The second author was supported in part by NIH Grant ???????????.
% \end{funding}

%%%%%%%%%%%%%%%%%%%%%%%%%%%%%%%%%%%%%%%%%%%%%%
%% Supplementary Material, including data   %%
%% sets and code, should be provided in     %%
%% {supplement} environment with title      %%
%% and short description. It cannot be      %%
%% available exclusively as external link.  %%
%% All Supplementary Material must be       %%
%% available to the reader on Project       %%
%% Euclid with the published article.       %%
%%%%%%%%%%%%%%%%%%%%%%%%%%%%%%%%%%%%%%%%%%%%%%
% \vspace{-1.5em}
\begin{supplement}
\stitle{Supplement to ``Bridging Theory and Practice: Statistical Inference of Latent Space Models for Networks''.} 
\sdescription{Additional numerical results and proofs are deferred to the Supplementary Material.}
\end{supplement}

\end{document}